\documentclass[twocolumn,sloppy]{article}
\usepackage{graphicx}

\usepackage{latexsym}
\usepackage{amsfonts}
\usepackage{amsthm}
\usepackage{amsmath}
\def\R{{\cal R}}
\def\E{{\cal E}}

\newtheorem{theorem}{Theorem}[section]
\newtheorem{lemma}[theorem]{Lemma}
\newtheorem{conjecture}[theorem]{Conjecture}
\newtheorem{corollary}[theorem]{Corollary}

\newtheorem{proposition}[theorem]{Proposition}

\theoremstyle{definition}

\newtheorem{remark}[theorem]{Remark}

\font\twrm=cmr8

\def\oldlabel#1{\relax}

\def\doct{\delta_{oct}}
\def\dtet{\delta_{tet}}
\def\pt{\hbox{\it pt}}
\def\Vol{\hbox{vol}}
\def\sol{\operatorname{sol}}
\def\quo{\operatorname{quo}}

\def\vor{\operatorname{vor}}
\def\sign{\operatorname{sign}}
\def\octavor{\operatorname{octavor}}
\def\dih{\operatorname{dih}}
\def\mn{{\operatorname{min}}}
\def\mx{{\operatorname{max}}}
\def\Adih{\operatorname{Adih}}
\def\arc{\operatorname{arc}}

\def\sc{{\operatorname{sc}}}
\def\tausc{{\tau\!\operatorname{sc}}}
\def\piF{{\pi_F}}
\def\xiG{\xi_\Gamma}
\def\piG{\pi_\Gamma}
\def\xiV{\xi_V}
\def\piV{\pi_V}
\def\tauLP{{\tau_{\hbox{\twrm LP}}}}
\def\tlp{\tau_{\hbox{\twrm LP}}}  
\def\sLP{\sigma_{\hbox{\twrm LP}}}  
\def\A{{\mathbf A}}
\def\squander{(4\pi\zeta-8)\,\pt}

\def\refno#1{\hbox{}\nobreak\hfill {\tt (#1)}}
\def\mark#1{\hbox{\tt #1}}
\def\Sfour{{{\cal\mathbf S}_4^+}}
\def\Sminus{{{\cal\mathbf S}_3^-}}
\def\Splus{{{\cal\mathbf S}_3^+}}

\def\R{{\mathbb R}}
\def\N{{\mathbb N}}

\def\tildeF{{\hbox{$\tilde F$}}}
\def\zloop{{z}_{loop}}
\def\dloop{{\delta}_{loop}}

\begin{document}

\title{The Kepler Conjecture}
\author{Thomas C. Hales}

\maketitle

\begin{abstract}
    We present the final part of the proof of the Kepler
    Conjecture.
\end{abstract}

\def\today{\ifcase\month\or
    January\or February\or March\or April\or May\or June\or
    July\or August\or September\or October\or November\or December\fi
    \space\number\day, \number\year}

\section{Overview}
\label{sec:overview}

This section describes the structure of the proof of the Kepler
Conjecture.

\begin{theorem}
\label{theorem:kepler}
(The Kepler Conjecture)  No packing of
congruent balls in Euclidean three space has density greater than
that of the face-centered cubic packing.
\end{theorem}

This density is $\pi/\sqrt{18}\approx 0.74.$

\begin{figure}[htb]
  \centering
  \includegraphics{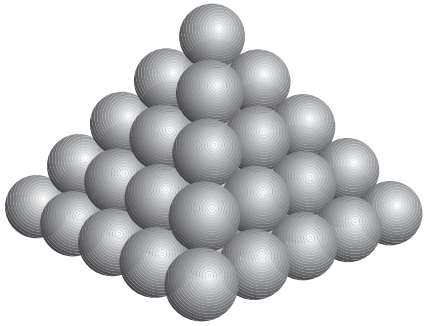}
  \caption{The face-centered cubic packing}
  \label{fig:fcc}
\end{figure}

The proof of this result is scattered throughout several papers.
(Every article in the bibliography will be required.) Here, we
describe the top-level outline of the proof and give references to
the sources of the details of the proof.

By a {\it packing}, we mean an arrangement of balls that are
non-overlapping in the sense that the interiors of the balls are
pairwise disjoint. Consider a packing of congruent balls in
Euclidean three space. There is no harm in assuming that all the
balls have unit radius. The density of a packing does not decrease
when balls are added to the packing. Thus, to answer a question
about the greatest possible density we may add non-overlapping
balls until there is no room to add further balls. Such a packing
will be said to be {\it saturated}.

Let $\Lambda$ be the set of centers of the balls in a saturated
packing.  Our choice of radius for the balls implies that any two
points in $\Lambda$ have distance at least $2$ from each other. We
call the points of $\Lambda$ {\it vertices}.  Let $B(x,r)$ denote
the ball in Euclidean three space at center $x$ and radius $r$.
Let $\delta(x,r,\Lambda)$ be the finite density, defined by the
ratio of $A(x,r,\Lambda)$ to the volume of $B(x,r)$, where
$A(x,r,\Lambda)$ is defined as the volume of the intersection with
$B(x,r)$ of the union of all balls in the packing. Set
$\Lambda(x,r) = \Lambda \cap B(x,r)$.

Recall that the Voronoi cell $\Omega(v)$ around a vertex $v\in
\Lambda$ is the set of points closer to $v$ than to any other ball
center. The volume of each Voronoi cell in the face-centered cubic
packing is $\sqrt{32}$.  This is also the volume of each Voronoi
cell in the hexagonal-close packing.

Let $a:\Lambda\to\R$ be a function.  We say that $a$ is {\it
negligible\/} if there is a constant $C_1$ such that  for all
$r\ge1$, we have
$$\sum_{v\in\Lambda(x,r)} a(v) \le C_1 r^2.$$
We say that the function $a$ is {\it fcc-compatible\/} if for all
$v\in\Lambda$ we have the inequality
$$\sqrt{32}\le \Vol(\Omega(v)) + a(v).$$

\begin{lemma}
\label{lemma:deltabound} If there exists a negligible
fcc-compatible function $a:\Lambda\to\R$ for a saturated packing
$\Lambda$, then there exists a constant $C$ such that for all
$r\ge1$, we have
$$\delta(x,r,\Lambda)
\le \pi/\sqrt{18} + C/r.$$
\end{lemma}

\begin{proof}
The numerator $A(x,r,\Lambda)$ of $\delta(x,r,\Lambda)$ is at most
the product of the volume of a ball $4\pi/3$ with the number
$|\Lambda(x,r+1)|$ of balls intersecting $B(x,r)$.  Hence
    \begin{equation}
    A(x,r,\Lambda) \le |\Lambda(x,r+1)| 4\pi/3.
    \label{eqn:Abound}
    \end{equation}

In a saturated packing each Voronoi cell is contained in a ball of
radius $2$ centered at the {\it center} of the cell.  The volume
of the ball $B(x,r+3)$ is at least the combined volume of Voronoi
cells lying entirely in the ball. This observation, combined with
fcc-compatibility and negligibility, gives
    \begin{equation}
    \begin{split}
    \sqrt{32}|\Lambda(x,r+1)|
    &\le \sum_{v\in\Lambda(x,r+1)} (a(v) +
    \Vol(\Omega(v))) \\
    &\le C_1 (r+1)^2 + \Vol\,B(x,r+3) \\
    &\le C_1 (r+1)^2 + (1+3/r)^3 \Vol\,B(x,r)
    \label{eqn:Bbound}
    \end{split}.
    \end{equation}
Divide through by $\Vol\,B(x,r)$ and eliminate $|\Lambda(x,r+1)|$
between Inequality (\ref{eqn:Abound}) and Inequality
(\ref{eqn:Bbound}) to get
    $$\delta(x,r,\Lambda)
        \le \frac{\pi}{\sqrt{18}} (1+3/r)^3 + C_1 \frac{(r+1)^2}{r^3\sqrt{32}}.
    $$
The result follows for an appropriately chosen constant $C$.
\end{proof}

\begin{remark} We take the precise meaning of the Kepler Conjecture to
be a bound on the essential supremum of the function $\delta(x,r)$
as $r$ tends to infinity. Lemma \ref{lemma:deltabound} implies
that the essential supremum of $\delta(x,r,\Lambda)$ is bounded
above by $\pi/\sqrt{18}$, provided a negligible fcc-compatible
function can be found.  The strategy will be to define a
negligible function, and then to solve an optimization problem in
finitely many variables to establish that it is fcc-compatible.
\end{remark}

The paper \cite{formulation} defines a compact topological space
$X$ and a continuous function $\sigma$ on that space.

The topological space $X$ is directly related to packings. If
$\Lambda$ is a saturated packing, then there is a geometric object
$D(v,\Lambda)$ constructed around each vertex $v\in\Lambda$.
$D(v,\Lambda)$ depends on $\Lambda$ only through the vertices in
$\Lambda$ at distance at most $4$ from $v$.  The objects
$D(v,\Lambda)$ are called {\it decomposition stars}, and the space
of all decomposition stars is precisely $X$.

The following constants arise in that paper.

Let $\dtet$ be the packing density of a regular tetrahedron.  That
is, let $S$ be a regular tetrahedron of edge length $2$.  Let $B$
the part of $S$ that lies within distance $1$ of some vertex. Then
$\dtet$ is the ratio of the volume of $B$ to the volume of $S$. We
have $\dtet = \sqrt{8} \arctan(\sqrt{2}/5)$.

Let $\doct$ be the packing density of a regular octahedron of edge
length $2$, again constructed as the ratio of the volume of points
within distance $1$ of a vertex to the volume of the octahedron.

The density of the face-centered cubic packing is a weighted
average of these two ratios
    $$\frac{\pi}{\sqrt{18}} = \frac{\dtet}{3} + \frac{2 \doct}{3}.$$
This determines the exact value of $\doct$ in terms of $\dtet$. We
have $\doct \approx 0.72$.

Let $\pt = -\pi/3 + \sqrt{2}\dtet\approx 0.05537$

The following conjecture is made in \cite{formulation}

\begin{conjecture}
\label{conjecture:sigma}
 The maximum of $\sigma$ on $X$ is the constant
$8\,\pt\approx 0.442989$.
\end{conjecture}

\begin{lemma}
\label{lemma:exista} An affirmative answer to Conjecture
\ref{conjecture:sigma} implies the existence of a negligible
fcc-compatible function for every saturated packing $\Lambda$.
\end{lemma}

\begin{proof}
For any saturated packing $\Lambda$ define a function
$a:\Lambda\to\R$ by
$$-\sigma(D(v,\Lambda))/(4\doct) + 4\pi/(3\doct) = \Vol(\Omega(v)) + a(v).$$
Negligibility follows from \cite[Prop. 3.14 (proof)]{formulation}.
The upper bound of $8\,\pt$ gives a lower bound $$-8\,\pt/(4\doct)
+ 4\pi/(3\doct) \le \Vol(\Omega(v)) + a(v).$$ The left-hand side
of this inequality evaluates to $\sqrt{32}$, and this establishes
fcc-compatibility.
\end{proof}

\begin{theorem}
\label{theorem:sigma} Conjecture \ref{conjecture:sigma} is true.
That is,
  the maximum of the function $\sigma$ on the
topological space $X$ of all decomposition stars is $8\,\pt$.
\end{theorem}

Theorem \ref{theorem:sigma}, Lemma \ref{lemma:exista}, and Lemma
\ref{lemma:deltabound} combine to give a proof of the Kepler
Conjecture \ref{theorem:kepler}

Let $t_0=1.255$ ($2t_0 = 2.51$).  This is a parameter that is used
for truncation throughout the series of papers on the Kepler
Conjecture.  The significance of this particular choice of
truncation parameter will be explained elsewhere.

Let $U(v,\Lambda)$ be the set of vertices in $\Lambda$ at distance
at most $2t_0$ from $v$.  From a decomposition star $D(v,\Lambda)$
it is possible to recover $U(v,\Lambda)$ (at least up to Euclidean
translation: $U\mapsto U+y$, for $y\in\R^3$).  We can completely
characterize the decomposition stars at which the maximum of
$\sigma$ is attained.

\begin{theorem}
\label{theorem:sharp} Let $D$ be a decomposition star at which
the maximum $8\,\pt$ is attained.  Then the set $U(D)$ of vectors
at distance at most $2t_0$ from the center has cardinality $12$.
Up to Euclidean motion, $U(D)$ is the kissing arrangement of the
$12$ balls around a central ball in the face-centered cubic
packing or hexagonal-close packing.
\end{theorem}

\subsection{Outline of proofs}

To prove Theorems \ref{theorem:sigma} and \ref{theorem:sharp}, we
wish to show that there is no counterexample.  That is, we wish to
show that there is no decomposition star $D$ with value $\sigma(D)
> 8\,\pt$.  We reason by contradiction, assuming the existence of
such a decomposition star.  With this in mind, we call $D$ a {\it
contravening decomposition star}, if
    $$\sigma(D)\ge 8\,\pt.$$
In much of what follows we will tacitly assume that every
decomposition star under discussion is a contravening one.  Thus,
when we say that no decomposition stars exist with a given
property, it should be interpreted as saying that no such
contravening decomposition stars exist.

To each contravening decomposition star, we associate a
(combinatorial) plane graph.  A restrictive list of properties of
plane graphs is described in Section \ref{definition:tame}.  Any
plane graph satisfying these properties is said to be {\it tame}.
All tame plane graphs have been classified. (There are several
thousand, up to isomorphism.) Theorem
\ref{theorem:classification}, asserts that the plane graph
attached to each contravening decomposition star is tame. By the
classification of such graphs, this reduces the proof of the
Kepler Conjecture to the analysis of the decomposition stars
attached to the finite explicit list of tame plane graphs.

A few of the tame plane graphs are of particular interest. Every
decomposition star attached to the face-centered cubic packing
gives the same plane graph (up to isomorphism).  Call it
$G_{fcc}$.  Likewise, every decomposition star attached to the
hexagonal-close packing gives the same plane graph $G_{hcp}$. Let
$X_{crit}$ be the set of decomposition stars $D$ such that the set
$U(D)$ of vertices is the kissing arrangement of the $12$ balls
around a central ball in the face-centered cubic or
hexagonal-close packing.  There are only finitely many orbits of
$X_{crit}$ under the group of Euclidean motions.
\begin{figure}[htb]
  \centering
  \includegraphics{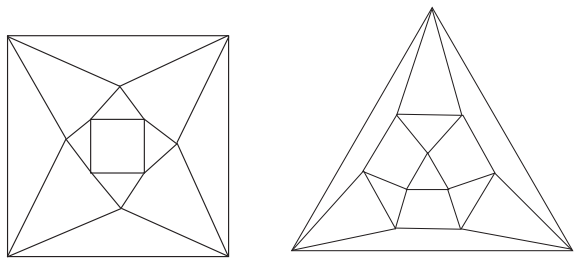}
  \caption{The plane graphs $G_{fcc}$ and $G_{hcp}$}
  \label{fig:gfcchcp}
\end{figure}

In \cite[Lemma 3.13]{formulation}, the necessary local analysis is
carried out to prove the following local optimality.

\begin{lemma}
A decomposition star whose plane graph is $G_{fcc}$ or $G_{hcp}$
has score at most $8\,\pt$, with equality precisely when the
decomposition star belongs to $X_{crit}$. \end{lemma}

In light of this result, we prove \ref{theorem:sigma} and
\ref{theorem:sharp} by proving that any decomposition star whose
graph is tame and not equal to $G_{fcc}$ or $G_{hcp}$ is not
contravening

\begin{figure}[htb]
  \centering
  \includegraphics{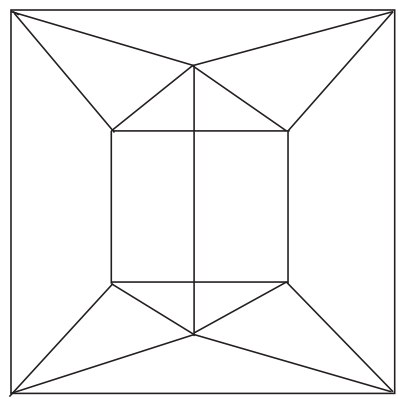}
  \caption{The plane graph $G_{pent}$}
  \label{fig:pentahedral}
\end{figure}

There is one more tame plane graph that is particularly
troublesome.  It is the graph $G_{pent}$ obtained from the
pictured configuration of twelve balls tangent to a given central
ball (Figure \ref{fig:pentahedral}). (Place a ball at the north
pole, another at the south pole, and then form two pentagonal
rings of five balls.) This case requires individualized attention.
S. Ferguson proves in \cite{thesis} that if $D$ is any
decomposition star with this graph, then $\sigma(D)< 8\,\pt$.

To eliminate the remaining cases, more-or-less generic arguments
can be used.  A linear program is attached to each tame graph $G$.
The linear program can be viewed as a linear relaxation of the
nonlinear optimization problem of maximizing $\sigma$ over all
decomposition stars with a given tame graph $G$. Because it is
obtained by relaxing the constraints on the nonlinear problem, the
maximum of the linear problem is an upper bound on the maximum of
the original nonlinear problem.  Whenever the linear programming
maximum is less than $8\,\pt$, it can be concluded that there is
no contravening decomposition star with the given tame graph $G$.
This linear programming approach eliminates most tame graphs.

When a single linear program fails to give the desired bound, it
is broken into a series of linear programming bounds, by branch
and bound techniques.  For every tame plane graph $G$ other than
$G_{hcp}$, $G_{fcc}$, and $G_{pent}$ , we produce a series of
linear programs that establish that there is no contravening
decomposition star with graph $G$.  When every face of the plane
graph is a triangle or quadrilateral, this is accomplished in
\cite{part3}.  The general case is completed in the final sections
of this paper.

\subsection{Organization of the Paper}

This paper has been written in such a way that several sections
can be skipped without disrupting the expository flow.  To get a
general overview of the paper, Sections \ref{sec:overview},
\ref{sec:tame}, \ref{sec:startame}, and \ref{sec:linearprogram}
may be read. To get a somewhat more detailed view of the results
and methods of proof,  Sections \ref{sec:overview},
\ref{sec:tame}, \ref{sec:classification}, \ref{sec:startame},
\ref{sec:contraproof}, \ref{sec:linearprogram},
\ref{sec:lpformulation}, and the first two appendices may be read.

The paper has been written so as to be substantially independent
of the other papers in the series.  An appendix gives a summary of
the results of \cite{part4} in a form that is suited to the
purposes of this paper.

A historical account of the Kepler Conjecture can be found at
\cite{historical}.  An introduction to the ideas of the proof can
be found in \cite{honeycombs}.  An introduction to the algorithms
can be found at \cite{algorithm}.  Speculation on a
second-generation design of a proof can be found in
\cite{algorithm} and \cite{arbeitstagung}.

\section{Tame Plane Graphs}
\label{sec:tame}

We give a succinct statement of properties that hold of all plane
graphs considered in the proof of the Kepler conjecture.  Such
graphs are said to be tame. A list of all tame graphs has been
generated by computer.

\subsection{Basic Definitions}

An {\it n-cycle\/} is a finite set $C$ of cardinality $n$,together
with a cyclic permutation $s$ of $C$.   We write $s$ in the form
$v\mapsto s(v,C)$, for $v\in C$.  The element $s(v,C)$ is called
the {\it successor\/} of $v$ (in $C$).  A {\it cycle\/} is a
$n$-cycle for some natural number $n$. By abuse of language, we
often identify $C$ with the cycle. The natural number $n$ is the
{\it length\/} of the cycle.

(We may assume that all vertices of all graphs lie in some large
finite set of vertices $\Omega$, if we wish to arrange that only
finitely many graphs occur in this discussion.)

Let $G$ be a nonempty finite set of cycles (called faces) of
length at least $3$. The elements of faces are called the {\it
vertices} of $G$. An unordered pair of vertices $\{v,w\}$ such
that one element is the successor of the other in some face is
called an {\it edge}. The vertices $v$ and $w$ are then said to be
{\it adjacent}. The set $G$ is a {\it plane graph} if three
conditions hold.
    \begin{enumerate}
    \item If an element $v$ has successor $w$ in some face $F$, then
there is a unique face (call it $s'(F,v)$) in $G$ for which $v$ is
the successor of $w$. (Thus, $v=s(w,s'(F,v))$, and each edge
occurs twice with opposite orientation.)
    \item  For each vertex $v$, the function $F\mapsto s'(F,v)$ is
    a cyclic permutation of the set of faces containing $v$.
    \item  Euler's formula holds relating the number of vertices $V$, the
number of edges $E$, and the number of faces $F$:
    $$V-E+F = 2.$$
    \end{enumerate}
(The set of vertices and edges of a plane graph form a planar
graph in the usual graph-theoretic sense of admitting an embedding
into the plane.)

Let $len$ be the length function on faces. Let $tri(v)$ be the
number of triangles containing a vertex $v$. Faces of length $3$
are called {\it triangles}, those of length $4$ are called {\it
quadrilaterals}, and so forth.  An face of length at least $5$ is
called an {\it exceptional\/} face.

Two plane graphs are {\it properly} isomorphic if there is a
bijection of vertices inducing a bijection of faces.  For each
plane graph, there is an opposite plane graph $G^{op}$ obtained by
reversing the cyclic order of vertices in each face.  A plane
graph $G$ is  isomorphic to another if $G$ or $G^{op}$ is properly
isomorphic to the other.

The {\it degree\/} of a vertex is the number of faces it belongs
to. An {\it $n$-circuit\/} in $G$ is a cycle $C$ in the vertex-set
of $G$, such that for every $v\in C$, it forms an edge in $G$ with
its successor: that is, $(v,s(v,C))$ is an edge of $G$.  In a
plane graph $G$ we have a combinatorial form of the Jordan curve
theorem: each $n$-circuit determines a partition of $G$ into two
sets of faces.

The {\it type\/} of a vertex is defined to be a triple of
non-negative integers $(p,q,r)$, where $p$ is the number of
triangles containing the vertex, $q$ is the number of
quadrilaterals containing it, and $r$ is the number of exceptional
faces. When $r=0$, we abbreviate the type to the ordered pair
$(p,q)$.

\subsection{Weight Assignments}

We call the constant $14.8$, which arises repeatedly in this
section, the {\it target}.

  Define $a:\N\to \R$ by
  $$a = \begin{cases}
    14.8 &n=0,1,2,\\
    1.4 & n=3,\\
    1.5 & n=4,\\
    0 & \text{otherwise.}
  \end{cases}
  $$
  Define $b:\N\times \N\to \R$ by $b(p,q)=14.8$,
  except for the values in the following table
  (with the understanding that $x=14.8$):
  $$\begin{matrix}  &q=0&1&2&3&4\\
           p=0&x&x&x&7.135&10.649\\
           1&x&x&6.95&7.135&x\\
           2&x&8.5&4.756&12.981&x\\
           3&x&3.642&8.334&x&x\\
           4&4.139&3.781&x&x&x\\
           5&0.55&11.22&x&x&x\\
           6&6.339&x&x&x&x
\end{matrix}$$
  Define $c:\N\to \R$ by
  $$c = \begin{cases}
    1 & n=3,\\
    0 & n=4,\\
    -1.03 &n=5,\\
    -2.06 &n=6,\\
    -3.03 &\text{otherwise.}
    \end{cases}
    $$
    Define $d:\N\to \R$ by
  $$d = \begin{cases}
    0 & n=3, \\
    2.378 & n=4, \\
    4.896 & n=5, \\
    7.414 & n=6, \\
    9.932 & n=7, \\
    10.916 & n=8,\\
    14.8 & \text{otherwise}.
  \end{cases}
  $$

A set $V$ of vertices is called a {\it separated\/} set of
vertices if the following four conditions hold.
    \begin{enumerate}
      \item For every vertex in $V$ there is an exceptional face
         containing it.
      \item No two
        vertices in $V$ are adjacent.
      \item No two vertices
        in $V$ lie on a common quadrilateral.
      \item Each vertex in $V$ has degree 5.
    \end{enumerate}

A {\it weight assignment\/} of a plane graph $G$ is a function
$w:G\to \R$ taking values in the set of non-negative real numbers.
A weight assignment is {\it admissible} if the following
properties hold:
\begin{enumerate}
  \item If the face $F$ has length $n$, then
        $w(F) \ge d(n)$
  \item If $v$ has type $(p,q)$, then
        $$\sum_{v\in F} w(F) \ge b(p,q).$$
        \label{admissible:b}
  \item Let $V$ be any set of vertices of type $(5,0)$.
        If the cardinality of $V$ is $k\le 4$,
        then
        $$\sum_{V\cap F\ne\emptyset} w(F) \ge 0.55 k.$$
  \item Let $V$ be any separated set of vertices.
        Then
        $$\sum_{V\cap F\ne\emptyset} (w(F) -d(len(F)))
            \ge \sum_{v\in V} a(tri(v)).$$
        \label{definition:admissible:excess}
\end{enumerate}

 The sum
$\sum_F w(F)$ is called the {\it total weight} of $w$.

\subsection{Plane Graph Properties}
\label{sec:graphproperty}

We say that a plane graph is {\it tame\/} if it satisfies the
following conditions.

\begin{enumerate}
    \label{definition:tame}
    \item The length of each face is (at least $3$ and) at most $8$.
    \label{definition:tame:length}

    \item Every $3$-circuit is a face or the opposite of a face.
    \label{definition:tame:3-circuit}

    \item Every $4$-circuit surrounds one of the cases illustrated in Figure
    \ref{fig:fourcircuit}.
    \label{definition:tame:4-circuit}
    \begin{figure}[htb]
        \centering
        \includegraphics{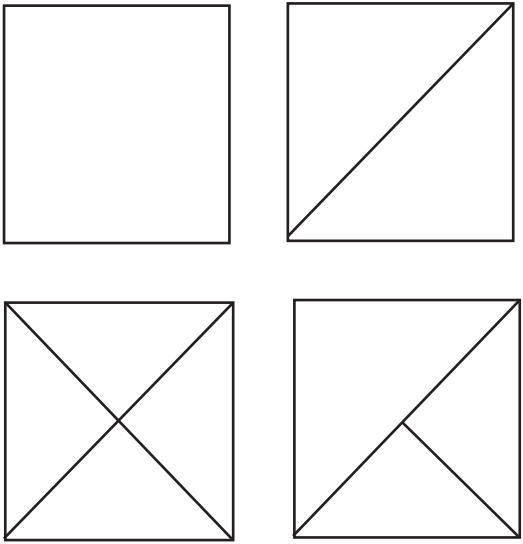}
        \caption{Tame $4$-circuits}
        \label{fig:fourcircuit}
    \end{figure}

    \item The degree of every vertex is (at least $2$ and) at most
    $6$.
    \label{definition:tame:degree}

    \item If a vertex is contained in an exceptional face,
        then the degree of the vertex is at most $5$.
    \label{definition:tame:degreeE}

    \item $$\sum_F c(len(F)) \ge 8,$$
    \label{definition:tame:score}

    \item There exists an admissible weight assignment
        of total weight less than the target, 14.8.
    \label{definition:tame:squander}

\end{enumerate}
It follows from the definitions that the abstract vertex-edge
graph of $G$ has no loops or multiple joins.  Also, by
construction, every vertex lies in at least two faces. Property
\ref{definition:tame:score} asserts that the graph has at least
        eight triangles.

\subsection{Classification of Tame Plane Graphs}

A list of several thousand plane graphs appears at \cite{web}.

\begin{theorem}
\label{theorem:classification} Every tame plane graph is
isomorphic to a plane graph in this list.
\end{theorem}

\section{Proof of Classification}
\label{sec:classification}

The results of this section are not needed except in the proof of
Theorem \ref{theorem:classification}.

Computers are used to generate a list of all tame plane graphs and
to check them against the archive of tame plane graphs.  We will
describe a finite state machine that produces all tame plane
graphs. This machine is not particularly efficient, and so we also
include a description of pruning strategies that prevent a
combinatorial explosion of possibilities.

\subsection{Basic Definitions}

A {\it partial plane graph} is a plane graph with additional data:
every face is marked as ``complete'' or ``incomplete.'' We call a
face {\it complete\/} or {\it incomplete\/} according to the
markings. We require the following condition.
\begin{itemize}
  \item {\it No two incomplete faces share an edge.}
  \label{definition:partial}
\end{itemize}

Each unmarked plane graph is identified with the marked plane
graph in which every face is complete. We represent a partial
plane graph graphically by deleting one face (the face at
infinity) and drawing the others and shading those that are
complete.

A {\it patch\/} is a partial plane graph $P$ with two
distinguished faces $F_1$ and $F_2$, such that the following hold.
\begin{itemize}
  \item Every vertex of $P$ lies in $F_1$ or $F_2$.
  \item The face $F_2$ is the only complete face.
  \item $F_1$ and $F_2$ share an edge.
  \item Every vertex of $F_2$ that is not in $F_1$
has degree $2$.
\end{itemize}

$F_1$ and $F_2$ will be referred to as the distinguished
incomplete and the distinguished complete faces, respectively.

\begin{figure}[htb]
  \centering
  \includegraphics{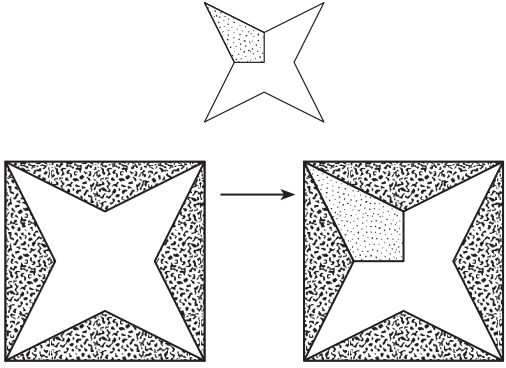}
  \caption{Patching a plane graph}
  \label{fig:patching}
\end{figure}

Patches can be used to modify a partial plane graph as follows.
Let $F$ be an incomplete face of length $n$ in a partial plane
graph $G$. Let $P$ be a patch whose incomplete distinguished face
$F_1$ has length $n$. Replace $P$ with a properly isomorphic patch
$P'$ in which the image of $F_1$ is equal to $F^{op}$ and in which
no other vertex of $P'$ is a vertex of $G$. Then
    $$ G' = \{F' \in G\cup P' : F'\ne F^{op}, F'\ne F\}$$
is a partial plane graph. Intuitively, we cut away the faces $F$
and $F_1$ from their plane graphs, and glue the holes together
along the boundary (Figure \ref{fig:patching}). (It is immediate
that the Condition \ref{definition:partial} in the definition of
partial plane graphs is maintained by this process.) There are $n$
distinct proper ways of identifying $F_1$ with $F^{op}$ in this
construction, and we let $\phi$ be this identification. The
isomorphism class of $G'$ is uniquely determined by the
isomorphism class of $G$, the isomorphism class of $P$, and $\phi$
(ranging over proper bijections $\phi:F_1\mapsto F^{op}$).

\subsection{A Finite State Machine}

For a fixed $N$ we define a finite state machine as follows. The
states of the finite state machine are isomorphism classes of
partial plane graphs $G$ with at most $N$ vertices. The
transitions from one state $G$ to another are isomorphism classes
of pairs $(P,\phi)$ where $P$ is a patch, and $\phi$ pairs an
incomplete face of $G$ with the distinguished incomplete face of
$P$. However, we exclude a transition $(P,\phi)$ at a state if the
resulting partial plane graphs contains more than $N$ vertices.
Figure \ref{fig:patching} shows two states and a transition
between them.

The initial states $I_n$ of the finite state machine are defined
to be the isomorphism classes of partial plane graphs with two
faces:
    $$\{(1,2,\ldots,n),(n,n-1,\ldots,1)\}$$
where $n\le N$, one face is complete, and the other is incomplete.
In other words, they are patches with exactly two faces.

A terminal state of this finite state machine is one in which
every face is complete.  By construction, these are (isomorphism
classes of) plane graphs with at most $N$ vertices.

\begin{lemma}
    \label{lemma:reachable}
Let $G$ be a plane graph with at most $N$ vertices. Then its state
in the machine is reachable from an initial state through a series
of transitions.
\end{lemma}

\begin{proof}  Pick an face in $G$ of length $n$ and identify it with
the complete face in the initial state $I_n$. At any stage at
state $G'$, we have an identification of all of the vertices of
the plane graph $G'$ with some of the vertices of $G$, and an
identification of all of the complete faces of $G'$ with some of
the faces of $G$ (all faces of $G$ are complete). Pick an
incomplete face $F$ of $G'$ and an oriented edge along that face.
We let $F'$ be the complete face of $G$ with that edge, with the
same orientation on that edge as $F$. Create a patch with
distinguished faces $F_1 = F^{op}$ and $F_2 = F'$. ($F_1$ and
$F_2$ determine the patch up to isomorphism.) It is immediate that
the conditions defining a patch are fulfilled. Continue in this
way until a graph isomorphic to $G$ is reached.
\end{proof}

\begin{remark}
It is an elementary matter to generate all patches $P$ such that
the distinguished faces have given lengths $n$ and $m$. Patching
is also entirely algorithmic, and thus by following all paths
through the finite state machine, we obtain all plane graphs with
at most $N$ vertices.
\end{remark}

\subsection{Pruning Strategies}

Although we reach all graphs in this manner, it is not
computationally efficient. We introduce pruning strategies to
increase the efficiency of the search. We can terminate our search
along a path through the finite state machine, if we can
determine:
\begin{enumerate}
  \item Every terminal graph
along that path violates one of the defining properties of
tameness, or
    \label{enum:not-tame}
  \item An isomorphic terminal graph will be reached by
some other path that will not be terminated early.
    \label{enum:branch}
\end{enumerate}

Here are some pruning strategies of the first type
(\ref{enum:not-tame}). They are immediate consequences of the
conditions of the defining properties of tameness.

\begin{itemize}
  \item If the current state contains an incomplete face of length 3,
    then eliminate all transitions, except for the transition
    that carries the partial plane graph to a partial plane graph that
    is the same in all respects, except that the face has
    become complete.
  \item If the current state contains an incomplete face of length 4,
    then eliminate all transitions except those that lead to
    the possibilities of Section \ref{sec:graphproperty}, Property \ref{definition:tame:4-circuit},
    where in Property \ref{definition:tame:4-circuit}
     each depicted face is interpreted
    as being complete.
  \item Remove all transitions with
    patches whose complete face has length greater than
    $8$.
  \item It is frequently possible to conclude from the examination of a partial
    plane graph that no matter what the terminal position,
    any admissible weight assignment will give total weight greater than
    the target $(14.8)$.  In such cases, the all transitions out of the
    partial plane graph can be pruned.
\end{itemize}

    To take a simple example of the last item, we observe that weights are always
    non-negative, and that the weight of a complete face of
    length $n$ is at least $d(n)$.  Thus, if there are complete
    faces $F_1,\ldots,F_k$ of lengths $n_1,\ldots,n_k$, then
     any admissible weight assignment has total weight at least
    $\sum_{i=1}^k d(n_i)$.  If this number is at least the
    target, then no transitions out of that state need be considered.

More generally, we can apply all of the inequalities in the
definition of admissible weight assignment to the complete portion
of the partial plane graph to obtain lower bounds.  However, we
must be careful, in applying Property
\ref{definition:admissible:excess} of admissible weight
assignments, because vertices that are not adjacent at an
intermediate state may become adjacent in the complete graph.
Also, vertices that do not lie together in a quadrilateral at an
intermediate state may do so in the complete graph.

\begin{figure}[htb]
  \centering
  \includegraphics{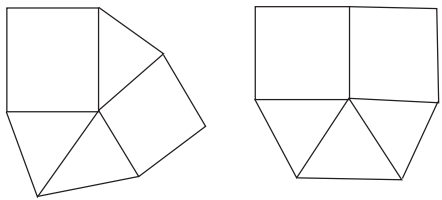}
  \caption{States $I_{3,2}$ and $I_{3,2}'$}
  \label{fig:states}
\end{figure}

Here are some pruning strategies of the second type
(\ref{enum:branch}).
\begin{itemize}
    \item At a given state it is enough to fix one incomplete face and
        one edge of that face and then to follow only the transitions that
        patch along that face and add a complete face along that
        edge. (This is seen from the proof of Lemma
        \ref{lemma:reachable}.)
    \item In leading out from the initial state $I_n$, it is enough
        to follow paths in which every added complete face has
        length at most $n$. (A graph with a face of length $m$,
        for $m>n$, will be also be found downstream from $I_m$.)
    \item Make a list of all type $(p,q)$ with $b(p,q)<14.8$.
    Remove the initial states $I_3$ and $I_4$, and create new initial states
    $I_{p,q}$  ($I_{p,q}'$, $I_{p,q}''$, etc.)
    in the finite state machine.  Define the state $I_{p,q}$ to be
    one consisting of $p+q+1$ faces, with $p$ complete triangles
    and $q$ complete quadrilaterals all meeting at a vertex (and one other incomplete face away from $v$).
    (If there is more than one way to arrange $p$ triangles and
    $q$ quadrilaterals, create states $I_{p,q}$, $I'_{p,q}$,
    $I''_{p,q}$, for each possibility.  See Figure \ref{fig:states}.)
    Put a linear order on states $I_{p,q}$.  In state transitions
    downstream from $I_{p,q}$ disallow any transition that creates
    a vertex of type $(p',q')$, for any $(p',q')$ preceding $(p,q)$
    in the imposed linear order.
\end{itemize}

This last pruning strategy is justified by the following lemma,
which classifies vertices of type $(p,q)$.

\begin{lemma}
Let $A$ and $B$ be triangular or quadrilateral faces that have at
least $2$ vertices in common in a tame graph. Then the faces have
exactly two vertices in common, and an edge is shared by the two
faces.
\end{lemma}

\begin{figure}[htb]
  \centering
  \includegraphics{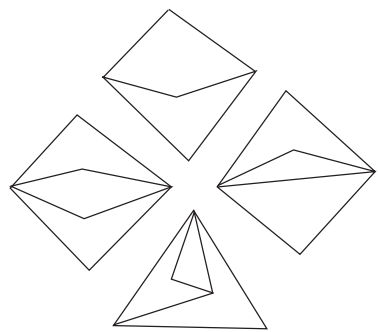}
  \caption{Some impossibilities}
  \label{fig:nonexistant}
\end{figure}

\begin{proof}
Exercise.  Some of the configurations that must be ruled out are
shown in Figure \ref{fig:nonexistant}.  Some properties that are
particularly useful for the exercise are Properties
\ref{definition:tame:3-circuit} and
\ref{definition:tame:4-circuit} of tameness, and Property
\ref{admissible:b} of admissibility.
\end{proof}

Once a terminal position is reached it is checked to see whether
it satisfies all the properties of tameness.

Duplication is removed among isomorphic terminal plane graphs. It
is not an entirely trivial procedure for the computer to determine
whether there exists an isomorphism between two plane graphs. This
is accomplished by computing a numerical invariant of a vertex
that depends only on the local structure of the vertex. If two
plane graphs are properly isomorphic then the numerical invariant
is the same at vertices that correspond under the proper
isomorphism.  If two graphs have the number of vertices with the
same numerical invariants, they become candidates for an
isomorphism.  All possible numerical-invariant preserving
bijections are attempted until an proper isomorphism is found, or
until it is found that none exist. If there is no proper
isomorphism, the same procedure is applied to the opposite plane
graph to find any possible orientation-reversing isomorphism.

This same isomorphism-producing algorithm is used to match each
terminal graph with a graph in the archive.  It is found that each
terminal graph matches with one in the archive. (The archive was
originally obtained by running the finite state machine and making
a list of all the terminal states up to isomorphism that satisfy
the given conditions.)

In this way Theorem \ref{theorem:classification} is proved.

\section{From Stars to Graphs}
\label{sec:startame}

\subsection{Plane Graphs}
\label{sec:stargraph}

A plane graph $G$ is attached to every contravening decomposition
star as follows. For simplicity, take the decomposition star to be
attached to a ball of the packing, centered at the origin.  From
the decomposition star $D$, it is possible to determine the
coordinates of the set $U(D)$ of vertices at distance at most
$2t_0 $ from the origin.

If we draw a geodesic arc on the unit sphere at the origin with
endpoints at the projections of $v_1$ and $v_2$ for every pair of
vertices $v_1$, $v_2\in U(D)$ such that $|v_1|, |v_2|,
|v_1-v_2|\le 2t_0 $, we obtain a plane graph that breaks the unit
sphere into regions called {\it standard regions}. (The arcs do
not meet except at endpoints \cite[Lemma 3.10]{part1}.)  Each
standard region is defined as the closure in the unit sphere of a
connected component of the unit sphere with all arcs removed.

For a given standard region, we consider the arcs forming its
boundary together with the arcs that are internal to the standard
region.  We consider the points on the unit sphere formed by the
endpoints of the arcs, together with the projections to the unit
sphere of vertices in $U$ whose projection lies in the interior of
the region.

\begin{figure}[htb]
  \centering
  \includegraphics{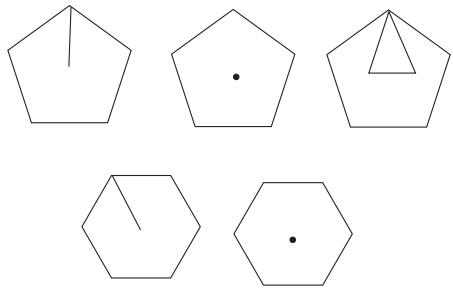}
  \caption{Non-polygonal standard regions}
  \label{fig:aggregates}
\end{figure}

\begin{remark}
The system of arcs and vertices associated with a standard region
in a contravening example must be a polygon, or one of the
following configurations of Figure \ref{fig:aggregates} (see
\cite[Corollary 4.4]{part4}).
\end{remark}

\begin{figure}[htb]
  \centering
  \includegraphics{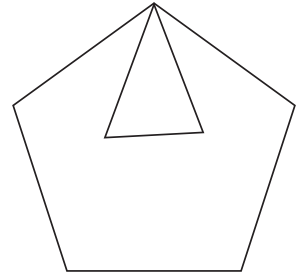}
  \caption{An aggregate forming a pentagon}
  \label{fig:tri-pent}
\end{figure}

\begin{remark}
\label{remark:tri-pent}
 Observe that one case is bounded by a
triangle and a pentagon, and that the others are bounded by a
polygon. Replacing the triangle-pentagon arrangement with the
bounding pentagon and replacing the others with the bounding
polygon, we obtain a partition of the sphere into simple polygons.
Each of these polygons is a single standard region, except in the
triangle-pentagon case (Figure \ref{fig:tri-pent}), which is a
union of two standard regions (a triangle and a eight-sided
region).
\end{remark}

\begin{figure}[htb]
  \centering
  \includegraphics{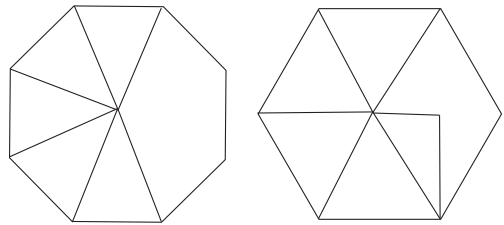}
  \caption{Degree $6$ aggregates}
  \label{fig:degree6}
\end{figure}

\begin{remark}
\label{remark:degree6}
 To simplify further, if we have an
arrangement of six standard regions around a vertex formed from
five triangles and one pentagon, we replace it with the bounding
octagon (or hexagon). See Figure \ref{fig:degree6}.  (It will be
shown in Lemma \ref{lemma:aggregate6} and Section
\ref{sec:impossible} that there is at most one such configuration
in the standard decomposition of a contravening decomposition
star, so we will not worry here about how to treat the case of two
overlapping configurations of this sort.)
\end{remark}

In summary, we have a plane graph that is approximately that given
by the standard regions of the decomposition star, but simplified
to a bounding polygon when one of the configurations of Remarks
\ref{remark:tri-pent} and \ref{remark:degree6} occur.  We refer to
the combination of standard regions into a single face of the
graph as {\it aggregation}.  We call it the plane graph $G = G(D)$
attached to a contravening decomposition star $D$.

Proposition \ref{prop:nonempty} will show the vertex set $U$ is
non-empty and that the graph $G(D)$ is non-empty.

When we refer to the plane graph in this manner, we mean the
combinatorial plane graph in the sense of Section \ref{sec:tame},
as opposed to the embedded metric graph on the unit sphere formed
from the system of geodesic arcs.  Given a vertex $v$ in $G(D)$,
there is a uniquely determined vertex $v(D)$ of $U(D)$ whose
radial projection to the unit sphere determines $v$.  We call
$v(D)$ the {\it corner} in $U(D)$ over $v$.

By definition, the plane graphs in Section \ref{sec:tame} do not
have loops or multiple joins.  The edges of $G(D)$ are defined by
triangles whose sides vary between lengths $2$ and $2t_0 $.  The
angles of such a triangle are strictly less than $\pi$.  This
implies that the edges of the metric graph on the unit sphere
always have arc-length strictly less than $\pi$.  In particular,
the endpoints are never antipodal.  A loop on the combinatorial
graph corresponds to a edge on the metric graph that is a closed
geodesic.  A multiple join on the combinatorial graph corresponds
on the metric graph to a pair of points joined by multiple minimal
geodesics, that is, a pair of antipodal points on the sphere.  By
the arc-length constraints on edges in the metric graph, there are
no loops or multiple joins in the combinatorial graph $G(D)$.

In Definition \ref{definition:tame}, a plane graph satisfying a
certain restrictive set of properties is said to be {\it tame}. If
a plane graph $G(D)$ is associated with a contravening
decomposition star $D$, we call $G(D)$ a contravening plane graph.

\begin{theorem} \label{theorem:contravene}
 Every contravening plane graph is tame.
\end{theorem}

This theorem is one of the main results of this paper.  Its proof
occupies Sections \ref{sec:contraproof} through
\ref{sec:fourcycle}.

In Theorem \ref{theorem:classification}, the tame graphs  are
classified up to isomorphism. As a corollary, we have an explicit
list containing all contravening plane graphs.

\section{Contravention is tame}
\label{sec:contraproof}

This section begins the proof of Theorem \ref{theorem:contravene}.

To prove Theorem \ref{theorem:contravene}, it is enough to show
that each defining property of tameness is satisfied.  This is the
substance of results in the next sections.  First, we prove the
promised non-degeneracy result.

\begin{proposition}
\label{prop:nonempty} The construction of Section
\ref{sec:stargraph} associates a (nonempty) plane graph with at
least two faces to every decomposition star $D$ with
$\sigma(D)>0$.
\end{proposition}

\begin{proof}
First we show that decomposition stars with $\sigma(D)>0$ have
non-empty vertex sets $U$. (Recall that $U$ is the set of vertices
of distance at most $2t_0$ from the center).  The vertices of $U$
are used in \cite{formulation} to create all of the structural
features of the decomposition star: quasi-regular tetrahedra,
quarters, and so forth.  If $U$ is empty, the decomposition star
is a solid containing the ball $B(t_0)$ of radius $t_0$, and its
score is just
    $$
    \begin{array}{lll}
    \sigma(D) &= \vor(D)\\
              & = -4\doct \Vol(D) + 4\pi/3 \\
              &< -4\doct\Vol(B(t_0)) + 4\pi/3
              &< 0.
    \end{array}
    $$
By hypothesis, $\sigma(D)>0$.  So $U$ is not empty.

The results of \cite[Eqn 3.11]{formulation} show that the function
$\sigma$ can be expressed as a sum of terms $\sigma_R$ indexed by
the standard regions $R$.  It is proved in \cite[Lemma 3.13
]{formulation} that $\sigma_R\le0$, unless $R$ is a triangle.
Thus, a decomposition star with positive score must have at least
one triangle. Its complement contains a second standard region.
Even after we form aggregates of distinct standard regions to form
the simplified plane graph (Remarks \ref{remark:tri-pent} and
\ref{remark:degree6}), there certainly remain at least two faces.
\end{proof}

\begin{proposition} The plane graph of a contravening decomposition
star satisfies Property \ref{definition:tame:length} of tameness:
The length of each face is (at least $3$ and) at most $8$.
\end{proposition}

\begin{proof} By the construction of the graph, each face has at
least three edges.  The upper bound of $8$ faces is
\cite[Corollary 4.4]{part4}.
\end{proof}

\begin{proposition} The plane graph of a contravening decomposition star
satisfies Property \ref{definition:tame:3-circuit} of tameness:
Every $3$-circuit is a face or the opposite of a face.
\end{proposition}

\begin{proof}
The simplifications of the plane graph in Remarks
\ref{remark:tri-pent} and \ref{remark:degree6}  do not produce any
new $3$-circuits. (See the accompanying figures.) The result is
proved in \cite[Lemma 3.7]{part1}. The proof is purely geometric:
it expresses a geometric impossibility true of all decomposition
stars.
\end{proof}

\begin{proposition} Contravening graphs satisfy Property
\ref{definition:tame:degree} of tameness: The degree of every
vertex is at least $2$ and at most $6$.
\end{proposition}

\begin{proof}
The statement that degrees are at least $2$ trivially follows
because each vertex lies on at least one polygon, with two edges
at that vertex.

The impossibility of a vertex of degree $7$ or more is found in
the proof of \cite[Lemma 6.2]{part3}.
\end{proof}

\bigskip
\section{More Tame Properties} 
\label{sec:moretame}

This section continues in the proof that all contravening plane
graphs are tame.

\subsection{Linear Programs} 
\label{sec:2.2}  To continue with the proof that contravening
plane graphs are tame, we need to introduce some more notation and
methods.

The function $\sigma$ can be written is a sum of functions on the
space of decomposition stars, indexed by standard regions.  (See
\cite[Eqn 3.11]{formulation}.)  Write $\sigma_R$ for the term
indexed by the standard region $R$.  Then
    $$\sigma(D) = \sum_R \sigma_R(D),$$
for every decomposition star.

Let $\zeta = 1/(2\arctan(\sqrt{2}/5)$. Let $\sol(R)$ denote the
solid angle of a region $R$.  We write $\tau_R$ for the following
modification of $\sigma_R$:
    $$\tau_R(D) = \sol(R)\zeta\pt - \sigma_R(D)$$
and
    $$\tau(D) = \sum\tau_R(D) = 4\pi\zeta\pt - \sigma(D).$$
Since $4\pi\zeta\pt$ is a constant, $\tau$ and $\sigma$ contain
the same information, but $\tau$ is often more convenient to work
with. A contravening decomposition star satisfies
    \begin{equation}
    \tau(D) \le 4\pi\zeta\pt -8\,\pt = (4\pi\zeta-8)\pt.
    \end{equation}
The constant $\squander$  (and its upper bound $14.8\,\pt$) will
occur repeatedly in the discussion that follows.

Let ${\cal R}(F)$ be the set of standard regions associated with a
face $F$ in the plane graph.  We occasionally combine the terms
$\sigma_R$ associated with the same face.  The set contains a
single standard region, unless multiple standard regions have been
aggregated
 according to the constructions of Remarks
\ref{remark:tri-pent} and \ref{remark:degree6}.
Set
    $$
    \sigma_F =\sum_{R\in {\cal R}(F)} \sigma_R
    $$
and
    $$
    \tau_F =\sum_{R\in{\cal R}(F)}\tau_R.
    $$

In \cite[Section 4.4]{part4} a natural number $n(R)$ is attached
to the standard region $R$.  Lower bounds on $\tau_R(D)$ are
obtained in \cite[Theorem 4.4]{part4} that depend only on $n(R)$.
Write $t_n$ for the constants in this theorem.  The same theorem
gives upper bounds $s_n$ for $\sigma_R(D)$ that depend only on
$n(R)$.

It is also helpful to bear in mind several bounds on angles in the
following discussion.  Every internal angle of every standard
region is at least $0.8638$ \cite[10.1.3]{part3}.  Every internal
angle of every standard region, except for triangles, is at least
$1.153$ \cite[4.3]{part3}.

A decomposition star $D$ determines a set of vertices $U(D)$ that
are of distance at most $2t_0$ from the center of $D$.  Three
consecutive vertices $p_1$, $p_2$, and $p_3$ of a standard region
are determined as the projections to the unit sphere of three
corners $v_1$, $v_2$, and $v_3$, respectively in $U(D)$. If the
internal angle of a standard region is less than $1.32$, then it
follows by an interval arithmetic calculation that
$|v_1-v_3|\le\sqrt{8}$, see \cite[Lemma 3.11.4]{part4}.

We are to the point in the proof of the Kepler Conjecture where it
becomes necessary to use inequalities proved by interval
arithmetic in a serious way.   To use these inequalities
systematically, we combine inequalities into linear programs and
solve the linear programs on computer.  At first, our use of
linear programs will be light, but our reliance will become
progressively strong as this paper develops.

To start out, we will make use of the inequalities
\cite[10.3.4-8]{part3} and \cite[4.1]{part3}.  These inequalities
give lower bounds on $\tau_R(D)$ when $R$ is a triangle or
quadrilateral. To obtain lower bounds through linear programming,
we take a linear relaxation.  Specifically, we introduce a linear
variable for each function $\tau_R$ and a linear variable for each
internal angle $\alpha_R$. We substitute these linear variables
for the nonlinear functions $\tau_R(D)$ and nonlinear internal
angle function into the given inequalities in \cite{part3}.  Under
these substitutions, the inequalities become linear.   Given $p$
triangles and $q$ quadrilaterals at a vertex, we have the linear
program to minimize the sum of the (linear variables associated
with) $\tau_R(D)$ subject to the constraint that the (linear
variables associated with the)  angles at the vertex sum to at
most $d$. Linear programming yields a lower bound $\tauLP(p,q,d)$
to this minimization problem.  This gives  a lower bound to the
corresponding constrained sum of nonlinear functions $\tau_R$.

Similarly, the inequalities from the same sections of \cite{part4}
yield upper bounds $\sLP(p,q,d)$ on the sum of $p+q$ functions
$\sigma_R$, with $p$ standard regions $R$ that are triangular, and
another $q$ that are quadrilateral.  These linear programs find
their first application in the proof of the following proposition.

\begin{proposition}
    The plane graph of a contravening decomposition star satisfies
    Property \ref{definition:tame:score} of tameness:
    $$\sum_F c(len(F)) \ge 8.$$
\end{proposition}

\begin{proof}
We will show that
    \begin{equation}
    c(len(F))\,\pt \ge \sigma_F(D)
    \label{eqn:sigma}
    \end{equation}
    Assuming this, the result follows for contravening
    stars $D$:
    $$
    \begin{array}{lll}
        \sum_F c(len(F)) \,\pt &\ge \sum_F \sigma_F(D) \\
            &= \sigma(D) \ge 8\,\pt.
    \end{array}
    $$

We consider three cases for Inequality \ref{eqn:sigma}. In the
first case, assume that the face $F$ corresponds to exactly one
standard region in the decomposition star.  In this case,
Inequality \ref{eqn:sigma} follows directly from the bounds
\cite[Theorem 4.4]{part4}.

In the second case, assume we are in the context of a pentagon $F$
formed in Remark \ref{remark:tri-pent}.  Then, again by
\cite[4.4]{part4}, we have
$$\sigma_F(D) \le (c(3)+c(8))\,\pt \le c(5)\,\pt.$$
(Just examine the constants $c(k)$.)

In the third case, we consider the situation of Remark
\ref{remark:degree6}.  The six standard regions score at most
$s_5+\sLP(5,0,2\pi-1.153)< c(8)\,\pt$.
\end{proof}

\begin{proposition}
    \label{proposition:wttau}
    Let $F$ be a face of a contravening plane graph $G(D)$.
    Then
    $$\tau_F(D) \ge d(len(F))\pt.$$
\end{proposition}

\begin{proof} Similar.
\end{proof}

\begin{figure}[htb]
  \centering
  \includegraphics{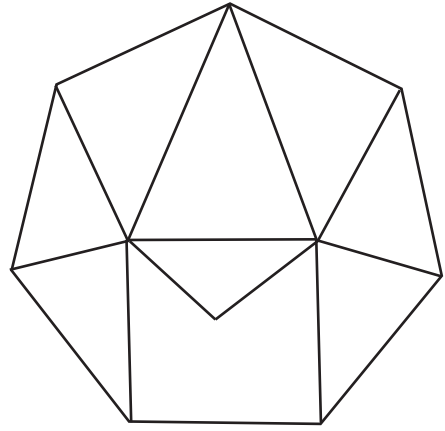}
  \caption{Adjacent vertices of degree $6$}
  \label{fig:doubledegree6}
\end{figure}

\begin{lemma}
    \label{lemma:aggregate6}
    Consider the aggregate in Remark \ref{remark:degree6}.
    \begin{enumerate}
        \item There are at most two vertices (of a standard region) of degree six.  If
        there are two, then they are adjacent vertices on a
        pentagon (Figure \ref{fig:doubledegree6}).
        \item If a vertex of a pentagonal standard region has degree six, then the
        aggregate $F$ of the six faces satisfies
            $$
            \begin{array}{lll}
            \sigma_F(D) < s_8,\\
            \tau_F(D) > t_8.
            \end{array}
            $$
    \end{enumerate}
\end{lemma}

\begin{proof}
We begin with the second part of the lemma. The sum  $\tau_F(D)$
over these six standard regions is at least
    $$t_5+\tauLP(5,0,2\pi-1.153)> t_8.$$
Similarly,
    $$s_5+\sLP(5,0,2\pi-1.153)<s_8.$$
We note that there can be at most one exceptional region with a
vertex of degree six.  Indeed, if there are two, then they must
both be vertices of the same pentagon:
    $$t_8+t_5>\squander.$$
Such a second vertex on the octagonal aggregate leads to one of
the follow constants greater than $\squander$.  These same
constants show that such a second vertex on a hexagonal aggregate
must share two triangular faces with the first vertex of degree
six.
$$\begin{array}{lll}
    t_8 &+\tauLP(4,0,2\pi-1.32-0.8638),\quad\text{or}\\
    t_8 &+1.47\,\pt+\tauLP(4,0,2\pi-1.153-0.8638) ,\\
    t_8 &+\tauLP(5,0,2\pi-1.153) .
\end{array}
$$
(The relevant inequalities involving the constants $1.47$ and
$1.32$ are found at \cite[Lemma 5.12.1]{part4}.)

\end{proof}

\begin{proposition} A contravening plane graph satisfies Property
\ref{definition:tame:degreeE} of tameness: If a vertex is
contained in an exceptional face,
        then the degree of the vertex is at most $5$.
\end{proposition}


\begin{proof}
Assume for a contradiction that $v$ is a vertex of degree six that
lies in an exceptional region.  There are several cases according
to the number $k$ of triangular regions at the vertex.

{\bf($k\le2$)} If there are at least four nontriangular regions at
the vertex, then the sum of interior angles around the vertex is
at least $4(1.153)+2(0.8638)>2\pi$, which is impossible.

{\bf($k=3$)} If there are three nontriangular regions at the
vertex, then $\tau(D)$ is at least
$2t_4+t_5+\tauLP(3,0,2\pi-3(1.153))>\squander$.

{\bf($k=4$)} If there are two exceptional regions at the vertex,
then $\tau(D)$ is at least
$2t_5+\tauLP(4,0,2\pi-2(1.153))>\squander$.

If there are two nontriangular regions at the vertex, then
$\tau(D)$ is at least  $t_5+\tauLP(4,1,2\pi-1.153)>\squander$.

We are left with the case of five triangular regions and one
exceptional region.

When there is an exceptional standard region at a vertex of degree
six, we claim that the exceptional region must be a pentagon. If
the region is a heptagon or more, then $\tau(D)$ is at least
$t_7+\tauLP(5,0,2\pi-1.153) > \squander$.

If the standard region is a hexagon, then $\tau(D)$ is at least
$t_6 + \tauLP(5,0,2\pi-1.153) > t_9$. Also,
$s_6+\sLP(5,0,2\pi-1.153) < s_9$. The aggregate of the six
standard regions is a nonagon. The argument of \cite[Section
4.6]{part4} extends to this context to give the bound of $8\,\pt$.

\end{proof}

\section{Weight Assignments}
\label{sec:weight}

The purpose of this section is to prove the existence of a good
admissible weight assignment for contravening plane graphs.

\begin{theorem}  Every contravening planar graph has an admissible
weight assignment of total weight less than $14.8$.
\end{theorem}

Given a contravening decomposition star $D$, we define a weight
assignment $w(F)$ by
    $$w(F) = \tau_F(D)/\pt.$$
Since $D$ contravenes,
    $$
    \begin{array}{lll}
    \sum_F w(F) &= \sum_F \tau_F(D)/\pt \\
            &= \tau(D)\le\squander/\pt \\
        &< 14.8.
    \end{array}
    $$
The challenge of the theorem will be to prove that $w(F)$, when
defined by this formula, is admissible.

The next three lemmas establish that this definition of $w(F)$ for
contravening plane graphs satisfies the first three defining
properties of an admissible weight assignment.

\begin{lemma}  Let $F$ be a face of length $n$ in a contravening plane graph.
Define $w(F)$ as above. Then
        $w(F) \ge d(n)$.
\end{lemma}

\begin{proof} This is Proposition \ref{proposition:wttau}.
\end{proof}

\begin{lemma} Let $v$ be a vertex of type $(p,q)$ in a
contravening plane graph.  Define $w(F)$ as above. Then
        $$\sum_{v\in F} w(F) \ge b(p,q).$$
\end{lemma}

\begin{proof} This is \cite[Proposition 5.2]{part3}.
\end{proof}

\begin{lemma} Let $V$ be any set of vertices of type $(5,0)$ in a
contravening plane graph.  Define $w(F)$ as above.
        If the cardinality of $V$ is $k\le 4$,
        then
        $$\sum_{V\cap F\ne\emptyset} w(F) \ge 0.55 k.$$
\end{lemma}

\begin{proof} This is \cite[Proposition 5.3]{part3}.
\end{proof}

The following proposition establishes the final property that
$w(F)$ must satisfy to make it admissible.

\begin{proposition}
        \label{proposition:excess}
        Let $V$ be any separated set of vertices in a contravening plane graph.
        Define $w(F)$ as above.
        Then
        $$\sum_{V\cap F\ne\emptyset} (w(F) -d(len(F)))
            \ge \sum_{v\in V} a(tri(v)),$$
        where $tri(v)$ denotes the number of triangles containing
        the vertex $v$.
\end{proposition}

The proof will occupy the rest of this section, as well as all of
Section \ref{sec:aggregate}. Since the degree of each vertex is
five, and there is at least one face that is not a triangle at the
vertex, the only constants $tri(v)$ that arise are
    $$tri(v) \in\{0,\ldots,4\}$$
We will prove that in a contravening plane graph that the
Properties (1) and (4) of a separated set are incompatible with
the condition $tri(v)\le 2$, for some $v\in V$.

This will allows us to assume that $$tri(v)\in\{3,4\},$$ for all
$v\in V$.  These cases will be treated in Sections \ref{sec:tri3}
and \ref{sec:tri4}.

The outline of the proof of Proposition \ref{proposition:excess}
is the following. Let $D$ be a contravening decomposition star.
First, assume that every face meeting $V$ corresponds to a single
standard region of $D$. (That is, we exclude the aggregates
created in Remark \ref{remark:tri-pent} and Remark
\ref{remark:degree6}.) Under this hypothesis, we show that
$tri(v)\le 2$ is incompatible with the hypothesis.  This allows us
to assume that $tri(v)\in\{3,4\}$.
 We express the desired inequality as a sum of
inequalities indexed by $V$, where the inequality indexed by $v\in
V$ depends only on (an explicit description of) the function
$\tau_R(D)$ for standard regions $R$ corresponding to a face
containing $v$.   These individual inequalities are proved first
for $tri(v)=3$, and then for $tri(v)=4$.

Finally in Section \ref{sec:aggregate}, we return to the
aggregates of Remark \ref{remark:tri-pent} and Remark
\ref{remark:degree6} and show that the conclusions can be extended
to aggregate faces as well.

\subsection{Proof that $tri(v)>2$}
\label{sec:2.4}

In this subsection $D$ is a contravening decomposition star with
associated graph $G(D)$.  Let $V$ be a set of vertices that
satisfies the conditions of Proposition \ref{proposition:excess}.
Let $v$ be a vertex in $V$ such that none of its faces is an
aggregate in the sense of Remarks \ref{remark:tri-pent} and
\ref{remark:degree6}.

\begin{lemma}  Under these conditions, for every $v\in V$,
$tri(v)>1$.
\end{lemma}

\begin{proof}
If there are $p$ triangles, $q$ quadrilaterals, and $r$ other
faces, then
    $$
    \begin{array}{lll}
    \tau(D) &\ge\sum_{v\in R}\tau_R(D)\\
        &\ge r\, t_5 + \tauLP(p,q,2\pi-r(1.153)).
    \end{array}
    $$ If there is a vertex $w$ that is
not on any of the faces containing $v$, then the sum of
$\tau_F(D)$ over the faces containing $w$ yield an additional
$0.55\,\pt$ by \cite[5.2]{part3}. We run the linear programs for
each $(p,q,r)$ and find that the bound is always greater than
$\squander$.  This implies that $D$ cannot be contravening.
$$\begin{matrix}
    (p,q,r)&\hbox{\it lower bound }\\
    &\\
    (0,5,0)&22.27\,\pt \text{ \cite[5.2]{part3}}\\
    (0,q,r\ge1)& t_5+4 t_4\approx 14.41\,\pt\\
    (1,4,0) &17.62\,\pt \text{ \cite[5.2]{part3}}\\
    (1,3,1) &t_5 + 12.58\,\pt \ (\tauLP)\\
    (1,2,2) &2t_5 + 7.53\,\pt \ (\tauLP)\\
    (1,q,r\ge3)& 3 t_5 + t_4\\
\end{matrix}
$$
\end{proof}

\begin{lemma} Under these same conditions, for every $v\in V$,
$tri(v)>2$.
\end{lemma}

\begin{proof}
Assume that $tri(v)=2$.  We will show that this implies that $D$
is not contravening.  Let $e$ be the number of exceptional faces
at $v$.

 The constants $0.55\,\pt$ and $4.52\,\pt$ used
throughout the proof come from \cite[5.3]{part3}; The constants
$t_n$ comes from \cite[4.4]{part4}.

 ($e=3$): First, assume that there are three
exceptional faces around vertex $v$. They must be pentagons
($2t_5+t_6>\squander$). The aggregate of the five faces is an
$m$-gon (some $m\le11$).  If there is a vertex not on this
aggregate, use $3t_5+0.55\,\pt>\squander$. So there are at most
nine triangles away from the aggregate.
    $$
    \sigma(D) \le 9\,\pt + (3 s_5+2\,\pt) < 8\,\pt.
    $$

The argument if there is a quad, pentagon, and hexagon is the same
$(t_4+t_6=2t_5,s_4+s_6=2s_5)$.

($e=2$): Assume next that there are two pentagons and a
quadrilateral around the vertex. The aggregate of the two
pentagons, quadrilateral, and two triangles is an $m$-gon (some
$m\le10$). There must be a vertex not on the aggregate of five
faces, for otherwise we have
    $$
    \sigma(D) \le 8\,\pt+(2s_5+2\,\pt)<8\,\pt.
    $$

The interior angle of one of the pentagons is at most $1.32$.  For
otherwise, $\tauLP(2,1,2\pi-2(1.32))+2t_5+0.55\,\pt>\squander$.

The calculation \cite[Lemma 5.12.2]{part4} shows that any pentagon
$R$ with an internal angle less than $1.32$ yields $\tau_R(D)\ge
5.66\,\pt$. If both pentagons have an interior angle $<1.32$ the
lemma follows easily from this calculation:
$2(5.66)\,\pt+\tauLP(2,1,2\pi-2(1.153))+0.55\,\pt>\squander$. If
there is one pentagon with angle $>1.32$, we then have
$5.66\,\pt+\tauLP(2,1,2\pi-1.153-1.32)+t_5+0.55\,\pt>\squander$.

($e=1$): Assume finally that there is one exceptional face at the
vertex. If it is a hexagon (or more), we are done
$t_6+\tauLP(2,2,2\pi-1.153)>\squander$. Assume it is a pentagon.
The aggregate of the five faces at the vertex is bounded by an
$m$-circuit (some $m\le9$). If there are no more than $11$
quasi-regular tetrahedra outside the aggregate, the score is at
most $(1+2(4.52))\,\pt+s_5+\sLP(2,2,2\pi-1.153)<8\,\pt$
(\cite[5.3]{part3}). So we may assume that there are at least
three vertices not on the aggregate.

If the interior angle of the pentagon is greater than $1.32$, we
have
$$\tauLP(2,2,2\pi-1.32) +3(0.55)\,\pt +t_5 > \squander;$$
and if it is less than $1.32$, we have by \cite[Lemma
5.12.1]{part4}
    $$
    \begin{array}{lll}
        \tauLP(2,2,2\pi-1.153)&+3(0.55)\pt+1.47\,\pt+t_5 \\
            &> \squander.
    \end{array}
    $$
\end{proof}

\begin{lemma} The bound $tri(v)>2$ holds if $v$ is a vertex
of an aggregate face.
\end{lemma}

\begin{proof}
The exceptional region enters into the preceding two proofs in a
purely formal way.  Pentagons enter through the bounds
    $$t_5,\ s_5,\ 1.47\,\pt,\ 5.66\,\pt$$
and angles $1.153$, $1.32$.  Hexagons enter through the bounds
    $$t_6,\ s_6$$
and so forth.  These bounds hold for the aggregate faces.  Hence
the proofs hold for aggregates as well.
\end{proof}

\subsection{Bounds when $tri(v)\in\{3,4\}$.  } 
\label{sec:2.7}

In this subsection $D$ is a contravening decomposition star with
associated graph $G(D)$.  Let $V$ be a separated set of vertices.
For every vertex  $v$ in $V$, we assume that none of its faces is
an aggregate in the sense of Remarks \ref{remark:tri-pent} and
\ref{remark:degree6}.  We assume that there are three or four
triangles containing $v$, for every $v\in V$.

To prove the Inequality \ref{definition:admissible:excess} in the
definition of admissible weight assignments, we will rely on the
following reductions. Define an equivalence relation on
exceptional faces by $F\sim F'$ if there is a sequence
$F_0=F,\ldots, F_r=F'$ of exceptional faces such that consecutive
faces share a vertex of type $(3,0,2)$.  Let ${\cal F}$ be an
equivalence class of faces.

\begin{lemma} Let $V$ be a separated set of vertices.  For every
equivalence class of exceptional faces $\cal F$, let $V({\cal F})$
be the subset of $V$ whose vertices lie in the union of faces of
${\cal F}$. Suppose that for every equivalence class $\cal F$, the
Inequality \ref{definition:admissible:excess} (in the definition
of admissible weight assignments) holds for $V({\cal F})$. Then
the Inequality holds for $V$.
\end{lemma}

\begin{proof}
By construction, each vertex in $V$ lies in some $F$, for an
exceptional face.  Moreover, the separating property of $V$
insures that the triangles and quadrilaterals in the inequality
are associated with a well-defined  ${\cal F}$. Thus, the
inequality for $V$ is a sum of the inequalities for each $V({\cal
F})$.
\end{proof}

\begin{lemma}
\label{lemma:split}
 Let $v$ be a vertex at which there are $p$
triangles, $q$ quadrilaterals, and $r$ other faces.  Suppose that
for some $p'\le p$ and $q'\le q$, we have
    $$\tauLP(p',q',\alpha) > ( p' d(3) + q' d(4) + a(p))\,\pt$$
for some upper bound $\alpha$ on the angle occupied by $p'$
triangles and $q'$ quadrilaterals at $v$.  Suppose further that
the Inequality \ref{definition:admissible:excess} (in the
definition of admissible weight assignments) holds for $V' =
V\setminus \{v\}$. Then the inequality holds for $V$.
\end{lemma}

\begin{proof}  Let $F_1,\ldots,F_m$, $m={p'+q'}$, be faces corresponding
to the triangles and quadrilaterals in the lemma.  The hypotheses
of the lemma imply that
    $$\sum_{1}^{m} (\tau_{F_i}(D) - d(len(F_i))) > a(p).$$
Clearly, the Inequality for $V$ is the sum of this inequality, the
inequality for $V'$, and $d(n)\ge0$.
\end{proof}

\subsection {$tri(v)=3$}
\label{sec:tri3}

 To continue further, we need to recall from
\cite{part4} how $\tau_R(D)$, for exceptional standard regions,
can be written as a sum of terms in a way that reflects the
internal geometry of $D$ over $R$.  We have prepared a summary of
the relevant material from \cite{part4} in Appendix
\ref{app:part4}. We recommend that the reader review Appendix
\ref{app:part4} before reading the next lemma.  We use terminology
from that appendix (penalties, loops, erasure, flat quarters,
$D(n,k)$, $\dloop(n,k)$, $\hat\sigma$).

\begin{lemma} Property \ref{definition:admissible:excess}  of admissibility holds when
$tri(v)=3$ for all $v\in V$.
\end{lemma}

\begin{proof}
When there are three triangles at a vertex of degree five along
with an exceptional face, the relevant constant from the
definition of admissibility is $a(3)=1.4$. Erase all upright
quarters except those in loops. The result follows from
\cite[Lemma 5.12.1]{part4}. and the subsequent remark if the
internal angle at the vertex of an exceptional region is
$\le1.32$.

Assume the angles on the exceptional regions are $\ge1.32$. If
there are three triangles, a quad, and an exceptional face at the
vertex, then we have $\tauLP(3,1,2\pi-1.32)>1.4\,\pt+t_4$.

 The final case is three triangles and
two exceptional faces at the vertex. There can be no heptagons and
at most one hexagon $2t_6 = t_5+t_7>\squander$. In these
calculations, the flat quarters are scored by $\hat\sigma$. By
Appendix \ref{app:part4}, each flat quarter contributes $D(3,1)$
to $t_n$. The constant $a(tri(v))\,\pt$ estimates what the flat
quarter and surrounding quasi-regular tetrahedra squander in
excess of $D(3,1)$.  Most of the calculations now follow directly
from the interval calculations of Appendix \ref{app:tri3}. A few
additional comments are needed if there is a loop masking a flat
quarter.  If the constant $\dloop(n,k)$ exceeds $r(1.4)\,\pt$,
where $r$ is the greatest possible number of masked flat quarters
within the loop, then the result follows.  This gives every case
except $(n,k)=(4,1),(5,0)$. If $(n,k)=(4,1)$, we erase and take
the penalty $\piF=2\xiV+\xiG$ whenever the internal angle is at
least $1.32$. This penalty appears in Appendix \ref{app:tri3}.

If $(n,k)=(5,0)$, there are various possibilities.  If there is
one masked flat quarter, use $\dloop(5,0)>1.4\,\pt$.  If there are
two masked flat quarters, erase with penalty $\piF=2\xiV+\xiG/2$
at each flat quarter.  If the internal angle is at most  $1.32$,
use
$$3.07\,\pt-D(3,1)>1.4\,\pt+\piF.$$
If the internal angle is at least $1.32$, use the interval
calculations Appendix \ref{app:tri3} with penalties.
\end{proof}

\subsection{ $tri(v) = 4$ } 
\label{sec:tri4}

If there are four triangles and one exceptional region at a
vertex, we justify the constant $a(tri(4))=1.5\,\pt$, in the
definition of admissible weight assignment. These calculations are
similar to those of Section \ref{sec:tri3} and are based on the
interval calculations Appendix \ref{app:tri4}. Details are left to
the reader. It is simpler because there is only one exceptional
region, and it is not necessary to separate the case according to
the size of the internal angle (at least $1.32$ and at most
$1.32$).

\section{The Aggregate Cases}
\label{sec:aggregate}

\subsection{The non-polygonal standard regions}

\begin{lemma}
Consider a separated set of vertices $V$ on an aggregated face $F$
as in Remark \ref{remark:tri-pent}.  The Inequality
\ref{definition:admissible:excess} holds (in the definition of
admissible weight assignments):
    $$\sum_{V\cap F\ne\emptyset} (w(F) -d(len(F)))
            \ge \sum_{v\in V} a(tri(v)).$$
\end{lemma}

\begin{proof}
We may assume that $tri(v)\in\{3,4\}$.

 First consider the
aggregate of Remark \ref{remark:tri-pent} of a triangle and
eight-sided region, with pentagonal hull $F$. There is no other
exceptional region in a contravening decomposition star with this
aggregate:
    $$t_8 + t_5 > \squander.$$
A separated set of vertices $V$ on $F$ has cardinality at most
$2$.  This gives the desired bound $$t_8 > t_5 + 2 (1.5)\,\pt.$$

Next, consider the aggregate of a hexagonal hull with an enclosed
vertex.  Again, there is no other exceptional face. If there are
at most $k\le 2$ vertices in a separated set, then the result
follows from
    $$t_8 > t_6 + k (1.5)\,\pt.$$
There are at most three vertices in $V$ on a hexagon, by the
non-adjacency conditions defining $V$.  A vertex $v$ can be
removed from $V$ if it does not lie on a flat quarter (Lemma
\ref{lemma:split} and Appendix \ref{app:tri4}).  However, if there
is an enclosed vertex $w$, it is impossible for there to be three
nonadjacent vertices, none lying along a flat quarter:
    $$\E(2,2,2,\sqrt8,\sqrt8,\sqrt8,2t_0,2t_0,2)>2t_0.$$
(The definition of $\E$ can be reviewed at
\cite[Section~1]{formulation}.)

Finally consider the aggregate of a pentagonal hull with an
enclosed vertex.  There are at most $k\le2$ vertices in a
separated set in $F$.  There is no other exceptional region:
    $$t_7 + t_5 > \squander.$$
The result follows from
    $$t_7 > t_5 + 2(1.5)\,\pt.$$
\end{proof}

\begin{lemma}
Consider a separated set of vertices $V$ on an aggregate face of a
contravening plane graph as in Remark \ref{remark:degree6}.  The
Inequality \ref{definition:admissible:excess} holds in the
definition of admissible weight assignments.
\end{lemma}

\begin{proof}
There is at most one exceptional face in the plane graph:
    $$t_8 + t_5 > \squander.$$
Assume first that aggregate face is an octagon (Figure
\ref{fig:degree6}). At each of the vertices of the face that lies
on a triangular standard region in the aggregate, we can remove
the vertex from $V$ using Lemma \ref {lemma:split} and the
estimate
    $$\tauLP(4,0,2\pi-2 (0.8638)) > 1.5\,\pt.$$
This leaves at most one vertex in $V$, and it lies on a vertex of
$F$ which is ``not aggregated,'' so that there are five standard
regions of the underlying decomposition star at that vertex, and
one of those regions is pentagonal.  The value $a(4)=1.5\,\pt$ can
be estimated at this vertex in the same way it is done for a
non-aggregated case. (We use the calculations of Appendix
\ref{app:tri4}.)

Now consider the case of an aggregate face that is a hexagon
(Figure \ref{fig:degree6}).  The argument is the same: we reduce
to $V$ containing a single vertex, and argue that this vertex can
be treated by the inequalities of Appendix \ref{app:tri4}.
\end{proof}

\section{Four-cycles} 
\label{sec:fourcycle}

It has now been verified that all the properties of tameness
except one have been verified.  This section verifies that
contravening plane graphs satisfy this final property, and in this
way we complete the proof that all contravening plane graphs are
tame.

\subsection{Possible four-circuits}

Every $4$-circuit divides a plane graph into two aggregates of
faces that we may call the interior and exterior.  We call
vertices of the faces in the aggregate that do not lie on the
$4$-cycle {\it enclosed vertices}.  Thus, every vertex lies in the
$4$-cycle, is enclosed over the interior, or is enclosed over the
exterior.

Proposition \cite[4.2]{part1} asserts that either the interior or
the exterior has at most $1$ enclosed vertex.  (The proof does not
require the hypothesis that the graph be contravening. It
expresses a geometric constraint holding for all decomposition
stars.)  When choosing which aggregate is to be called the
interior, we may make our choice so that the interior contains at
most $1$ enclosed vertex.  With this choice, we have the following
proposition.

\begin{proposition}
Let $D$ be a contravening plane graph. The interior of a
$4$-circuit surrounds is one of the aggregates of faces shown in
Property
    \ref{definition:tame:4-circuit} of tameness.
\end{proposition}

\begin{proof}
If there are no enclosed vertices, then the only possibilities are
for it to be a single quadrilateral face or a pair of adjacent
triangles.

Assume there is one enclosed vertex $v$.  If $v$ is connected to
$3$ or $4$ vertices of the quadrilateral, then that possibility is
listed as part of the conclusion.

If $v$ is connected to $2$ opposite vertices in the $4$-cycle,
then the vertex $v$ has type $(0,2)$ and the bounds of
\cite[Proposition 5.2]{part3} show that the graph cannot be
contravening.

If $v$ is connected to $2$ adjacent vertices in the $4$-cycle,
then we appeal to Lemma \ref{lemma:nobad4} to conclude that the
graph is not contravening.

If $v$ is connected to $0$ or $1$ vertices, then we appeal to
\cite[Lemma 2.2]{part3}, which shows that this creates an
impossible geometric situation.  This completes the proof.
\end{proof}

\subsection{A non-contravening Four-circuit}
\label{sec:impossible}

\begin{figure}[htb]
  \centering
  \includegraphics{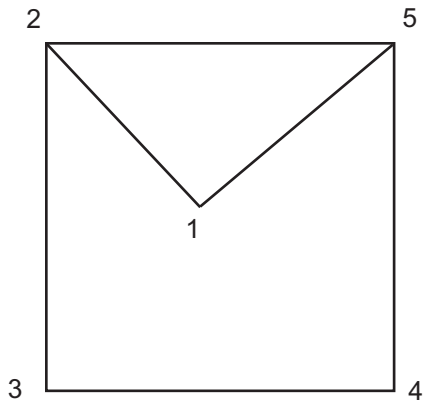}
  \caption{A non-contravening $4$-circuit}
  \label{fig:no4circuit}
\end{figure}

This subsection rules out the existence of a particular
four-circuit on a contravening plane graph.  The interior of the
circuit consists of two faces: a triangle and a pentagon.  The
circuit and its enclosed vertex are show in Figure
\ref{fig:no4circuit} with vertices marked $p_1,\ldots,p_5$.  The
vertex $p_1$ is the enclosed vertex, the triangle is
$(p_1,p_2,p_5)$ and the pentagon is $(p_1,\ldots,p_5)$.

Suppose that $D$ is a decomposition star whose associate graph
contains such a $4$-circuit.  Recall that $D$ determines a set
$U(D)$ of vertices in Euclidean $3$-space of distance at most
$2t_0$ from the center (the origin) of $D$, and that each vertex
$p_i$ can be realized geometrically as a point on the unit sphere
at the origin, obtained as the radial projection of some $v_i\in
U(D)$.

 In the rest of this subsection, we
assume that the reader is familiar with techniques in
\cite{part4}. (In particular, the reader will need to be familiar
with {\it geometric considerations} and the details about how the
function $\tau_R(D)$ is estimated for a pentagonal standard region
$R$. This is reviewed in Appendix \ref{app:part4}.)

\begin{lemma}  One of the edges $(v_1,v_3)$, $(v_1,v_4)$ has
length less than $2\sqrt{2}$.  Both of the them have lengths less
than $3.02$. Also, $|v_1|\ge2.3$.
\end{lemma}

\begin{proof}
This is a standard exercise in geometric considerations
(\cite[Section 1]{formulation} and \cite[Section 2]{part3}). (The
reader should review these sections for the framework of the
following argument.)  We deform the figure using pivots to a
configuration $v_2,\ldots,v_5$ at height $2$, and $|v_i-v_j|=2t_0
$, $(i,j)=(2,3),(3,4),(4,5),(5,2)$. We scale $v_1$ until
$|v_1|=2t_0 $. We can also take the distance from $v_1$ to $v_5$
and to $v_2$ to be $2$.  If we have $|v_1-v_3|\ge 2\sqrt{2}$, then
we stretch the edge $|v_1-v_4|$ until $|v_1-v_3|=2\sqrt{2}$. The
resulting configuration is rigid.  Pick coordinates to find that
$|v_1-v_4|<2\sqrt{2}$. If we have $|v_1-v_3|\ge 2t_0 $, follow a
similar procedure to reduce to the rigid configuration
$|v_1-v_3|=2t_0$, to find that $|v_1-v_4|<3.02$. The estimate
$|v_1|\ge2.3$ is similar.
\end{proof}

There are restrictive bounds on the dihedral angles of the
simplices $(0,v_1,v_i,v_j)$ along the edge $(0,v_1)$. The
quasi-regular tetrahedron has a dihedral angle of at most $1.875$
(\cite[10.1.2]{part3}).  The dihedral angles of the simplices
$(0,v_1,v_2,v_3)$, $(0,v_1,v_5,v_4)$ adjacent to it are at most
$1.63$. The dihedral angle of the remaining simplex
$(0,v_1,v_3,v_4)$ is at most $1.51$. (These constants come from
$1.624$ and $1.507$ in Appendix \ref{app:dih}.)  This leads to
lower bounds as well. The quasi-regular tetrahedron has a dihedral
angle that is at least $2\pi - 2(1.63)-1.51 > 1.51$.  The dihedral
angles adjacent to the quasi-regular tetrahedron is at least
$2\pi- 1.63-1.51-1.875> 1.26$. The remaining dihedral angle is at
least $2\pi-1.875-2(1.63) > 1.14$.

\begin{lemma}
\label{lemma:nobad4}
 A decomposition star that contains this
configuration is not contravening: $\tau(D) >\squander$.
\end{lemma}

\begin{proof}
We will show that the sum of $\tau_R(D)$, as $R$ ranges over the
pentagonal standard region and triangular regions, is at least
$11.16\,\pt$. Let $P$ denote the aggregate formed by these
standard regions.

First we show how the lemma follows from this bound. There are no
other exceptional faces, because $11.16\,\pt+t_5>\squander$. Every
vertex not on $P$ has type $(5,0)$, by \cite[Proposition
5.2]{part3}. In particular, there are no quadrilateral regions.
There are at most $4$ triangles at every vertex of $P$, because
the $\tauLP(5,0,2\pi-1.32)>6.02\,\pt$. (The $1.32$ comes from that
fact that $P$ has no flat quarters because $v$ is too short to be
enclosed over one \cite[1.3]{formulation} \cite[3.11.4]{part4}.)
There are at least $3$ triangles at every vertex of $P$, otherwise
we contradict a property of tameness (Property
\ref{definition:tame:3-circuit} and Property
\ref{definition:tame:4-circuit}).

The only triangulation with these properties is obtained by
removing one edge from the icosahedron (Exercise).  This implies
that there are two opposite corners of $P$ each having four
quasi-regular tetrahedra. Since the diagonals of $P$ have lengths
greater than $2\sqrt{2}$, the results of Section \ref{sec:tri4}
show that these eight quasi-regular tetrahedra squander at least
$2(1.5)\,\pt$.  There are two additional vertices of type $(5,0)$
whose tetrahedra are distinct from these eight quasi-regular
tetrahedra. They give an additional $2(0.55)\,\pt$.  Now
$(11.16+2(1.5)+2(0.55))\,\pt>\squander$.  The result follows.

Return to the bound $11.16\,\pt$. (A review of Appendix
\ref{app:part4} is recommended before continuing. The definitions
of $\tau$, $\hat\tau$ and $\tau_0$ are reviewed in Appendix
\ref{app:constants}.) Let $\dih$ denote the dihedral angle of a
simplex along a given edge. Let $S_{ij}$ be the simplex
$(0,v_1,v_i,v_j)$, for $(i,j)=(2,3),(3,4), (4,5),(2,5)$. We have
$\sum_{(4)}\dih(S_{ij}) = 2\pi$. Suppose one of the edges
$(v_1,v_3)$ or $(v_1,v_4)$ has length $\ge2\sqrt2$. Say
$(v_1,v_3)$.

We have (by Appendix \ref{app:group1})
$$
\begin{array}{lll}
\tau(S_{25}) &- 0.2529\dih(S_{25}) > -0.3442,\\
\tau_0(S_{23}) &- 0.2529\dih(S_{23}) > -0.1787,\\
\hat\tau(S_{45}) &- 0.2529\dih(S_{45}) > -0.2137,\\
\tau_0(S_{34}) &- 0.2529\dih(S_{34}) > -0.1371.\\
\end{array}
$$
We have by \cite[$\A_{10}$]{part4} that
    $$
    \begin{array}{lll}
        \tau(D) &\ge \sum_{(4)}\tau_x(S_{ij}) - 5\xiG\\
                &>2\pi(0.2529)-0.3442\\
                &\qquad -0.1787-0.2137-0.1371-5\xiG\\
                &>11.16\,\pt,
    \end{array}
    $$
where $\tau_x=\tau,\hat\tau,\tau_0$ as appropriate.

Now suppose $(v_1,v_3)$ and $(v_1,v_4)$ have length $\le2\sqrt2$.
If there is an upright diagonal that is not enclosed over either
flat quarter, the penalty is at most $3\xiG+2\xiV$. Otherwise, the
penalty is smaller: $4\xiG'+\xiV$. We have
    $$
    \begin{array}{lll}
    \tau(D)
    &\ge \sum_{(4)}\tau(S_{ij})-(3\xiG+2\xiV)\\
    &>2\pi(0.2529)-0.3442\\
    &\qquad -2(0.2137)-0.1371 -(3\xiG+2\xiV)\\
    &>11.16\,\pt.\\
    \end{array}
    $$
\end{proof}

The proof that contravening plane graphs are tame is complete.

\newpage
\section{Linear Programs} 
\label{sec:linearprogram}

We have completed the first half of the proof of the Kepler
Conjecture by proving that every contravening plane graph is tame.

The second half of the proof of the Kepler Conjecture consists in
showing that tame graphs are not contravening, except for the
isomorphism class of graphs isomorphic to $G_{fcc}$ and $G_{hcp}$
associated with the face-centered cubic and hexagonal close
packings.

The paper \cite{part3} treats all contravening tame graphs in
which every face has length three or four, except for three cases
$G_{fcc}$, $G_{pent}$, and $G_{hcp}$. The two cases $G_{fcc}$ and
$G_{hcp}$ are treated in \cite[3.13]{formulation}, and the case
$G_{pent}$ is treated in \cite{thesis}.

We complete the proof of the Kepler Conjecture by analyzing all
contravening plane graphs that contain at least one exceptional
face.  The primary tool that will be used is linear programming.
The linear programs are obtained as relaxations of the original
nonlinear optimization problem of maximizing $\sigma(D)$ over all
decomposition stars whose associated graph is a given tame graph
$G$.  The upper bounds obtained through relaxation are upper
bounds to the nonlinear problem.

To eliminate a tame graph, we must show that it is not
contravening. By definition, this means we must show that
$\sigma(D) < 8\,\pt$.  When a single linear program does not yield
an upper bound under $8\,\pt$, we branch into a sequence of linear
programs that collectively imply the upper bound of $8\,\pt$. This
will call for a sequence of increasingly complex linear programs.

\subsection{Relaxation}

(NLP) Let $f:P\to\R$ be a function on a nonempty set $P$. Consider
the nonlinear maximization problem
    $$\max_{p\in P} f(p).$$

(LP): Consider a linear programming problem
    $$\max\, c\cdot x$$ such that $A\, x\le b$, where $A$ is a matrix,
    $b$, $c$ are vectors of real constants and $x$ is a vector of
    variables $x = (x_1,\ldots,x_n)$.  We write the linear
    programming problem as
    $$\max(c\cdot x : A\,x\le b).$$

An {\it interpretation\/} $I$ of a linear programming problem (LP)
is a set $|I|$, together with an assignment $x_i\mapsto x_i^I$ of
functions $x^I_i:|I|\to\R$ to variables $x_i$.  We say the
constraints $A\,x\le b$ of the linear program are {\it
satisfied\/} under the interpretation $I$ if for all $p\in |I|$,
we have $$A\, x^I(p) \le b.$$ The interpretation $I$ is said to be
a {\it relaxation\/} of the nonlinear program (NLP), if the
following three conditions hold.
    \begin{enumerate}
    \item  $P=|I|$.
    \item The constraints are satisfied under the interpretation.
    \item $f(p)\le c\cdot x^I(p)$, for all $p\in|I|$.
    \end{enumerate}

\begin{lemma}
\label{lemma:bound} Let (LP) be a linear program with relaxation
$I$ to (NLP). Then (LP) has a feasible solution.  Moreover, if
(LP) is bounded above by a constant $M$, then $M$ is an upper
bound on the function $f:|I|\to\R$.
\end{lemma}

\begin{proof}
A feasible solution is $x_i = x_i^I(p)$, for any $p\in |I|$. The
rest is clear.
\end{proof}

\begin{remark}  In general, it is to be expected that the
interpretations $A\,x^I \le b$ will be nonlinear inequalities on
the domain $P$.  In our situation, satisfaction of the inequality
will be proved by interval arithmetic.  Thus, the construction of
an upper bound to (NLP) breaks into two tasks: to solve the linear
programs and to prove the nonlinear inequalities required by the
interpretation.
\end{remark}

There are many nonlinear inequalities that enter into our
interpretation, which have been proved by interval arithmetic on
computer.  These inequalities are contained in Appendix
\ref{app:calculation}.

\begin{remark}
\label{remark:derived} There is a second method of establishing
the satisfaction of inequalities under an interpretation. Suppose
that we wish to show that the inequality $e\cdot x\le b'$ is
satisfied under the interpretation $I$. Suppose that we have
already established that a system of inequalities $A\,x\le b$ is
satisfied under the interpretation $I$.  We solve the linear
programming problem $\max(e\cdot x : A\,x\le b)$.  If this maximum
is at most $b'$, then the inequality $e\cdot x\le b'$ is satisfied
under the interpretation $I$.  We will refer to $e\cdot x\le b'$
as an {\it LP-derived inequality} (with respect to the system
$A\,x\le b$).
\end{remark}

\subsection{The Linear Programs}

Let $G$ be a tame plane graph containing an exceptional face. Let
$X_G$ be the space of all decomposition stars whose associated
plane graph is isomorphic to $G$.

\begin{theorem} For every tame plane graph $G$ containing at least
one exceptional face, there exists a finite sequence of linear
programs with the following properties.
    \begin{enumerate}
    \item Every linear program has an admissible solution and its solution is strictly
        less than $8\,\pt$.
    \item For every linear program in this sequence,
        there is an interpretation $I$ of the linear program that is a
        relaxation of the nonlinear optimization problem
    $$\sigma:|I|\to \R,$$
        where $|I|$ is a subset of $X_G$.
    \item The union of the subsets $|I|$,
        as we run over the sequence of linear programs, is $X_G$.
    \end{enumerate}
\end{theorem}

The proof is constructive.  For every tame plane graph $G$ a
sequence of linear programs is generated by computer and solved.
The optimal solutions are all bounded above by $8\,\pt$. It will
be clear from construction of the sequence that the union of the
sets $|I|$ exhausts $X_G$.  We estimate that nearly $10^5$ linear
programs are involved in the construction.

\begin{remark} The paper \cite[Section 3.1.1]{algorithm} shows how the
linear programs that arise in connection with the Kepler
Conjecture can be formulated in such a way that they always have a
feasible solution and so that the optimal solution is bounded.  We
assume that all our linear programs have been constructed in this
way.
\end{remark}

\begin{corollary} If a tame graph $G$ contains an exceptional face,
then it is not contravening.
\end{corollary}

\begin{proof}  This follows immediately from the theorem and
Lemma \ref{lemma:bound}.
\end{proof}

\section{LP Formulation}
\label{sec:lpformulation}

\subsection{Basic Linear Programs}

To describe the basic linear program, we need the following
indexing sets.

Let \mark{VERTEX} be the set of all vertices in $G$. Let
\mark{FACE} be the set of all faces in $G$.  (Recall that by
construction each face $F$ of the graph carries an orientation.)
Let \mark{ANGLE} be the set of all angles in $G$, defined as the
set of pairs $(v,F)$, where the vertex $v$ lies in the face $F$.
Let \mark{DIRECTED} be the set of directed edges. It consists of
all ordered pairs $(v,s(v,F))$, where $s(v,F)$ denotes the
successor of the vertex $v$ in the oriented face $F$.  Let
\mark{TRIANGLES} be the subset of \mark{FACE} consisting of those
faces of length $3$.  Let \mark{UNDIRECTED} be the set of
undirected edges.  It consists of all unordered pairs
$\{v,s(v,F)\}$, for $v\in F$.

We introduce variables indexed by these sets.  Following AMPL
notation, we write for instance
    $y\{\mark{VERTEX}\}$
to declare a collection of variables $y[v]$ indexed by vertices
$v$ in \mark{VERTEX}.  With this in mind, we declare the variables
    $$\begin{array}{lllll}
        &\alpha\{\mark{ANGLE}\},\
        &y\{\mark{VERTEX}\},\
        &e\{\mark{UNDIRECTED}\},\\
        &\sigma\{\mark{FACE}\},\
        &\tau\{\mark{FACE}\},\
        &\sol\{\mark{FACE}\}.
    \end{array}
    $$

We obtain an interpretation $I$ on the compact space $X_G$. First,
we define an interpretation at the level of indexing sets.  A
decomposition star determines the set $U(D)$ of vertices of height
at most $2t_0$ from the origin of $D$.  Each decomposition star
$D\in X_G$ determines a (metric) graph with geodesic edges on the
surface of the unit sphere, which is isomorphic to $G$ as a
(combinatorial) plane graph.  There is a map from the vertices of
$G$ to $U(D)$ given by $v\mapsto v^I$, if the radial projection of
$v^I$ to the unit sphere at the origin corresponds to $v$ under
this isomorphism. Similarly, each face $F$ of $G$ corresponds to a
set $F^I$ of standard regions.   Each edge $e$ of $G$ corresponds
to a geodesic edge $e^I$ on the unit sphere.

Now we give an interpretation $I$ to the linear programming
variables at a decomposition star $D$.   As usual, we add a
superscript $I$ to a variable to indicate its interpretation. Let
$\alpha[v,F]^I$ be the sum of the internal angles at $v^I$ of the
metric graph in the standard regions $F^I$. Let $y[v]^I$ be the
length $|v^I|$ of the vertex $v^I\in U(D)$ corresponding to $v$.
Let $e[v,w]^I$ be the length $|v^I-w^I|$ of the edge between $v^I$
and $w^I\in U(D)$. Let
    $$\begin{array}{lll}
    \sigma[F]^I &= \sigma_F(D),\\
    \sol[F]^I &= \sol(F^I),\\
    \tau[F]^I &= \tau_F(D).
    \end{array}
    $$

The objective function for the optimization problems is
    $$\max:\quad \sum_{F\in\mark{FACE}} \sigma_F.$$
Its interpretation under $I$ is the score $\sigma(D)$.

We can write a number of linear inequalities  that will be
 satisfied under our interpretation.  For example, we
have the bounds
    $$
    \begin{array}{lll}
    &0\le y[v]\le 2t_0, &v\in \mark{VERTEX}\\
    &0\le e[v,w]\le 2t_0 &(v,w)\in\mark{EDGE}\\
    &0\le \alpha[v,F]\le 2\pi &(v,F)\in\mark{ANGLE}\\
    &0\le \sol[F]\le 4\pi &F\in\mark{FACE}\\
    \end{array}
    $$
There are other linear relations that are suggested directly by
the definitions or the geometry.  Here, $v$ belongs to
$\mark{VERTEX}$.
    $$
    \begin{array}{lll}
    \tau[F] &= \sol[F] \zeta\pt - \sigma[F],\\
    2\pi &= \sum_{F: v\in F} \alpha[v,F],\\
    \sol[F] &= \sum_{v\in F}\alpha[v,F] - (len(F)-2)\pi\\
    \end{array}
    $$
There are long lists of additional inequalities that come from
interval arithmetic verifications. For example, the interval
verifications of \cite[Propositions 4.1-4.3]{formulation} supply
us with many linear inequalities relating the variables
    $$
    \begin{array}{lll}
    \sigma[F],&\tau[F],&\alpha[v,F],\\
    \sol[F],&y[v],&e[v,w]
    \end{array}
    $$
whenever $F^I$ is a single standard region having three sides.
Similarly, \cite[Section 4]{part3} gives inequalities for
$\sigma[F]$ and related variables, when the length of $F$ is four.
A complete list of inequalities that are used for triangular and
quadrilateral faces is found in \cite[Section 9]{part3}.

For exceptional faces, we have an admissible weight function
$w(F)$.  According to definitions $w(F) = \tau[F]/pt$, so that the
inequalities for the weight function can be expressed in terms of
the linear program variables.

When the exceptional face is not an aggregate, then it also
satisfies the inequalities of \cite[Theorem 4.4]{part4}.

\subsection{Face Refinements}
\label{sec:facerefinement}

The linear program formulated in the previous section does not
carry enough information in general to obtain the desired bounds.
The paper \cite{part4} shows that exceptional standard regions
have a great deal of internal structure. Until now in this paper,
we have mostly ignored this internal structure. We will now make a
detailed examination of this internal structure.

A {\it refinement\/} $\tildeF$ of a face $F$ of a plane graph $G$
is a set $\tildeF$ of faces  such that
    \begin{enumerate}
    \item The intersection of
the vertex set of $G$ with that of $\tildeF$ is the set $F$.
    \item $\tildeF \cup\{F^{op}\}$ is a plane graph.
    \end{enumerate}
We use refinements of faces to describe the internal structure of
faces.

We introduce indexing sets $\mark{FACE-\tildeF }$,
$\mark{VERTEX-\tildeF }$, $\mark{ANGLE-\tildeF }$,
$\mark{EDGE-\tildeF }$, the sets of faces, vertices, angles, and
edges in $\tildeF$, respectively, analogous to those introduced
for $G$.

We create variables
    $\pi[\tildeF]$, and indexed variables
    $$
    \begin{array}{lll}
    \sol\{\mark{FACE-\tildeF }\}, & \sc\{\mark{FACE-\tildeF }\} &
    \tausc\{\mark{FACE-\tildeF }\},\\
    \alpha\{\mark{ANGLE-\tildeF }\}, & y\{\mark{VERTEX-\tildeF }\}, &
    e\{\mark{EDGE-\tildeF }\}.
    \end{array}
    $$
(Variables with names ``$y[v]$'' and ``$e[v,w]$'' were already
created for some $v,w\in \mark{VERTEX-\tildeF }\cap
\mark{VERTEX}$.  In these cases, we use the variables already
created.)

Each vertex $v$ in the refinement will be interpreted either as a
vertex $v^I\in U(D)$, or as the endpoint of an upright diagonal
lying over the standard region $F^I$. We will interpret the faces
of the refinement in terms of the geometry of the decomposition
star $D$ variously as flat quarters, upright quarters, anchored
simplices, and the other constructs of \cite{part4}. This
interpretation depends on the context, and will be described in
greater detail below.

Once the interpretation of faces is fixed, the interpretations are
as before for the variable names introduced already: $y$, $e$,
$\alpha$, $\sol$.  The lower and upper bounds for $\alpha$ and
$\sol$ are as before.  The lower and upper bounds for $y[v]$ are
$2$ and $2t_0$ if $v^I\in U(D)$, but if $(0,v^I)$ is an upright
diagonal, then the bounds are $[2t_0,2\sqrt2]$. The lower and
upper bounds for $e$ will depend on the context.

\subsection{Variables related to score}
\label{sec:variable}

The variables $\sc$ are a stand-in for the score $\sigma$ on a
face. We do not call them $\sigma$ because the sum of these
variables will not in general equal the variable $\sigma[F]$, when
$\tildeF$ is a refinement of $F$:
    $$[\sum_{F'\in \tildeF }\sc[F'] \ne \sigma[F]].$$
We will use have a weaker relation:
    $$\sigma[F] \le \sum_{F'\in \tildeF }\sc[F'] + \pi[\tildeF].$$
The variable $\pi[\tildeF]$ is called the penalty associated with
the refinement $\tildeF$.  (Penalties are discussed at length in
Appendix \ref{app:part4}.)   The interpretations of $\sc$ and
$\pi[\tildeF]$ are rather involved, and will be discussed on a
case-by-case basis below. The interpretation of $\tausc$ follows
from the identity:
    $$\tausc[F'] = \sol[F']\zeta\pt - \sc[F'],\quad\forall
    F'\in \tildeF .$$

The interpretation of variables that follows might appear to be
hodge-podge at first.  However, they are obtained in a systematic
way.  We analyze the proofs and approximations in \cite{part4},
and define $\sc[F]^I$ as the best penalty-free scoring
approximation that is consistent with the given face refinement.
here are the details.

If the subregion is a flat quarter, the interpretation of $\sc[F]$
is the function $\hat\sigma$, defined in Appendix \ref{app:part4}.
If the subregion is an upright quarter $Q$, the interpretation of
$\sc[F]$ is the function $\sigma(Q)$ from \cite{formulation}. If
the subregion is an anchored simplex that is not an upright
quarter, $\sc[F]$ is interpreted as the analytic Voronoi function
$\vor$ if the simplex has type $C$ or $C'$, and as $\vor_0$
otherwise. (The types $A$, $B$, $C$ and $C'$ are defined in
\cite[Sections~2.5--2.10]{part4}.) Whether or not the simplex has
type $C$, the inequality $\sc[F]\le0$ is satisfied. In fact, if
$\vor_0$ scoring is used, we note that there are no quoins, and
$\phi(1,t_0)<0$.

If the subregion is triangular, if no vertex represents an upright
diagonal, and if the subregion is not a quarter, then $\sc[F]$ is
interpreted as  $\vor$ or $\vor_0$ depending on whether the
simplex has type $A$.  In either case, the inequality
$\sc[F]\le\vor_0$ is satisfied.

In most other cases, the interpretation of $\sc[F]$ is $\vor_0$.
However, if $R$ is a heptagon or octagon, and $F$ has $\ge4$
sides, then $\sc[F]$ is interpreted as $\vor_0$ except on
simplices of type $A$, where it becomes the analytic Voronoi
function.

If $R$ is a pentagon or hexagon, and $F$ is a quadrilateral that
is not adjacent to a flat quarter, then the interpretation of
$\sc[F]$ is the actual score of the subregion over the subregion.
In this case, the score $\sigma_R$ has a well-defined meaning for
the quadrilateral, because it is not possible for an upright
quarter in the $Q$-system to straddle the quadrilateral region and
an adjacent region. Consequently, any erasing that is done can be
associated with the subregion without ambiguity. By the results of
\cite{part2}, we have $\sc[F]\le0$. We also have
$\sc[F]\le\vor_0$.

One other bound that we have not explicitly mentioned is the bound
$\sigma_R(D)< s_n$.  For heptagons and octagons that are not
aggregates, this is a better bound than the one used in the
definition of tameness (Property \ref{definition:tame:score}). In
heptagons and octagons that are not aggregates, if we have a
subregion with four or more sides, then $\sc[F]< Z(n,k)$ and
$\tausc[F]>D(n,k)$. (See \cite[5.5.1]{part4},
\cite[5.5.2]{part4}.)

The variables are subject to a number of compatibility relations
that are evident from the underlying definitions and geometry.
    $$
    \begin{array}{lll}
    \sol[F'] = \sum_{v\in F'}\alpha[v,F'] - (len[F']-2)\pi,&\forall F'\\
    \sum_{F':v\in F',F'\in \mark{FACE-\tildeF }} \alpha[v,F'] =
    \alpha[v,F], &\forall v
    \end{array}
    $$

Assume that a face $F_1\in \tildeF$ has been interpreted as a
subregion $R=F_1^I$ of a standard region.  Assume that each vertex
of $F_1$ is interpreted as a vertex in $U(D)$ or as the endpoint
of an upright diagonal over $F^I$.  One common interpretation of
$\sc$ is $\vor_{0,F}(U(D))$, the truncated Voronoi function. When
this is the interpretation, we introduce further variables:
    $$
    \begin{array}{lll}
    \quo[v,s(v,F_1)] & \forall v\in F_1,\\
    \quo[s(v,F_1),v] & \forall v\in F_1,\\
    \Adih[v,F_1] & \forall v\in F_1,\\
    \end{array}
    $$
We interpret the variables as follows. If $w = s(v)$, and the
triangle $(0,v^I,w^I)$ has circumradius $\eta$ at most $t_0$, then
    $$
    \begin{array}{lll}
    \quo[v,w]^I &= \quo(R(|v^I|/2,\eta,t_0)),\\
    \quo[w,v]^I &= \quo(R(|w^I|/2,\eta,t_0)).
    \end{array}
    $$
If the circumradius is greater than $t_0$, we take
    $$\quo[v,w]^I =\quo[w,v]^I =0.$$
The variable $\Adih$ has the following interpretation:
    $$\Adih[v,F_1]^I = \begin{cases}
        A(|v^I|/2)\alpha(v^I,F_1^I) & |v^I|\le 2t_0,\\
        0 & \text{otherwise.}\\
        \end{cases}
    $$
Under these interpretations, the following identity is satisfied:
    $$
    \begin{array}{lll}
    \sc[F_1] &= \sol[F_1]\phi_0 +\sum_{v\in F_1} \Adih[v,F_1]\\
            &\quad - 4\doct \sum \quo[v,w].
    \end{array}
    $$
This relation and the constants that appear in it are based on
\cite[Appendix A.4]{part3}.  The final sum runs over all pairs
$(v,w)$, where $v=s(w,F_1)$ or $w=s(v,F_1)$.

For this to be useful, we need good inequalities governing the
individual variables.  Such inequalities for $\Adih[v,F]$ and
$\quo[v,w]$ are found in \cite[A.4]{part3}.  To make of these
inequalities, it is necessary to have lower and upper bounds on
$\alpha[v,F]$ and $y[v]$.  We obtain such bounds as LP-derived
inequalities in the sense of Remark \ref{remark:derived}.

\subsection{Triangle and Quad Branching}
\label{sec:tribranch}

The paper \cite[Appendix A.7]{part3} divides triangular faces into
two cases:
    $$
    \begin{array}{lll}
        y[v_1]+y[v_2]+y[v_3] &\le 6.25.\\
        y[v_1]+y[v_2]+y[v_3] &\ge 6.25.\\
    \end{array}
    $$
whenever sufficiently good bounds are not obtained as a single
linear program.  That paper also divides quadrilateral faces into
four cases: two flat quarters, two flat quarters with diagonal
running in the other direction, four upright quarters forming an
octahedron, no quarters.  In general, if there are $r_1$ triangles
and $r_2$ quadrilaterals, we obtain as many as $2^{r_1 + 2 r_2}$
cases by breaking the various triangles and quadrilaterals into
subcases.

We break triangular faces and quadrilaterals into subcases as
needed in the linear programs that follow without further comment.

\section{Elimination of Aggregates}
\label{sec:noaggregate}

The proof of the following theorem occupies the entire section. It
eliminates all the pathological cases that we have had to carry
along until now.  The terminology in the next lemma comes from
\cite{part4} (anchors, $\Sminus$, $\Splus$, and so forth).  This
terminology is reviewed in Appendix \ref{app:part4}.

\begin{theorem}
\label{theorem:noaggregate}
 Let $D$ be a contravening
decomposition star, and let $G$ be its tame graph.  Every face of
$G$ corresponds to exactly one standard region of $D$.  No
standard region of $D$ has any enclosed vertices from $U(D)$.
(That is, enclosed vertices have height at least $2t_0$.  The
decomposition star $D$ does not contain any $\Sminus$ or $\Sfour$
configurations.  If an upright diagonal with at least $5$ anchors
appears, the $5$ anchors are the five corners of $U(D)$ lying over
a pentagonal standard region.
\end{theorem}

\subsection{ A pentagonal hull with $n=8$ } 
\label{sec:pentagonal}
 \label{sec:3.6}

The next few sections treat the nonpolygonal standard regions
illustrated in Remark \ref{remark:tri-pent}. In this subsection,
there is an aggregate of the octagonal region and a triangle has a
pentagonal hull. Let $P$ denote this aggregate.

\begin{lemma}
\label{lemma:not301} Let $G$ be a  contravening plane graph with
this aggregate. Some vertex on the pentagonal face has type not
equal to $(3,0,1)$.
\end{lemma}

\begin{proof} If every vertex on the pentagonal face has type
$(3,0,1)$, then at the vertex of the pentagon meeting the
aggregated triangle, the four triangles together with the octagon
give
    $$t_8+\sum_{(4)}\tauLP(4,0,2\pi-2(1.153)) > \squander,$$
so that the graph does not contravene.
\end{proof}

For a general contravening plane graph with this aggregate,
 we have bounds
    $$
    \begin{array}{lll}
    \sigma_F(D)&\le\pt+s_8,\\
    \tau_F(D)&\ge t_8.
    \end{array}
    $$
We add the inequalities $\tausc[F]>t_8$ and $\sc[F]< \pt+s_8$ to
the exceptional face. There is no other exceptional face, because
$t_8+t_5>\squander$. We run the linear programs for all tame
graphs with the property asserted by Lemma \ref{lemma:not301}.
Every upper bound is less than $8\,\pt$, so that there are no
contravening decomposition stars with this configuration.

\subsection{ $n=8$, hexagonal hull} 
\label{sec:hexagonal}
 \label{3.7}

We treat the two cases from Remark \ref{remark:tri-pent} that have
a hexagonal hull  (Figure \ref{fig:aggregates}). One can be
described as a hexagonal region with an enclosed vertex that has
height at most $2t_0$ and distance at least $2t_0$ from each
corner over the hexagon.  The other is described as a hexagonal
region with an enclosed vertex of height at most $2t_0$, but this
time with distance less than $2t_0$ from one of the corners over
the hexagon.

The argument for the case $n=8$ with hexagonal hull is similar to
the argument of Section \ref{sec:pentagonal}. Add the inequalities
$\tau(R)>t_8$ and $\sigma(R)<s_8$ for each hexagonal region. Run
the linear programs for all tame graphs, and check that these
additional inequalities yield linear programming bounds under
$8\,\pt$.

\subsection{$n=7$, pentagonal hull} 
\label{sec:3.8}

We treat the two cases illustrated in Figure \ref{fig:aggregates}
that have a pentagonal hull. One can be described as a pentagon
with an enclosed vertex that has height at most $2t_0$ and
distance at least $2t_0$ from each corner of the pentagon. The
other is described as a pentagon with an enclosed vertex of height
at most $2t_0$, but this time with distance less than $2t_0$ from
one of the corners of the pentagon.

In discussing various maps, we let $v_i$ be the corners of the
regions, and we set $y_i = |v_i|$ and $y_{ij}= |v_i-v_j|$. The
subscript $F$ is dropped, when there is no great danger of
ambiguity.

Add the inequalities $\tau[F]>t_7$, $\sigma[F]<s_7$ for the
pentagonal face. There is no other exceptional region, because
$t_5+t_7> \squander$. With these changes, of all the tame plane
graphs with a pentagonal face, all but one of the linear programs
give a bound under $8\,\pt$.

The plane graph $G_{15}$ that remains is easy to describe.  It is
the plane graph with $11$ vertices, obtained by removing from an
icosahedron a vertex and all the edges that meet at that vertex.

We treat the case $G_{15}$.  Familiarity with Appendix
\ref{app:part4} will be assumed in the following passage.   Let
$v_{12}$ be the vertex enclosed over the pentagon. We let
$v_1,\ldots,v_5$ be the five corners of $U(D)$ over the pentagon.
Break the pentagon into five simplices along $(0,v_{12})$:  $S_i =
(0,v_{12},v_i,v_{i+1})$. We have LP-derived bounds (in the sense
of Remark \ref{remark:derived}) $y[v_i]\le2.168$, and
$\alpha[v_i,F]\le2.89$, for $i=1,2,3,4,5$. In particular, the
pentagonal region is convex.

Another LP-derived inequality is
    $$\sigma[F] > -0.2345.$$
Since $-0.4339$ is less than this the lower bound, the
configuration $\Sminus$ does  not occur. Similarly, since $-0.25$
is less than the lower bound, $\Sfour$ does not occur
(\cite[3.7~and~3.8]{part4}.

Suppose that there is a loop in context $(n,k)=(4,2)$. Again the
score is less than the LP-lower bound, showing that this does not
occur:
$$\sigma_R(D) < s_7+ \zloop(4,2) < -0.2345.$$ The constants come from Table
\cite[Table 5.11]{part4} and \cite[Theorem 4.4]{part4}.

We have the LP-derived inequality
$$\tau[F] < 0.644.$$
If we branch and bound on the triangular faces, this LP-derived
inequality can be improved to
    $$\tau[F] < 0.6079.$$
If there is a loop other than $(4,2)$ and $(4,1)$, the linear
program becomes infeasible:
    $$\tau[F] < 0.644 < t_7 + \dloop(n,k) < \tau[F].$$
We conclude that all loops have context $(n,k)=(4,1)$.

Case 1.  {\it The vertex $v_{12}$ has distance at least $2t_0$
from the five corners of $U(D)$ over the pentagon.}

The interval calculations relevant in Case 1 appear in
Appendix~\ref{A.3.8}.

The penalty to switch the pentagon to a pure $\vor_0$ score is at
most $5\xiG$ (see \cite[$\A_{10}$]{part4}). There cannot be two
flat quarters because then
$$|v_{12}|>\E(S(2,2,2,2t_0,2\sqrt2,2\sqrt2),2t_0,2t_0,2t_0)>2t_0.$$

(Case 1-a)
Suppose there is one flat quarter, $|v_1-v_4|\le2\sqrt2$.
There is a lower bound of 1.2 on the dihedral angles of the
simplices $(0,v_{12},v_i,v_{i+1})$.  This is
obtained as follows.  The proof relies on the convexity of the
quadrilateral region.  We leave it to the reader to verify that
the following pivots can be made to preserve convexity.  Disregard
all vertices except $v_1,v_2,v_3,v_4,v_{12}$.  We give the argument
that $\dih(0,v_{12},v_1,v_4)>1.2$.  The others are similar.
Disregard the length $|v_1-v_4|$.  We show that
    $$
    \begin{array}{lll}
        sd &:=\dih(0,v_{12},v_1,v_2)+\dih(0,v_{12},v_2,v_3)\\
           &+\dih(0,v_{12},v_3,v_4) < 2\pi-1.2.
    \end{array}
    $$
Lift $v_{12}$ so $|v_{12}|=2t_0$. Maximize $sd$ by taking
$|v_1-v_2|=|v_2-v_3|=|v_3-v_4|=2t_0$.  Fixing $v_3$ and $v_4$,
pivot $v_1$ around $(0,v_{12})$ toward $v_4$, dragging $v_2$
toward $v_{12}$ until $|v_2-v_{12}|=2t_0$.  Similarly, we obtain
$|v_3-v_{12}|=2t_0$. We now have $sd\le 3(1.63)< 2\pi-1.2$, by an
interval calculation.

Return to the original figure and move $v_{12}$ without increasing
$|v_{12}|$ until each simplex $(0,v_{12},v_i,v_{i+1})$ has an edge
$(v_{12},v_j)$ of length $2t_0$.
  Interval calculations show that
the four simplices around $v_{12}$ squander
    $$2\pi(0.2529)-3(0.1376)-0.12 > \squander + 5\xiG.$$

(Case 1-b) Assume there are no flat quarters.  An LP-derived bound
on the perimeter $\sum|v_i-v_{i+1}|$ is 1.407. We have
$\arc(2,2,x)'' = 2x/(16-x^2)^{3/2} >0$. The arclength of the
perimeter is therefore at most
$$2\arc(2,2,2t_0) + 2\arc(2,2,2) + \arc(2,2,2.387) <  2\pi.$$
There is a well-defined interior of the spherical pentagon, a
component of area $<2\pi$.  If we deform by decreasing the
perimeter, the component of area $<2\pi$ does not get swapped with
the other component.

Disregard all vertices but $v_1,\ldots,v_5,v_{12}$.  If a vertex
$v_i$ satisfies  $|v_i-v_{12}|>2t_0$, deform $v_i$ as in
\cite[Section 4.9]{part4} until $|v_{i-1}-v_{i}|=|v_i-v_{i+1}|=2$,
or $|v_i-v_{12}|=2t_0$. If at any time, four of the edges realize
the bound $|v_i-v_{i+1}|=2$, we have reached an impossible
situation, because it leads to the contradiction (see Appendix
\ref {A.3.8})
    $$2\pi = \sum^{(5)}\dih < 1.51 + 4 (1.16).$$
(This inequality relies on the observation, which we leave to the
reader, that in any such assembly, pivots can by applied to bring
$|v_{12}-v_i|=2t_0$ for at least one edge of each of the five
simplices.)

The vertex $v_{12}$ may be moved without increasing $|v_{12}|$ so
that eventually by these deformations (and reindexing if
necessary) we have $|v_{12}-v_i|=2t_0$, $i=1,3,4$. (If we have
$i=1,2,3$, the two dihedral angles along $(0,v_2)$ satisfy
$<2(1.51)<\pi$, so the deformations can continue.)

There are two cases. In both cases $|v_i-v_{12}|=2t_0$, for
$i=1,3,4$.
$$
\begin{array}{lll}
(i)\quad &|v_{12}-v_2|=|v_{12}-v_5|=2t_0,\\
(ii)\quad &|v_{12}-v_2|=2t_0,\quad |v_4-v_5|=|v_5-v_1|=2,\\
\end{array}
$$
Case (i) follows from interval calculations
$$
\sum\tau_0 \ge 2\pi(0.2529) - 5 (0.1453) > 0.644+7\xiG.
$$
In case (ii), we have again
    $$2\pi(0.2529)-5 (0.1453).$$
In this interval calculation we have assumed that $|v_{12}-v_5|<3.488$.
Otherwise, setting $S=(v_{12},v_4,v_5,v_1)$, we have
    $$\Delta(S) < \Delta(3.488^2,4,4,8,(2t_0)^2,(2t_0)^2)<0,$$
and the simplex does not exist.  ($|v_4-v_1|\ge2\sqrt2$ because
there are no flat quarters.) This completes Case 1.

\medskip

Case 2. {\it The vertex $v_{12}$ has distance at most $2t_0$ from
the vertex $v_1$ and distance at least $2t_0$ from the others.}

Let $(0,v_{13})$ be the upright diagonal of a loop $(4,1)$.  The
vertices of the loop are not $v_2,v_3,v_4,v_5$ with $v_{12}$
enclosed over $(0,v_2,v_5,v_{13})$ by \cite[Lemma 3.6]{part4}.  The
vertices of the loop are not $v_2,v_3,v_4,v_5$ with $v_{12}$
enclosed over $(0,v_1,v_2,v_5)$ because this would lead to a
contradiction
$$y_{12}\ge \E(S(2,2,2,2t_0,2t_0,3.2),2t_0,2t_0,2)>2t_0,$$
or
$$y_{12}\ge \E(S(2,2,2,2t_0,2t_0,3.2),2,2t_0,2)>2t_0.$$
We get a contradiction for the same reasons
 unless $(v_1,v_{12})$ is an edge of some
upright quarter of every loop of type $(4,1)$.

We consider two cases.  (2-a) There is a flat quarter along an
edge other than $(v_1,v_{12})$.  That is, the central vertex is
$v_2$, $v_3$, $v_4$, or $v_5$.  (Recall that the {\it central
vertex} of a flat quarter is the vertex other than the origin that
is not an endpoint of the diagonal.) (2-b) Every flat quarter has
central vertex $v_1$.

Case 2-a.  We erase all upright quarters including those in loops,
taking penalties as required.
There cannot be two flat quarters by geometric considerations
$$
\begin{array}{lll}
\E(S(2,2,2,2\sqrt2,2\sqrt2,2t_0),2t_0,2t_0,2)&>2t_0\\
\E(S(2,2,2,2\sqrt2,2\sqrt2,2t_0),2,2t_0,2t_0)&>2t_0\\
\end{array}
$$

The penalty is at most $7\xiG$.  We show that the region (with
upright quarters erased) squanders $>7\xiG+0.644$.  We assume that
the central vertex is $v_2$ (case 2-a-i) or $v_3$ (case 2-a-ii).
In case 2-a-i, we have three types of simplices around $v_{12}$,
characterized by the bounds on their edge lengths.  Let
$(0,v_{12},v_1,v_5)$ have type A, $(0,v_{12},v_5,v_4)$ and
$(0,v_{12},v_4,v_3)$ have type B, and let $(0,v_{12},v_3,v_1)$
have type C.  In case 2-a-ii there are also three types.  Let
$(0,v_{12},v_1,v_2)$ and $(0,v_{12},v_1,v_5)$ have type A,
$(0,v_{12},v_5,v_4)$ type B, and $(0,v_{12},v_2,v_4)$ type D.
(There is no relation here between these types and the types of
simplices $A$, $B$, $C$ defined in Appendix \ref{sec:flat} and
\ref{sec:typeAB}.) Upper bounds on the dihedral angles along the
edge $(0,v_{12})$ are given in Appendix \ref {A.3.8}. These upper
bounds come as a result of a pivot argument similar to that
establishing the bound 1.2 in Case 1-a.

These upper bounds imply the following lower bounds.  In case 2-a-i,
$$
\begin{array}{lll}
\dih &> 1.33 \quad(A),\\
\dih &> 1.21 \quad(B),\\
\dih &> 1.63 \quad(C),\\
\end{array}
$$
and in case 2-a-ii,
$$
\begin{array}{lll}
\dih &> 1.37 \quad(A),\\
\dih &> 1.25 \quad(B),\\
\dih &> 1.51 \quad(D),\\
\end{array}
$$
In every case the dihedral angle is at least $1.21$. In case
2-a-i, the interval calculations of Appendix \ref{A.3.8} lead to a
lower bound on what is squandered by the four simplices around
$(0,v_{12})$. Again, we move $v_{12}$ without decreasing the score
until each simplex $(0,v_{12},v_i,v_{i+1})$ has an edge satisfying
$|v_{12}-v_j|\le2t_0$.  Interval calculations give
    $$
    \begin{array}{lll}
    \sum_{(4)}\tau_0 &> 2\pi (0.2529) - 0.2391-2(0.1376)-0.266\\
        &>0.808.
    \end{array}
    $$
In case 2-a-ii, we have
    $$
    \begin{array}{lll}
    \sum_{(4)}\tau_0 &> 2\pi (0.2529) - 2(0.2391)-0.1376-0.12\\
        &>0.853.
    \end{array}
    $$
So we squander more than $7\xiG+0.644$, as claimed.

Case 2-b.  We now assume that there are no flat quarters with central
vertex $v_2,\ldots,v_5$.
We claim
 that $v_{12}$ is not enclosed over $(0,v_1,v_2,v_3)$ or
$(0,v_1,v_5,v_4)$.
In fact, if $v_{12}$ is enclosed over $(0,v_1,v_2,v_3)$, then
we reach the contradiction
    $$
    \begin{array}{lll}
    \pi&<\dih(0,v_{12},v_1,v_2)+\dih(0,v_{12},v_2,v_3)\\
        &< 1.63+1.51.
    \end{array}
    $$

We claim
 that $v_{12}$ is not enclosed over $(0,v_5,v_1,v_2)$.
Let $S_1=(0,v_{12},v_1,v_2)$, and $S_2=(0,v_{12},v_1,v_5)$.  We
have the linear programming bound
$$y_4(S_1)+y_4(S_2) = |v_1-v_2|+|v_1-v_5|< 4.804.$$
An interval calculation gives (Appendix \ref{A.3.8bis})
    $$
    \begin{array}{lll}
    \sum_{(2)}\dih(S_i) &\le \sum_{(2)}
    \left(\dih(S_i)+0.5(0.4804/2-y_4(S_i))\right)\\
    &<\pi.
    \end{array}
    $$
So $v_{12}$ is not enclosed over $(0,v_1,v_2,v_5)$.

Erase all upright quarters, taking penalties as required.  Replace
all flat quarters with $\vor_0$-scoring taking penalties as
required. (Any flat quarter has $v_1$ as its central vertex.) We
move $v_{12}$ keeping $|v_{12}|$ fixed and not decreasing
$|v_{12}-v_1|$.  The only effect this has on the score comes
through the quoins along $(0,v_1,v_{12})$.  Stretching
$|v_{12}-v_1|$ shrinks the quoins and increases the score.  (The
sign of the derivative of the quoin with respect to the top edge
is computed in the proof of \cite[Lemma 4.9]{part4}.)

If we stretch $|v_{12}-v_1|$ to length $2t_0$, we are done by case
1 and case 2-a. (If deformations produce a flat quarter, use case
2-a, otherwise use case 1.) By the claims, we can eventually
arrange (reindexing if necessary) so that
$$
\begin{array}{lll}
(i)&\quad |v_{12}-v_3|=|v_{12}-v_4|=2t_0,\quad\text{or}\\
(ii)&\quad |v_{12}-v_3|=|v_{12}-v_5|=2t_0.
\end{array}
$$
We combine this with the deformations of \cite[Section 4]{part4}
so that in case (i) we may also assume that if
$|v_5-v_{12}|>2t_0$, then $|v_4-v_5|=|v_5-v_1|=2$ and that if
$|v_2-v_{12}|>2t_0$, then $|v_1-v_2|=|v_2-v_3|=2$. In case (ii) we
may also assume that if $|v_4-v_{12}|>2t_0$, then
$|v_3-v_4|=|v_4-v_5|=2$ and that if $|v_2-v_{12}|>2t_0$, then
$|v_1-v_2|=|v_2-v_3|=2$.

Break the pentagon into subregions by cutting along the edges
$(v_{12},v_i)$ that satisfy $|v_{12}-v_i|\le2t_0$. So for example
in case (i), we cut along $(v_{12},v_3)$, $(v_{12},v_4)$,
$(v_{12},v_1)$, and possibly along $(v_{12},v_2)$ and
$(v_{12},v_5)$.  This breaks the pentagon into triangular and
quadrilateral regions.

In case (ii), if $|v_4-v_{12}|>2t_0$, then the argument used in
Case 1 to show that $|v_4-v_{12}|<3.488$ applies here as well. In
case (i) or (ii), if $|v_{12}-v_2|>2t_0$, then for similar
reasons, we may assume
    $$\Delta(|v_{12}-v_2|^2,4,4,8,(2t_0)^2,|v_{12}-v_1|^2)\ge0.$$
This justifies the hypotheses for the interval calculations in
Appendix \ref{A.3.8}.  We conclude that
$$\sum\tau_0 \ge 2\pi (0.2529) - 3 (0.1453) - 2 (0.2391) > 0.6749.$$
Recall that the LP-upper bound on what is squandered is $0.6079$.
If the penalty is less than $0.067=0.6749-0.6079$, we are done.

We have ruled out the existence of all loops except $(4,1)$. Note
that a flat quarter with central vertex $v_1$ gives penalty at
most $0.02$ by \cite[Lemma 3.11.3]{part4}.
  If there is at most one
such a flat quarter and at most one loop,
we are done:
$$3\xiG + 0.02 < 0.067.$$
Assume there are two loops of context $(n,k)=(4,1)$.  They both
lie along the edge $(v_1,v_{12})$, which precludes any unmasked
flat quarters. If one of the upright diagonals has height
$\ge2.696$, then the penalty is at most $3\xiG+3\xiV< 0.067$.
Assume both heights are at most $2.696$. The total internal angle
of the exceptional face at $v_1$ is at least four times the
dihedral angle of one of the flat quarters along $(0,v_1)$, or
$4(0.74)$ by an interval calculation. This is contrary to the
linear programming upper bound $(0.289)$ on the dihedral angle
along $(0,v_1)$.  This completes Case 2. This shows that heptagons
with pentagonal hulls do not occur.

\subsection{$\Sminus$, $\Sfour$} 
\label{sec:3.9}

We use the bounds on the score and on what is squandered from
\cite[Sections~3.7~and~3.8]{part4}: $\Sfour$ ($\sigma<-0.25$,
$\tau>-0.4$), $\Sminus$ ($\sigma<-0.4339$, $\tau>0.5606$). If we
have the configuration $\Sminus$, there is only one exceptional
region ($0.5606+t_5>\squander$). The configuration with six
standard regions around a vertex from Section \ref{sec:stargraph}
does not occur because a calculation with $\tauLP$ (see Section
\ref{sec:2.2}) show that the five quasi-regular tetrahedra
 in the configuration squander $>6\,\pt$,
giving $6\,\pt +0.5606>\squander$.

We add the inequalities $\tau>0.5606$ and $\sigma< -0.4339$ to the
exceptional region. All linear programming bounds drop under
$8\,\pt$ when these changes are made.

The configuration $\Sfour$ requires more work.  We begin with a
lemma.

\begin{lemma}
Let $\alpha$ be the dihedral angle along the large gap in an
$\Sfour$ configuration.  Let $v_1$ and $v_2$ be the anchors of
$U(D)$ lying along the large gap.  If $|v_1|+|v_2| < 4.6$, then
$\alpha >1.78$ and the score of the four quarters is at most
$-0.31547$.
\end{lemma}

\begin{proof}
The bound $\alpha>1.78$ is an interval calculation (Appendix
\ref{app:178}). The upper bound on the score is a linear
programming calculation involving the inequality $\alpha > 1.78$
and the known inequalities on the score of an upright quarter.
\end{proof}

 Add the inequalities $\sigma<-0.25$, $\tau>0.4$ at the
exceptional regions. The configuration $\Sfour$ does not appear in
a pentagon. Run the linear programs for all tame plane graphs with
an exceptional region that is not a pentagon. If this linear
program fails to produce a bound of $8\,\pt$, we use the lemma to
branch into two cases: either $y[v_1]+y[v_2]\ge4.6$, or
$\sigma[R]<-0.31547$. In every case the bound drops below
$8\,\pt$.

\subsection{ Type $(n,k)=(5,1)$} 
\label{sec:3.10}

We return briefly to the case of six standard regions around a
vertex discussed in Section \ref{sec:stargraph}.  In the plane
graph they are aggregated into an octahedron.  We take each of the
remaining cases with an octahedron, and replace the octahedron
with a pentagon and six triangles around a new vertex.  There are
eight ways of doing this. All eight ways in each of the cases
gives an LP bound under $8\,\pt$.  This completes this case.

\subsection{ Four anchors}

In the context $(n,k)=(4,3)$, the standard region $R$ must have at
least seven sides $n(R)\ge7$.   Then
    $$
    \begin{array}{lll}
    \tau(D)&\ge t_7+\dloop(4,3)\\
            &>\squander.
    \end{array}
    $$
Thus, we may assume that this context does not occur.

\subsection {Five Anchors}

Now turn to the decomposition stars with an upright diagonal with
five anchors. Five quarters around a common upright diagonal in a
pentagonal region can certainly occur.  We claim that any other
upright diagonal with five anchors leads to a decomposition star
that does not contravene.  In fact, the only other possible
context is $(n,k)=(5,1)$ (see \cite[5.11]{part4}).

If this appears in an octagon, we have
$$\tau(D) >\dloop(5,1)+t_8>\squander.$$
 If this appears in a heptagon, we
have
$$\tau(D) >\dloop(5,1) + t_7+ 0.55\,\pt > \squander,$$
because there must be a vertex that is not a corner of the
heptagon. It cannot appear on a pentagon.  If it appears on a
hexagonal region $R$, we have
$$\tau(D) > \dloop(5,1)+t_6$$
 and
$\sigma(D) < \zloop(5,1) + s_6$. Add these constraints to the
linear program of the tame graphs with a hexagonal face.
  The LP-bound on $\sigma(D)$
with these changes is $<8\,\pt$.

\section{Branching and Bound} 
\label{sec:branchbound}

In Section \ref{sec:tribranch} we mentioned that when the linear
programs do not give sufficiently good bounds, we branch into
various subcases depending on conditions triangles and
quadrilaterals in the graph.

This section describes additional branching for exceptional
regions.  We assume the results from the Section
\ref{sec:noaggregate} that eliminate the most unpleasant types of
configurations.

\subsection{Pent and Hex Branching}

In Section \ref{sec:facerefinement} we discussed face refinements
in abstract terms.  We now describe these refinements in detail.

The possibilities are listed in the diagram only up to symmetry by
the dihedral group action on the polygon.  We do not prove the
completeness of the list, but its completeness can be seen by
inspection, in view of the comments that follow here and in
Section \ref{sec:4.2}.


\begin{figure}[htb]
  \centering
  \includegraphics{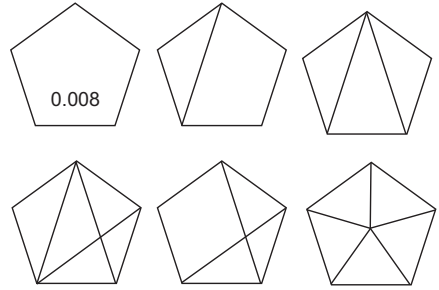}
  \caption{Pentagonal Face Refinements}
  \label{fig:penthex}
\end{figure}

\begin{figure}[htb]
  \centering
  \includegraphics{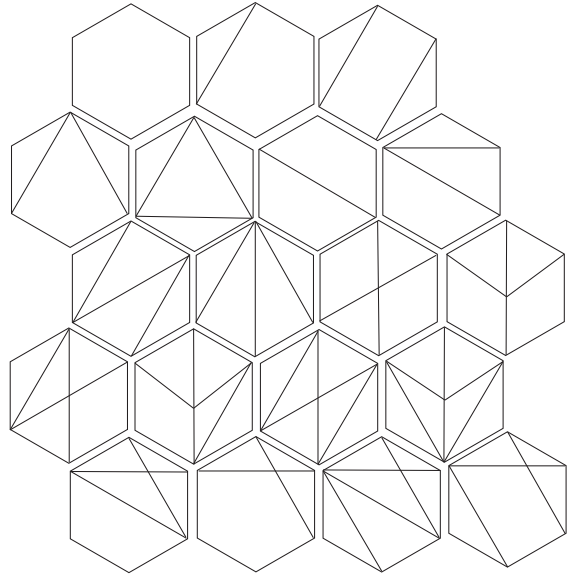}
  \caption{Hexagonal Face Refinements.  The only figures with a penalty are the
  first two on the top row and those on the bottom row.
  The first two on the top row have penalties $2(0.008)$ and $0.008$.
  Those on the bottom row have penalties
  $3\xiG$, $3\xiG$, $\xiG+2\xiV$, and $\xiG+2\xiV$.}
  \label{fig:hexrefine}
\end{figure}

The conventions for generating the possibilities are different for
the pentagons and hexagons than for the heptagons and octagons. We
describe the pentagons and hexagons first.  We erase $\Splus$. If
there is one loop we leave the loop in the figure. If there are
two loops (so that both  necessarily have context $(n,k)=(4,1)$),
we erase one and keep the other.  We draw all edges from an
enclosed upright diagonal to its anchors, and all edges of length
$[2t_0,2\sqrt2]$ that are not masked by upright quarters. Since
the only remaining upright quarters belong to loops, the four
simplices around a loop are anchored simplices and the edge
opposite the diagonal has length at most $3.2$.

\subsection{Branching on pentagonal faces}

The first three inequalities of Appendix \ref{app:pent} and the
first inequality of Appendix \ref{app:pentB}  have been designed
with pentagons in mind. The additional inequalities in Appendix
\ref{app:hexA}, \ref{app:hexB}, \ref{app:hexC}, and the last two
inequalities of \ref{app:pent} have been designed for subregions
in hexagonal regions.

\subsection{Hept and Oct Branching}

When the figure is a heptagon or octagon, we proceed differently.
We erase all $\Splus$ configurations and all loops (either context
$(n,k)=(4,1)$ or $(4,2)$) and draw only  the flat quarters.  An
undrawn diagonal of the polygon has length at least $2t_0$.  Thus,
in these cases much less internal structure is represented.

\begin{figure}[htb]
  \centering
  \includegraphics{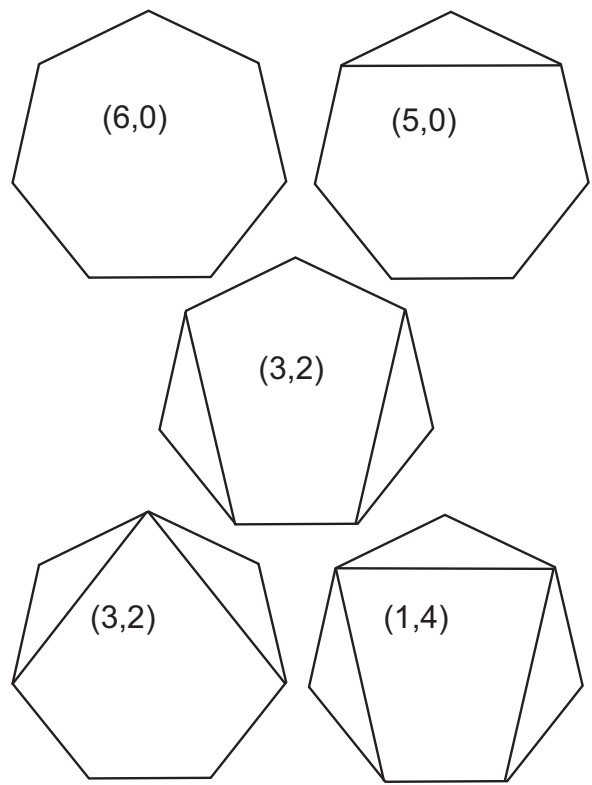}
  \caption{Hept Face Refinements. The pair $(a,b)$ represents the penalty $a\xiG+b\xiV$.}
  \label{fig:heprefine}
\end{figure}

\begin{figure}[htb]
  \centering
  \includegraphics{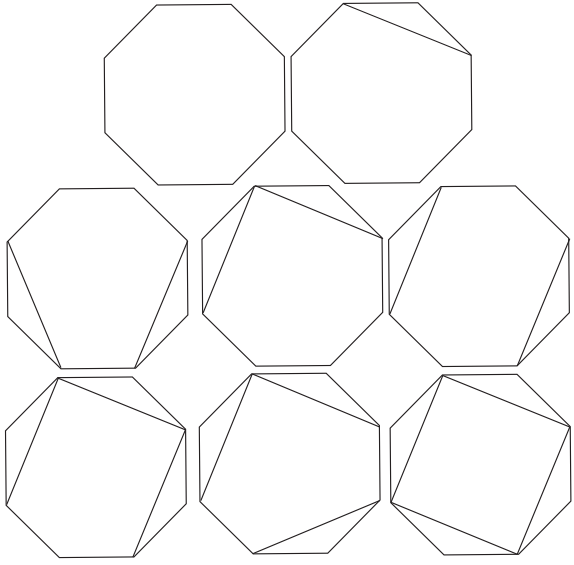}
  \caption{Oct Face Refinements.  The penalties reading from left to right
    are $6\xiG$ and $6\xiG$ on the top row, $4\xiG+2\xiV$, $4\xiG+2\xiV$ and
        $4\xiG+2\xiV$ on the middle row, and $2\xiG+4\xiV$, $2\xiG+4\xiV$, and $0$
        on the bottom row.}
  \label{fig:octrefine}
\end{figure}

In the cases where $\Splus$ or loops have been erased, a number
indicating a penalty accompanies the diagram.  These penalties are
derived in Section \ref{sec:4.2}.

 Here are some special arguments that are
used for heptagons and octagons.

Suppose that the face refinement gives two subregions: a flat
quarter and one with parameters $(n,k)=(n-1,1)$, with $n=6$ or
$7$. We have the inequality:
    $$\sigma_R(D) < \hat\sigma + Z(n-1,1) +\xiG+2\xiV.$$
The penalty term $\xiG+2\xiV$ comes from a possible anchored
simplex that masks the flat quarter.   Let $v$ be the central
vertex of the flat quarter.  Let $(v_1,v_2)$ be its diagonal.
Masked flat quarters satisfy restrictive edge constraints.  Thus,
we have one of the following three possibilities:
    \begin{enumerate}
    \item $y[v] \ge 2.2$,
    \item $y[v_1,v_2] \ge2.7$,
    \item $\sigma_R(D) < \hat\sigma + Z(n-1,1).$
    \end{enumerate}

We proceed similarly, if the face refinement has three subregions:
two flat quarters $R_1$, $R_2$ and one ($R_3$) with parameters
$(n,k)=(n-2,2)$. We have the cases:
    \begin{enumerate}
    \item The height of a central vertex is at least $2.2$.
    \item The diagonal of a flat quarter is at least $2.7$.
    \item
        $$
        \begin{array}{lll}
        \sigma_R(D) &< \hat\sigma_1 +\hat\sigma_2+Z(n-2,2)\\
        \tau_R(D) &>\hat\tau_1+\hat\tau_2+D(n-2,2).
        \end{array}$$
    \end{enumerate}

With heptagons, it is helpful on occasion to use an upper bound on
the penalty of $3\xiG = 0.04683$. This bound holds if neither flat
quarter is masked by a loop. For this, it suffices to show that
first two of the given three cases do not hold.

\subsubsection{Octagon Branching}.

If there is a loop of context $(n,k)=(4,2)$, we have, by Lemma
\ref{lemma:dloop}, the upper bound
    $$s_n+\zloop(4,2)$$
on the score of exceptional standard regions $R$, with $n = n(R)$.

If there is no loop of context $(n,k)=(4,2)$, and if there are two
flat quarters in the face refinement, we have the upper bound
    $$\hat\sigma_1 + \hat\sigma_2 + D(n-2,2) + 2(\xiG+2\xiV).$$

\subsection{ Penalties} 
\label{sec:4.2}

The reader should review Appendix \ref{app:part4} before reading
this section.

Erasing an upright quarter of compression type gives a penalty of
at most $\xiG$ and one of Voronoi type gives at most $\xiV$. We
take the worst possible penalty.  It is at most $n\xiG$ in an
$n$-gon. If there is a masked flat quarter, the penalty is at most
$2\xi_V$ from the two upright quarters along the flat quarter.  We
note in this connection that both edges of a polygon along a flat
quarter lie on upright quarters, or neither does.

If an upright diagonal appears enclosed over a flat quarter, the
flat quarter is part of a loop with context $(n,k)=(4,1)$, for a
penalty at most $2\xi'_\Gamma+\xi_V$.  This is smaller than the
bound on the penalty obtained from a loop with context
$(n,k)=(4,1)$, when the upright diagonal is not enclosed over the
flat quarter:
    $$\xi_\Gamma + 2\xi_V.$$
So we calculate the worst-case penalties under the assumption
that the upright diagonals are not enclosed over flat quarters.

A loop of context $(n,k)=(4,1)$ gives $\xi_\Gamma+2\xi_V$ or
$3\xi_\Gamma$.  A loop of context $(n,k)=(4,2)$ gives
$2\xi_\Gamma$ or $2\xi_V$.

  If we erase $\Splus$,
there is a penalty of $0.008$ (or 0 if it masks a flat quarter.)
This is dominated by the penalty $3\xi_\Gamma$ of context
$(n,k)=(4,1)$.

Suppose we have an octagonal standard region.  We claim that a
loop does not occur in context $(n,k)=(4,2)$. If there are at most
three vertices that are not corners of the octagon, then there are
at most 12 quasi-regular tetrahedra, and the score is at most
$$s_8 + 12\,\pt<8\,\pt.$$
Assume there are more than three vertices that are not corners
over the octagon. We squander
$$t_8+ \dloop(4,2)+4\tlp(5,0) > \squander.$$
As a consequence, context $(n,k)=(4,2)$ does not occur.

So there are at most $2$ upright diagonals and at most $6$ quarters,
and the penalty is at most
$6\xi_\Gamma$. Let $f$ be the number of flat quarters
This leads to
    $$
    \piF = \begin{cases} 6\xiG, & f=0,1,\\
                   4\xiG+2\xiV, & f=2,\\
                    2\xiG+4\xiV, & f=3,\\
                    0, & f=4.
            \end{cases}
    $$
The 0 is justified by a parity argument.  Each upright quarter occurs
in a pair at each masked flat quarter.  But there is an odd
number of quarters along the upright diagonal, so no penalty at
all can occur.

Suppose we have a heptagonal standard region.  Three loops are
a geometric impossibility. Assume there are at most two upright
diagonals.
 If there is no context $(n,k)=(4,2)$,
 then we have the following bounds on the penalty
    $$
    \piF = \begin{cases} 6\xiG, & f=0,\\
                 4\xiG+2\xiV, & f=1,\\
                3\xiG, & f=2,\\
                \xiG+2\xiV, & f=3.
            \end{cases}
    $$
If an upright diagonal has context $(n,k)=(4,2)$, then
    $$
    \piF = \begin{cases} 5\xiG, & f=0,1,\\
                3\xiG + 2\xiV, & f=2,\\
                \xiG + 4 \xiV, &f = 3.\\
            \end{cases}
    $$
This gives the bounds used in the diagrams of cases.

\subsection{Branching on Upright Diagonals} 
\label{sec:4.12}

We divide the upright simplices into two domains depending on the
height of the upright diagonal, using $y_1=2.696$ as the break
point. We break the upright quarters into cases:
    \begin{enumerate}
    \item The upright diagonal has height at most $2.696$.
    \item The upright diagonal $(0,v)$ has height at least $2.696$,
        and some anchor $w$ along the flat quarter satisfies
        $|w|\ge 2.45$ and $|v-w|\ge2.45$.  (There is a separate
        case here for each anchor $w$.)
    \item The upright diagonal $(0,v)$ has height at least
    $2.696$,
        and every anchor $w$ along the flat quarter satisfies
        $|w|\le 2.45$ and $|v-w|\le2.45$.
    \end{enumerate}
Many inequalities have been specially designed to hold on these
smaller domains.  They are included into the linear programming
problems as appropriate.

When all the upright quarters can be erased, then the case for
upright quarters follows from some other case without the upright
quarters.  An upright quarter can be erased in the following
situations.  If the upright quarter is compression type and
$\ge2.696$, then
    $$\nu_\Gamma<\vor_0$$
(\cite[$\A_{10}$]{part4}). (If there are masked flat quarters,
they become scored by $\hat\sigma$.) If upright quarters is of
Voronoi type and the anchors $w$ satisfy $|w|\le 2.45$ and
$|v-w|\le2.45$, then the quarter can be erased:
    $$\nu<\vor_0$$
In general, we only have the weaker inequality
(\cite[$\A_{11}$]{part4})
    $$\nu < \vor_0+\xiV, \quad \xiV = 0.003521.$$

\subsection{Branching on Upright Quarters} 
\label{sec:6}

In a pentagon or hexagon, consider an upright diagonal with three
upright quarters, that is, context $(n,k)=(4,1)$. If the upright
diagonal is at most $2.696$, and if an upright quarter shares both
faces along the upright diagonal with other upright quarters, then
we may assume that the upright quarter has compression type. For
otherwise, there is a face of circumradius at least $\sqrt2$, and
hence two upright quarters of Voronoi type.  The inequality
    \begin{equation} \octavor <
    \octavor_0 - 0.008,
    \label{eqn:6.1}
    \end{equation}
if $y_1\in[2t_0,2.696]$, and $\eta_{126}\ge\sqrt2$ shows that the
upright quarters can be erased without penalty because $\xiG <
2(0.008)$.  If erased, the case is treated as part of a different
case.

This allows the inequalities for $\nu_\Gamma$ in Appendix
\ref{app:nugamma1} and \ref{app:nugamma} to be used. Furthermore,
it can often be concluded that all three upright quarters have
compression type. For this, we use the Inequalities \ref{app:471}
and \ref{app:455}, which can often be used to show that if the
anchored simplex has a face of circumradius at least $\sqrt2$,
then the linear programming bound on the score is less than
$8\,\pt$.

\subsection{Branching on Flat Quarters}

There are a few  other interval-based inequalities that are used
in particular cases. Before an inequality that uses compression
scoring is used on a flat quarter, it is necessary to verify that
the quarter satisfies
    $y_1\le 2.2, y_4\le2.7, \eta_{234},\eta_{456}\le\sqrt2$,
to insure that the flat quarter has compression type. The
circumradius is not a linear-programming variable, so its upper
bound must be deduced from edge-length information.

Information about the internal structure of an exceptional face
gives improvements to the constants $1.4\,\pt$ and $1.5\,\pt$ of
Property \ref{definition:admissible:excess} in the definition of
admissible weight assignments. (The bounds remain fixed at
$1.4\,\pt$ and $1.5\,\pt$, but these arguments allow us to specify
more precisely what simplices contribute to these bounds.) These
constants contribute to the bound on the $\tau(D)$ through the
admissible weight assignment. Assume that at the vertex $v$ there
are four tetrahedra and an exceptional face, and that the
exceptional face has a flat quarter with central vertex $v$. The
calculations of Section \ref{sec:tri4} show that the four quarters
and exceptional squander at least $1.5\,\pt$. If there is no flat
quarter (masked or unmasked) whose first edge lies along $(0,v)$,
then the four quasi-regular tetrahedra at $v$ squander at least
$1.5\,\pt$.  We can make similar improvements when $tri(v)=3$.

If we have LP-derived inequalities on a flat quarter whose corners
$v_i$ have height at most 2.14, and if the diagonal has length
less than $2.77$, then the circumradius of the face containing the
origin and diagonal is at most $\eta(2.14,2.14,2.77)<\sqrt2$. This
allows us to combine the cases defining $\hat\sigma$ into three
cases.
    \begin{enumerate}
    \item The simplex has compression type.
    \item The diagonal has length $\le2.7$ and $\eta_{456}\ge\sqrt2$.
    \item The diagonal has length $\ge2.7$.
    \end{enumerate}
In the last two cases, the scoring is by $\vor_0$.  We use the
specially tailored Inequalities \ref {app:A.7.1} and \ref
{app:A.7.2}.

\subsection{Branching on Type A simplices} 

Suppose that the vertices of a triangular subregion come by
projection from a simplex with vertices at the origin and  in
$U(D)$. Suppose that the simplex has two edges of length in
$[2t_0,2\sqrt2]$, say $y_5,y_6\in[2t_0,2\sqrt2]$,
$y_4\in[2,2t_0]$. These may be simplices of type $A$. Suppose that
one edge is at most $2.77$. Such simplices fall into the following
domains.
    \begin{enumerate}
    \item
    The simplex has type $A$, so that $y_5,y_6\in[2t_0,2.77]$,
    and the scoring is $\vor$, the analytic Voronoi function.
    \item An edge (say $y_6$) of the face shared with the flat quarter
    has length $y_6\ge2.77$.
    \item  The edges have lengths $y_5,y_6\in[2t_0,2.77]$, and $\eta_{456}\ge\sqrt2$.
    \end{enumerate}
In the last two cases, the scoring is by the truncated Voronoi
function $\vor_0$. We use the specially tailored Inequalities
\ref{app:17}.

If there is an enclosed vertex of height at most $\sqrt2$, then we
can use Inequalities \ref{app:A.7.4}.

\subsection{Branching on Quadrilateral subregions}

The inequalities of  \ref{app:471} hold for a quadrilateral
subregion, if certain conditions are satisfied.  One of the
conditions is $y_4\in[2\sqrt2,3.0]$, where $y_4$ is a diagonal of
the subregion. Since this diagonal is not one of the linear
programming variables, these bounds cannot be verified directly
from the linear program.  Instead we use the interval calculation
\ref{ineq:1.678} which relates the desired bound $y_4\le3$ to the
linear programming variables $\alpha[v,F]$, $y_2$, $y_3$, $y_5$,
and $y_6$.

\subsection{Conclusion}

By combinations of branching along the lines set forth in the
preceding sections, a sequence of linear programs is obtained that
establishes that  $\sigma(D)$ is less than $8\,\pt$.  For details
of particular cases, the interested reader can consult the log
files in \cite{web}, which record which branches are followed for
any given tame graph.  (For most tame graphs, a single linear
program suffices.)

This completes the proof of the Kepler conjecture.


\newpage

\section{Appendix. Colored Spaces}
\label{app:coloredspace}

It is claimed at various places that the optimization problem is
that of a continuous function on a compact topological space. That
is, the function $\sigma$ on the space $X$ of decomposition stars
is continuous. The purpose of this appendix is to describe the
sense in which this is true.

\begin{figure}[htb]
  \centering
  \includegraphics{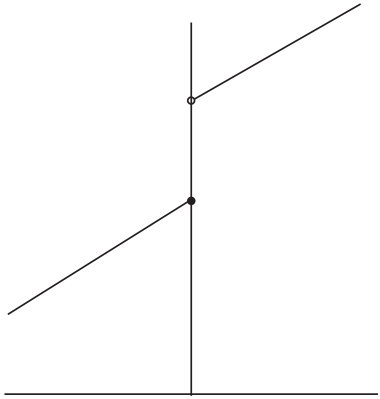}
  \caption{A piecewise linear function}
  \label{fig:discontinuous}
\end{figure}

We begin with an example that should make things clear.  Suppose
that we have a discontinuous piecewise linear function on the unit
interval $[-1,1]$, as in Figure \ref{fig:discontinuous}.  It is
continuous, except at $x=0$.

We break the interval in two at $x=0$, forming two compact
intervals $[-1,0]$ and $[0,1]$.  We have a continuous functions
$f_-:[-1,0]\to \R$ and $f_+:[0,1]$, such that
    $$
    f(x)=
    \begin{cases}
        f_-(x) & x\in[-1,0],\\
        f_+(x) & \text{otherwise}.
    \end{cases}
    $$
We have replaced the discontinuous function by a pair of
continuous functions on smaller intervals, at the expense of
duplicating the point of discontinuity $x=0$.  We view this pair
of functions as a single function $F$ on the compact topological
space with two components
    $$[-1,0]\times\{-\} \text{\ and\ } [0,1]\times\{+\}.$$
where $F(x,a) = f_a(x)$, and $a\in\{-,+\}$.

This is the approach that we follow in general with the Kepler
conjecture.  The function $\sigma$ is defined by a series of case
statements, and the function does not extend continuously across
the boundary of the cases.  However, in the degenerate cases that
land precisely between two or more cases, we form multiple copies
of decomposition star for each case, and place each case into a
separate compact domain on which the function $\sigma$ is
continuous.

This can be formalized as a {\it colored space}.  A colored space
is a topological space $X$ together with an equivalence relation
on $X$ with the property that no point $x$ is equivalent to any
other point in the same connected component as $x$.   We refer to
the connected components as colors, and call the points of $X$
{\it colored points}.  We call the set of equivalence classes of
$X$ the underlying uncolored space of $X$.   Two colored points
are equal as uncolored points if they are equivalent under the
equivalence relation.

In our example, there are two colors ``$-$'' and ``$+$.''  The
equivalence class of $(x,a)$ is the set of pairs $(x,b)$ with the
same first coordinate.  Thus, if $x\ne0$, the equivalence class
contains one element $(x,\sign(x))$, and in the boundary case
$x=0$ there are two equivalent elements $(0,-)$ and $(0,+)$.

In our treatment of the Kepler Conjecture, there are various
cases: whether an edge has length less than or greater than
$2t_0$, less than or greater than $\sqrt8$, whether a face has
circumradius less than or greater than $\sqrt2$, and so forth. By
duplicating the degenerate cases (say an edge of exact length
$2t_0$), creating a separate connected component for each case,
and expressing the optimization problem on a colored space, we
obtain a continuous function $\sigma$ on a compact domain $X$.

The colorings have in general been suppressed from the notation of
the paper.  To obtain consistent results, a statement about
$x\in[2,2t_0]$ should be interpreted as having an implicit
condition saying that $x$ has the coloring induced from the
coloring on the component containing $[2,2t_0]$.  A later
statement about $y\in[2t_0,\sqrt8]$ deals with $y$ of a different
color, and no relation between $x$ and $y$ of different colors is
assumed at the endpoint $2t_0$.

In general, to the greatest extent possible, we express each
inequality as a strict inequality on a compact domain.  There are,
however, a few inequalities that are not strict, such as the bound
of $1\,\pt$ on the score of a quasi-regular tetrahedron or the
bound of $0$ on the score of a quad cluster.  (These particular
sharp bounds appear in the proof of the local optimality of the
face-centered cubic and hexagonal close packings.)

\section{Appendix.  Constants and Functions}
\label{app:constants}

This appendix gives a summary of various constants and functions
that appear in the proof.

The {\it density} of a regular tetrahedron and Rogers's bound on
the density of a packing is
    $$\dtet = \sqrt8 \arctan(\sqrt2/5).$$
The density of a regular octahedron is
    $$\doct = \frac{\pi}{\sqrt8} - \frac{\dtet}{2}\approx 0.72.$$
The unit of measure in calculations of $\sigma(D)$ is a point:
    $$\pt = -\pi/3 + \sqrt2\dtet \approx 0.05537.$$
A contravening decomposition star $D$ satisfies
    $$\sigma(D)\ge 8\,\pt$$
or equivalently,
    $$\tau(D) \le \squander,$$
where
    $$\zeta = 1/(2\arctan(\sqrt2/5)).$$

Objects are truncated using the value $t_0=1.255$.  Centers of the
packing at distance at most $2t_0$ from the origin form a set
$U(D)$.

Internal angles of triangular standard regions are at least
$0.8638$, and those of other standard regions are at least
$1.153$.  The internal angles of standard regions are generally
denoted with variables $\alpha$.  The dihedral angles of simplices
$S$ are denoted $\dih(S)$.  If the edge of a simplex has index
$i$, then we write $\dih_i(S)$ for the dihedral angle along the
edge with index $i$.

The circumradius of a triangle is denoted $\eta$.  If there is a
simplex with indexed edges, then $\eta_{ijk}$ denotes the
circumradius of the face of the simplex whose edges have indices
$i$, $j$, and $k$.

In general, the edges of simplices are indexed $1,\ldots,6$ so
that the first three edges meet at a vertex, and so that opposite
edges have indices that are congruent mod $3$. We write
$y_1,\ldots,y_6$ for the edge lengths and $x_i=y_i^2$ as their
squares.  In general, the vertex where the first three edges meet
is fixed as a distinguished vertex of the simplex.  We write
$\sol(S)$ for the solid angle of the simplex at that vertex.

The function $\sigma$ on the topological space of decomposition
stars is the function optimized in this paper.  It is closely
related to
    $$\tau(D) = \sol(D)\zeta\pt - \sigma(D).$$
The function $\sigma$ is expressed as a sum of functions
$\sigma_R$, indexed by standard regions.  These sums are in turn
expressed as finer sums.

The formulas for $\sigma_R$ and their upper bounds contain
functions $\Gamma$, $\mu$, $\nu$, $\quo$, $\Adih$, $\vor$,
$\vor_0$.  We do not repeat these definitions (they are found in
\cite{formulation}), because they appear in this paper in a purely
formal way.  (We make use of interval arithmetic inequalities that
these functions satisfy, but we do not make use of anything deeper
about them.)  The function $\hat\sigma$  is defined in
\cite[Section 3.11]{part4}. The function  $\vor_x$ is a function
on anchored simplices that are not upright quarters.  It is $\vor$
or $\vor_0$ depending on whether the upright simplex has type $C$.
Also, we set
    $$\vor_x=
        \begin{cases}
        \vor& \text{type $A$, $C$, $C'$}\\
        \vor_0& \text{otherwise}
        \end{cases}
    $$

Corresponding to these functions related to $\sigma$ are functions
related to $\tau$:
    $$
    \begin{array}{lll}
        \tau_x &= \sol\zeta\pt - \vor_x\\
        \tau_R &= \sol\zeta\pt - \sigma_R\\
        \hat\tau &= \sol\zeta\pt - \hat\sigma\\
        \tau &= \sol\zeta\pt - \sigma\\
        \tau_0 &= \sol\zeta\pt - \vor_0\\
        \tau_\Gamma &=\sol\zeta\pt - \Gamma\\
        \tau_\nu &=\sol\zeta\pt-\nu\\
        \tau_\mu &=\sol\zeta\pt -\mu\\
        \tau_V &=\sol\zeta\pt -\vor\\
        \tau(\ ,t) &= \sol\zeta\pt - \vor(\ ,t).
    \end{array}
    $$
In the last case, the truncation is at $t=\sqrt2$, $1.385$, $t_0$,
and so forth.  In each case, $\sol$ is the function that measures
the solid angle of the object chosen from the domain of the
function. Each function on the left-hand side is referred to as a
{\it squander}, and each function that are subtracted on the
right-hand side is a {\it score}.  When the function $\tau_*$ is
related to the scoring function in this way, we say that it is
{\it adapted}.

The function $\vor_0$ is used in many contexts, and we sometimes
take liberties with the notation for this function.  For example,
we write $\vor_{0,R}(U(D))$ to indicate that the function depends
only on the set $U(D)$ and the combinatorial information encoded
by the region (or subregion) $R$.

The polynomial $\Delta$ of six variables has the property that
$\Delta(x_1,\ldots,x_6)<0$, then the simplex with squared edge
lengths $x_i$ does not exist.  In the case of equality $\Delta=0$,
the simplex degenerates to a planar quadrilateral.


\section{Appendix.  A Review of \cite{part4}}
\label{app:part4}

\subsection{Standard regions}

This paper relies heavily on results about the internal structure
of standard regions that is developed in \cite{part4}.  We make a
brief review of these structures.  We restrict our attention to
exceptional standard regions.  A {\it quarter} is a simplex, whose
vertices are vertices $v\in\Lambda$ in the packing, such that five
edges lengths are in $[2,2t_0]$ and the sixth edge length lies in
the interval $[2t_0,\sqrt8]$.  The sixth edge of the quarter is
called its {\it diagonal}.

We fix our attention on the decomposition star $D$ centered at a
given vertex in the packing, which we take to be located at the
origin.  The decomposition star determines the set $U(D)$ of all
vertices of $\Lambda$ at distance at most $2t_0$ from the origin.

From the decomposition star the set of quarters with a vertex at
the origin can be also determined.  There are two types. If the
diagonal of the quarter has an endpoint at the origin, it is said
to be {\it upright}.  Otherwise, it is said to be {\it flat}.  A
flat quarter has its vertices at the origin and three points of
the set $U(D)$.  An upright quarter has one of its vertices $v$
outside $U(D)$.  The edge $(0,v)$ is an {\it upright\/} diagonal.

The decomposition star determines a distinguished set of quarters
with vertex at the origin.  These distinguished quarters are said
to belong to the $Q$-{\it system}.  If one quarter along a
diagonal lies in the $Q$-system, then every quarter along that
diagonal lies in the $Q$-system.  Quarters in the $Q$-system do
not overlap (in the sense that they have disjoint interiors). The
cones over the quarters at the origin in the $Q$-system do not
overlap.

Each quarter lies entirely in the cone over a single standard
region.  We fix our attention on a single standard region and the
quarters that are associated with it in this way.  We write
$U_R(D)$ for the subset of $U(D)$ consisting of points whose
projection to the unit sphere lands in the (closed) region $R$.

There is a continuous function $\sigma$ defined on the topological
space of all decomposition stars.  It can be expressed as a sum of
terms
    $$\sigma = \sum_R \sigma_R,$$
indexed by the set of standard regions.  We will describe the
structural properties of $\sigma_R$ for a given standard region
$R$.  There is a closely related function
    $$\tau_R = \sol(R)\zeta\pt - \sigma_R.$$
The theory develops in parallel for these two functions.


\subsection{Flat quarters}
\label{sec:flat}

Flat quarters appear in several ways.  We color (in the sense of
Appendix \ref{app:coloredspace}) each quarter with identifying
information that marks how that flat quarter arises.  We will
define functions
    $\hat\sigma$ and $\hat\sigma_{adj}$
on the disjoint union of these different types of flat quarters.

If a flat quarter is in the $Q$-system, then we set
    $$\hat\sigma(S) = \hat\sigma_{adj}(S) = \mu(S).$$
Where $\mu(S)$ is the function defined in
\cite[Formula~3.8]{formulation}.

Another type of flat quarter are those appearing in {\it isolated
pairs} (see \cite[Figure 1.10]{formulation}).   A special subset
of such flat quarters are called simplices of type $B$. They have
the following properties.
\begin{enumerate}
    \item The simplex is a flat quarter.
    \item The diagonal has length at least $2.77$.
    \item The edges running between the endpoints of the diagonal
    and the origin have length at most $2.23$.
\end{enumerate}
For simplices of type $B$, set
    $$\hat\sigma(S)=\hat\sigma_{adj}(S) = \vor(S,1.385).$$
For simplices that are part of an isolated pair, but that are not
of type $B$, we set
    $$\hat\sigma(S)=\hat\sigma_{adj}(S)= \vor_0(S).$$

There are flat quarters that enclose an upright diagonal along an
upright diagonal in the $Q$-system.  We call these {\it masked
flat quarters of the first type}.  In this case,  there are two
upright quarters $Q_1$ and $Q_2$ over the flat quarter, as well as
a third simplex $S$. We set
    $$\hat\sigma_{adj} = \nu(Q_1)+\nu(Q_2)+\vor_x(S),$$
where $\vor_x=\vor$ if $S$ has type $C$, and $\vor_x=\vor_0$
otherwise. We note that $\eta_{456}(S)\ge\sqrt2$.  Also, set
    $$\hat\sigma(S) = \min(\vor_0(S),0).$$
We have
    $$\hat\sigma_{adj}\le \hat\sigma(S).$$
    \cite[Lemma 3.11.1]{part4}

Let $(0,v,v_1,v_2)$ be a flat quarter that is not any of those
discussed above.  Let $v$ be the central vertex, and let
$(v_1,v_2)$ be its diagonal.  (Recall that the central vertex is
defined as the vertex on a flat quarter that is not the origin and
is not an endpoint of the diagonal.)  These remaining flat
quarters have the properties:
    \begin{enumerate}
    \item $|v_1-v_2|\ge 2.6$.
    \item $|v_1-v_2|\ge 2.7$ or $|v|\ge2.2$.
    \end{enumerate}
We call these quarters {\it masked flat quarters} of the second
type. Set
    $$
    \epsilon=
    \begin{cases}
        0.0063 & |v|\ge2.2,\\
        0.0114 & |v|\le2.2.
    \end{cases}
    $$
    $$
    \begin{array}{lll}
    \hat\sigma_{adj}(S) &= \vor_0(S) - \epsilon.\\
    \hat\sigma(S) &=\vor_0(S).\\
    \end{array}
    $$

For all types of flat quarters $S$, set
    $$\hat\tau(S) = \sol(S)\zeta\pt -\hat\sigma(S).$$
We have
    $$
    \begin{array}{lll}
    \hat\sigma_{adj}&\le\hat\sigma,\\
    \hat\sigma_{adj}&\le Z(3,1),\\
    \hat\tau_{adj}&\ge D(3,1),
    \end{array}
    $$
where $Z(3,1) = 0.00005$ and $D(3,1)=0.06585$. (See
\cite[5.5]{part4}.)

\subsection{Type A}
\label{sec:typeAB}

 A simplex $S$ of type $A$ is a simplex with a
vertex at the origin and at three other points $v_1,v_2,v_3$ in
$U_R(D)$ such that
\begin{enumerate}
    \item The circumradius of the triangle with vertices
    $v_1,v_2,v_3$ is less than $\sqrt2$.
    \item Two of the edges $(v_i,v_j)$ have edge length in the
    interval $[2t_0,2.77]$.
    \item The third edge has length in the interval $[2,2t_0]$.
\end{enumerate}

\subsection{Loops}

Let $(0,v)$ be an upright diagonal in the cone over a standard
region $R$. We assume throughout this discussion that there is an
upright quarter in the $Q$-system with diagonal $(0,v)$.  A vertex
$w$ (distinct from $0$ and $v$) of distance at most $2t_0$ from
both $0$ and $v$ is called an {\it anchor\/} (of $v$). By
definition, we have $w\in U(D)$.

\begin{lemma}
    The projection of $v$ to the unit sphere does not lie in the
    subregions associated with any flat quarter in the $Q$-system,
    any simplex of type $A$ or any simplex of type $B$.  The projection $e$ of
    an edge $(v,w)$ does not meet the interior of a subregion associated
    with any flat quarter in the $Q$-system, a simplex of type
    $A$, or a simplex of type $B$.
\end{lemma}

The edges $(v,w)$ from $v$ to its anchors $w$ emanate like spokes
from $v$.  We order the anchors of $v$ according to this cyclic
order into a cycle and write $s_v(w)$ for the successor of the
anchor $w$ in this ordering.  If the angle swept out along $(0,v)$
in passing from one spoke to the next is greater than $\pi$, we
say that the angle around $(0,v)$ between $w$ and $s_v(w)$ is {\it
concave}.  Otherwise, we say it the angle is convex.

We say that a simplex $(0,v,w,s_v(w))$ along $(0,v)$ is an {\it
anchored simplex}, if the angle is convex and if
$|w-s_v(w)|\le3.2$.  We say that $(0,v,w,s_v(w))$ is a {\it large
gap}, if the angle is convex and if $|w-s_v(w)|\ge 3.2$.

Two other anchored simplices $S(y_1,\ldots,y_6)$ that arise are
those of type $C$ or $C'$. Write the anchored simplex $S$ with
edge $y_1$ representing the length of the upright diagonal. We say
the anchored simplex has type $C$ if $y_4\in[2t_0,\sqrt8]$ and
    \begin{enumerate}
    \item $y_4\le 2.77$, or
    \item Both face circumradii
    of $S$ along the fourth edge are
    at most $\sqrt{2}$.
    \end{enumerate}
We say that the anchored simplex has type $C'$ if
    \begin{enumerate}
    \item $y_4\in[2t_0,\sqrt8]$,
    \item $y_2,y_6\in[2.45,2t_0]$,
    \item The fourth edge is a diagonal of a flat quarter that
        shares a face with $S$.
    \end{enumerate}

If the anchored simplices around $(0,v)$ are not contained in a
half-space bounded by a plane through $(0,v)$ and if there are no
large gaps, then the anchored simplices around $v$ are said to
form a {\it loop} around $v$.  By construction, the anchored
simplices around $v$ do not overlap. Anchored simplices (even
those around different vertices) in a contravening decomposition
star do not overlap. Anchored simplices do not overlap flat
quarters in the $Q$-system, nor do they overlap simplices of types
$A$ and $B$.

A {\it special simplex} is a simplex $(0,v_1,v_2,v_3)$ with
vertices in $U_R(D)$ such that for some ordering of the indices,
$|v_1-v_2|\le 2t_0$, $|v_2-v_3|\le 2t_0$, and
$|v_1-v_3|\in[\sqrt8,3.2]$.  A special simplex $(0,v_1,v_2,v_3)$
is said to {\it lie along\/} an anchored simplex $(0,v,w,s_v(w))$
if the simplices can be indexed so that
$\{v_1,v_3\}=\{w,s_v(w)\}$, with $|v_1-v_3|\ge\sqrt8$.  Since
anchored simplices do not overlap, each special simplex is related
to at most one upright diagonal $(0,v)$ in this way.  A special
simplex attached to an upright diagonal in this way does not
overlap any of the preceding objects (other special simplices so
attached, other anchored simplices, simplices of types $A$ and
$B$, and flat quarters in the $Q$-system).

\subsection{Bounds on $\sigma_R$ and $\tau_R$.}
\label{sec:tnsn}

Let $R$ be a standard region. More generally, we can partition $R$
into subregions by considering three types of bounding edges:
\begin{enumerate}
    \item a bounding edge of a standard region,
    \item the projection (to the unit sphere)
        of a face (along the diagonal) of an upright
        quarter in the $Q$-system,
    \item the projection of a face formed by $0$ and $v_1,v_2\in U(D)$,
        where $|v_1-v_2|\le\sqrt8$, and additionally if it bounds a
            flat quarter and that flat quarter is not in the $Q$-system then
            the face also bounds a simplex of type $A$.
\end{enumerate}
Let $R$ be one of these subregions. We have
    $$\sigma(D) = \sum_R\sigma_R(D).$$
If $R$ is the projection of a simplex $Q$ in the $Q$-system, then
$\sigma_R(D) = \sigma(Q)$, for some function $\sigma$ depending
only on the simplex. Otherwise $\sigma_R(D)$ has the form
$$\sigma_R(D) = \sigma_R(V_D),$$
where $V_D$ is the $V$-cell  defined in
\cite[Section~2]{formulation}. The $V$-cell can be truncated at
$\sqrt2$, which leads to the upper bound
    $$\sigma_R(V_D)\le \vor_R(V_D,\sqrt2)\le0.$$
Each term indexed by quarters satisfies
    $$\sigma(Q)\le0.$$
In particular, if we have any subset ${\cal Q}$ of the set of
quarters over a standard region $R$, then
    $$\sigma_R(D) \le \sum_{\cal Q} \sigma(Q).$$

There is a second partition of subregions that is used more
frequently.  Let $R$ be a subregion of the preceding paragraph. We
break any anchored simplices that might occur from it (provided
that the subregion does not come from a special simplex that lies
along the anchored simplex).   If this involves breaking the
subregion into smaller pieces, we switch from the function
$\sigma_R$ to the upper bound $\vor_0$ (together with the volumes
$\delta_P(v)$ of \cite[2.11]{part4} . We also break off any flat
quarters that are not masked (we switch to the upper bound
$\hat\sigma_{adj}$ on the flat quarter and to $\vor_0$ on the rest
of the subregion).

According to a construction that extends over several sections of
\cite{part4} with a summary in \cite[Section 3.10]{part4}, it is
often possible to eliminate the upright quarters around an upright
diagonal by replacing the function $\sigma_R$ with an upper bound
$\vor_0$ over those pieces.  This process is called {\it erasing}
the upright quarters.  At other times, we get a weaker statement,
by replacing the functions $\sigma_R$ for regions $R$ at an
upright diagonal $(0,v)$ with
    $$\sum_R \sigma_R(D) \le \sum_R \vor_{0,R}(U(D))+\pi(v,D).$$
The constant $\pi(v,D)$ is called a penalty.

 There are other upper bounds for $\sigma_R$ of
interest. If we do not make approximations anywhere that produces
penalties, we have a rather lengthy expression for $\sigma_R$:
        \begin{eqnarray}
        \sigma_R(D) =
                &\ \sum_Q \hat\sigma(Q)  \\
                &+\sum_A \vor(S)  \\
                &+\sum_B\vor(S) \\
                &+\sum_{\text{isolated, not }B}\vor_0(Q)\\
                &+\sum_{\text{loop}:v} \sum_{S\mapsto v} \vor_x(S)\\
                &+\sum_{\text{non-loop:v}}\left(\delta_v+ \sum_{S\mapsto v}
                    \vor_x(S)\right)\\
                &+\sum_{\text{type} 1}\hat\sigma_{adj}\\
                &+\sum \vor_{0,R'}(U_R(D)).
                \label{eqn:partition}
        \end{eqnarray}
The sum indexed by $Q$ runs over flat quarters in the $Q$-system.
The sum indexed by $A$ runs over simplices of type $A$ over $R$.
The sum indexed by $B$ runs over simplices of type $B$ over $R$.
The next sum runs over isolated flat quarters that do not have
type $B$.  The next sum runs over the upright diagonals $(0,v)$ of
loops.  (We restrict the sum to those loops that have four or five
anchored simplices around it.)  Its inner sum runs over the
anchored simplices around $(0,v)$, as well as the special
simplices that lie along these anchored simplices.

The next sum runs over upright diagonals $(0,v)$ of upright
quarters in the $Q$-system that are not part of a loop.  The next
sum runs over masked flat quarters of type 1 for which there
enclosed upright diagonal $(0,v)$ bounds exactly two quarters. The
final sum is parameterized by subregions $R'$ that have not
otherwise been included, so that the subregions in the sum give a
partition of $R$.

Each summand in this formula for $\sigma_R$ defines a geometric
object whose cone cuts out a subregion $R'$ on the unit sphere. We
write $\sigma^+_{R'}(D)$ for the summand.  The summand is a single
term except for those indexed by upright diagonals.  The summands
indexed by upright diagonals include the full inner sum over
anchored simplices and so forth. We attach constants $(n,k)$ to
each such subregion.

\begin{figure}[htb]
  \centering
  \includegraphics{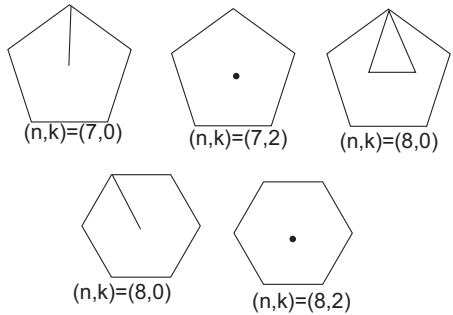}
  \caption{Parameters $(n,k)$ for non-polygonal regions}
  \label{fig:nonpolynk}
\end{figure}

To each of these subregions, we associate a natural number $n(R)$.
For polygons, it is the number of edges.  The edge is to be
counted with multiplicity two if it does not bound the standard
region.  We connect an internal dot that is not connected to any
other vertices to another vertex by an imaginary edge of length
greater than $2t_0$. Such an imaginary edge is counted with
multiplicity two.  The numbers $n(R)$ are indicated in Figure
\ref{fig:nonpolynk}. This constant is defined for the subregions
indexing the sums of Identity \ref{eqn:partition}.

A second parameter is $k(R)$ the number of sides of length at
least $2t_0$, counted with multiplicity $2$ if it is internal.
(Count the imaginary edge if it exists, with multiplicity two.)

Define constants
    $$
    \begin{array}{lll}
    t_4 = 0.1317 & t_5 = 0.27113 & t_6 = 0.41056,\\
    t_7 = 0.54999 & t_8 = 0.6045,\\
    s_5 = -0.05704 & s_6 = -0.11408 & s_7 = -0.17112,\\
    s_8 =-0.22816.
    \end{array}
    $$
Define constants
    $$
    D(n,k) = t_{n+k} - 0.06585 k,
    $$
for $0\le k\le n$ and $t_8\ge n+k\ge 4$. It is readily verified
that $$D(n_1,k_1)+D(n_2,k_2) \ge D(n_1+n_2-2,k_1+k_2-2).$$ Define
constants $\epsilon = 0.00005$ and
    $$
    Z(n,k) = s_{n+k} - k \epsilon,
    $$
for $(n,k)\ne(3,1)$ and $Z(3,1)=\epsilon$.  It is readily verified
that $$Z(n_1,k_1)+Z(n_2,k_2) \le Z(n_1+n_2-2,k_1+k_2-2).$$

\begin{theorem} For subregions $R'$ associated with summands as above, we have
    $$
    \begin{array}{lll}
    \sigma^+_{R'}&\le Z(n(R'),k(R')),\\
    \sol(R')\zeta\pt-\sigma^+_{R'}&\ge D(n(R'),k(R')).
    \end{array}
    $$
\end{theorem}

\subsection{Bounds for Standard Regions}

 Every standard region of a contravening decomposition star is a polygon with
between $3$ and $8$ sides or one of the  cases of Figure
\ref{fig:nonpolynk}. In these figures, we draw a dot to indicate
the projection of a vertex in $U(D)$ that does not lie on the
boundary.  These enclosed dots do not lie in the cone over a
simplex of type $A$, nor in the cone over a quarter.   We draw an
edge between vertices if the corresponding vertices in $U(D)$ have
distance at most $2t_0$ from each other. (This figure is a copy of
Figure \ref{fig:aggregates}.)

\begin{corollary} Define constants as in Section \ref{sec:tnsn}.
For every standard region $R$ such that $n=n(R)\ge5$:
    $$
    \begin{array}{lll}
    \sigma_{R}&\le s_n,\\
    \tau_{R}&\ge t_n.
    \end{array}
    $$
\end{corollary}

The corollary is an immediate consequence of the additivity
properties of $\sigma^+_{R'}$, the sub and superadditivity of $Z$
and $D$, and finally the geometric interpretation of sub and
superadditivity as gluing subregions together to form a larger
subregion.  See \cite[Section 4.5]{part4}.

\begin{lemma}
\label{lemma:dloop}
Let $R$ be a subregion cut out by a loop group.
Assume that the subregion is not that defined by a flat quarter.
Let $n$ be the number of anchors in the loop. Let $k$ be the
number of anchored simplices, with an edge (opposite the upright
diagonal) of length at least $2t_0$. Let $s$ be the number of
special simplices. In a contravening plane graph $(n,k)$ is one of
$(4,0)$, $(4,1)$, $(4,2)$, $(4,3)$, $(5,0)$, $(5,1)$. Furthermore,
    $$
    \sol(R)\zeta\pt-\sigma^+_R(D) \ge D(n+s,k-s) +\dloop(n,k),
    $$
where $\dloop$ is defined by the following constants. Similarly,
    $$
    \sigma^+_R(D)\le Z(n+s,k-s) + \zloop(n,k).
    $$
\end{lemma}

$$
\begin{matrix}
(n,k)   &   \dloop\\
(4,0)   &   0.0045\\
(4,1)   &   0.00272   \\
(4,2)   &   0.12034  \\
(4,3)   &  0.29426  \\
(5,0)   &   0.09537\\
(5,1)   &  0.24939\\
\end{matrix}
$$

$$
\begin{matrix}
(n,k) & \zloop\\
(4,2) & -0.08582\\
\end{matrix}
$$
(see \cite[Section 5.11]{part4}.)

\subsection{Remaining subregions}

A {\it masked\/} flat quarter is a flat quarter that is not in the
$Q$-system and that overlaps some anchored simplex around an
upright diagonal $(0,v)$.  The vertices of a masked flat quarter
consist of three consecutive anchors of $v$.

There are two cases to consider.  Either the upright diagonal
$(0,v)$ lies in the cone over the masked flat quarter or it does
not.    If it lies in the cone over the masked flat diagonal, then
the following facts are known.   (Recall that we are assuming
throughout this appendix that the standard region is not a
quadrilateral.)  As above, we call these masked flat quarters of
the {\it first type}.
    \begin{enumerate}
    \item The circumradius of the top face of the flat quarter is at least $\sqrt2$.
    \item The anchored simplices form a loop.
    \item The vertex $v$ has three or four anchors.
    \item If the vertex has three anchors, the loop consists of
    three anchored simplices that cut out the same region on the unit sphere as the
    masked flat quarter.
    \item If the vertex has four anchors, the loop consists of
    three upright quarters and one other that is not an upright
    quarter.
    \end{enumerate}

If the upright diagonal does not lie in the cone over the masked
flat quarter, then the following facts are known.  Call these
masked flat quarters of the {\it second type}.
    \begin{enumerate}
    \item The flat quarter has a diagonal of length at least
    $2.6$.
    \item If the flat quarter has a diagonal of length less than
    $2.7$, then the edge opposite the diagonal has length at least
    $2.2$.
    \end{enumerate}

\subsection{Penalties}

Set
    $$\piG(y) =
        \begin{cases}
        \xiG=0.01561 & 2t_0 \le y\le 2.57\\
        \xiG'=0.00935 &  2.57 \le y \le \sqrt8.\\
        \end{cases}
    $$
and
    $$\piV(y) =
        \begin{cases}
        0 & 2t_0\le y\le 2.57,\\
        \xiV=0.003521 &  2.57 \le y \le \sqrt8.\\
        \end{cases}
    $$

If $R$ is one of the subregions of Identity \ref{eqn:partition},
then we have a two functions
    $\vor_{0,R}(U(D))$ and $\sigma^+_R(D)$.
For example, $R$ may be the projection of the set of anchored
simplices around an upright diagonal $(0,v)$, together with the
set of special simplices that lie along the anchored simplices. In
this situation, we write $\pi(v,D)$ for the difference:
    $$\sigma^+_R(D) = \vor_{0,R}(U(D)) + \pi(v,D).$$
If $(0,v)$ is associated with a loop with $n$ anchored simplices,
then we have
    $$\pi(v,D)\le n\piG(|v|).$$

If there is a loop around $v$ with $n$ anchored simplices and with
$k$ masked flat quarters of the second type, then
    $$\pi(v,D) \le (n-2k)\piG(|v|) + 2k \piV(|v|).$$
We recall that for masked flat quarters $Q$ of the second type,
with projection to the subregion $R$, we have
    $$\hat\sigma(Q)=\vor_0(Q).$$

If there is a loop around $v$ with a masked flat quarter of the
first type, then the context is $(n,k)=(3,0)$ or $(4,1)$.  In the
case of context $(4,1)$, we have
    $$\pi(v,D) \le \piV(|v|).$$

The constants $\pi(v,D)$ are called penalties.  When we
approximate $\sigma_R(D)$ by means of an upper bound on
$\pi(v,D)$, the upright diagonal is said to be {\it erased}.

\subsection{Anchored Simplices}

Often, $\pi(v,D)\le0$.  When this is the case, we generally
eliminate (that is, erase) the upright diagonal and all the
anchored simplices around it. (This has been done implicitly
already.)  When we do this, we use the function $\vor_0$ on the
region that becomes ``exposed'' in the process.  It is possible to
eliminate the upright diagonals in all cases, except
    \begin{enumerate}
        \item loops,
        \item a configuration $\Splus$,
        \item a configuration $\Sminus$,
        \item a configuration $\Sfour$.
    \end{enumerate}

\subsubsection{The configuration $\Splus$}

The configuration $\Splus$ around an upright diagonal $(0,v)$
consists of two upright quarters, one other anchored simplex, and
one large gap. The configuration $\Splus$ never masks a flat
quarter.  If the other anchored simplex is not an upright quarter,
then the penalty satisfies
    $$\pi(v,D)\le 0.00222.$$
If all three anchored simplices are upright quarters, then the
penalty satisfies
    $$\pi(v,D)\le 0.008.$$

\subsubsection{The configuration $\Sminus$}

The configuration $\Sminus$ around an upright diagonal $(0,v)$
consists of three upright quarters that all lie in a common
half-space bounded by a plane through $(0,v)$.   If the
configuration masks a flat quarter, then the penalty satisfies
    $$\pi(v,D)\le \piG(|v|) + 2\piV(|v|).$$
If it does not mask a flat quarter, then
    $$\pi(v,D)\le 3\piG(|v|).$$
There is at most one $\Sminus$ configuration in a contravening
decomposition star.  Each upright quarter $Q$ is assigned a score
$\sigma(Q)$.  The three quarters satisfy
    $$\begin{array}{lll}
        \sum_{(3)} \sigma(Q) &< -0.4339,\\
        \sum_{(3)} \tau(Q) &> 0.5606.
    \end{array}
    $$

\subsubsection{The configuration $\Sfour$}

The configuration $\Sfour$ around an upright diagonal $(0,v)$
consists of four upright quarters, and a large gap.  The
configuration masks at most $k\le 2$ flat quarters.  The penalty
satisfies
    $$\pi(v,D) \le 2k \piV(|v|) + (4-2k)\piG(|v|).$$
There is at most one $\Sfour$ configuration in a contravening
decomposition star.  It does not appear in a contravening
decomposition star with $\Sminus$.  Each upright quarter $Q$ is
assigned a score $\sigma(Q)$.  The four quarters satisfy
    $$\begin{array}{lll}
        \sum_{(4)} \sigma(Q) &< -0.25,\\
        \sum_{(4)} \tau(Q) &> 0.4.
    \end{array}
    $$

\subsection{Compatibility Notes}

We conclude this section with a few comments about differences in
notation between this article and the others.

In \cite{part4}, simplices of types $A$, $B$, $C$, and $C'$ are
called $S_A$, $S_B$, $S_C$, and $S_D$, respectively.  The
functions $\sigma_R(D)$ and $\tau_R(D)$ are denoted $\sigma(R)$
and $\tau(R)$, respectively, in other articles. The truncations
$\vor_{0,R}(D)$ and $\tau_{0,R}(D)$ are denoted $\vor_0(R)$ and
$\tau_0(R)$, respectively, in other articles.

\newpage
\hbox{}
\newpage

\section{Appendix. Calculations}
\label{app:calculation}


The proof of the Kepler Conjecture relies on a large number of
inequalities that have been established by interval arithmetic
with a computer. Generally, speaking there is no geometry of
interest associated with these inequalities.   The methods that
are used in verifying these conjecture are described in
\cite{part1} and \cite{algorithm}.  (The paper describing the
algorithms goes beyond what is used in these verifications,
because it also proposes some new unimplemented techniques that
should simplify the verification process.)

Interval calculations are arranged according to the section in
which they appear. Each inequality is accompanied by one or more
reference numbers. These identification numbers are needed to find
further details about the calculation in \cite{part1}. Most of the
verifications were completed by Samuel Ferguson. \footnote{I
express my sincere thanks to Samuel Ferguson for his fundamental
contributions to this project.} The verifications that he made are
marked with a dagger (\dag). The computer-based proofs are carried
out by interval arithmetic as in our earlier papers.  A few
inequalities are not standard in the sense that they have
substantially higher dimension that the others. Linear programming
methods are used to break these problems into six-dimensional
pieces according to the suggestion of \cite[Appendix 2]{part3}.

  Edge lengths whose bounds are not specified are assumed
to be between 2 and $2t_0$.  The first edge of an upright quarter
is its diagonal.  The fourth edge of a flat quarter is its
diagonal.

$\hat\sigma$ is the function of \cite[Section 3.11]{part4} and
$\hat\tau=\sol\zeta\pt-\hat\sigma$.  The function definitions are
those appearing in Appendix \ref{app:constants}.  The other
functions are described in Appendix \ref{app:constants}.

\begin{figure}[htb]
  \centering
  \includegraphics{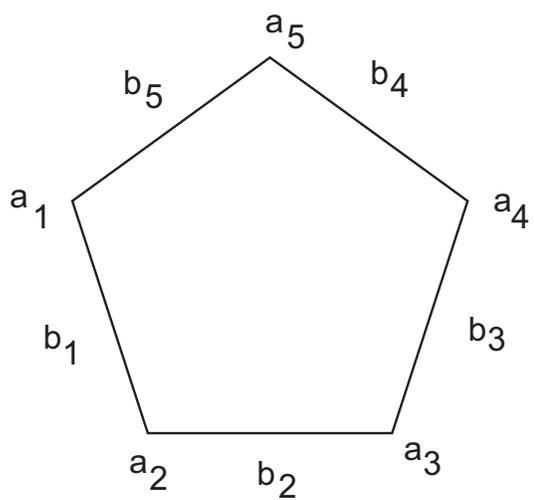}
  \caption{The labels on an $n$-gon.}
  \label{fig:labeledpent}
\end{figure}

We introduce some notation for the heights and edge lengths of a
polygon.  The heights will generally be between $2$ and $2t_0$,
the edge lengths between consecutive corners will generally be
$2$, $2t_0$, or $2\sqrt{2}$.  We represent the edge lengths by a
vector
    $$(a_1,b_1,a_2,b_2,\ldots,a_n,b_n),$$
if the corners of an $n$-gon, ordered cyclically have heights
$a_i$ and if the edge length between corner $i$ and $i+1$ is
$b_i$.  Figure \ref{fig:labeledpent} illustrates the case $n=5$.
We say two vectors are equivalent if they are related by a
different cyclic ordering on the corners of the polygon, that is,
by the action of the dihedral group.

\vfill

\parindent=0pt

\newpage\onecolumn

\subsection{Group 1}
\label{app:group1} \oldlabel{A.3.1}

$\tau - 0.2529\dih > -0.3442$,
    if $y_1\in[2.3,2t_0]$, and $\dih\ge1.51$.
    \refno{572068135}

$\tau_0  - 0.2529\dih > -0.1787$,
    if $y_1\in[2.3,2t_0]$, $y_6\in[2\sqrt2,3.02]$, $1.26\le\dih\le 1.63$.
    \refno{723700608}

$\hat\tau- 0.2529\dih_2 > -0.2137$,
    if $y_2\in[2.3,2t_0]$, $y_4\in[2t_0,2\sqrt2]$,
    \refno{560470084}

$\tau_0 - 0.2529\dih > -0.1371$,
    if $y_1\in[2.3,2t_0]$, $y_5,y_6\in[2t_0,3.02]$,
    $1.14\le\dih\le 1.51$.
    \refno{535502975}

\subsection{  Group 2\dag} 
\oldlabel{A.3.8} \label{A.3.8}

A.\quad $\dih < 1.63$, if $y_6\ge2t_0$, $y_2,y_3\in[2,2.168]$.
    \refno{821707685}

B.\quad $\dih< 1.51$, if $y_5=2t_0$, $y_6\ge 2t_0$,
$y_2,y_3\in[2,2.168]$.
    \refno{115383627}

C.\quad $\dih < 1.93$, if $y_6\ge2t_0$, $y_4=2\sqrt2$,
$y_2,y_3\in[2,2.168]$.
    \refno{576221766}

D.\quad $\dih < 1.77$, if $y_5=2t_0$, $y_6\ge 2t_0$,
$y_4=2\sqrt2$,
    $y_2,y_3\in[2,2.168]$.
    \refno{122081309}

$\tau_0 -0.2529\dih > -0.2391$, if $y_6\ge2t_0$, $\dih\ge1.2$,
    $y_2,y_3\in[2,2.168]$.
    \refno{644534985}

$\tau_0 -0.2529\dih > -0.1376$, if $y_5=2t_0$, $y_6\ge2t_0$,
$\dih\ge1.2$,
    and $y_2,y_3\in[2,2.168]$.
    \refno{467530297}

$\tau_0 -0.2529\dih > -0.266$, if $y_6\ge2t_0$,
$y_4\in[2t_0,2\sqrt2]$,
    $\dih\ge1.2$,
    $y_2,y_3\in[2,2.168]$.
    \refno{603910880}

$\tau_0 -0.2529\dih > -0.12$, if $y_5=2t_0$,
    $y_6\ge2t_0$, $y_4\in[2t_0,2\sqrt2]$,
    $\dih\ge1.2$,
    $y_2,y_3\in[2,2.168]$.
    \refno{135427691}

$\dih<1.16$, if $y_5=2t_0$, $y_6\ge 2t_0$, $y_4=2$,
$y_2,y_3\in[2,2.168]$.
    \refno{60314528}

$\tau_0 -0.2529\dih > -0.1453$, if $y_2,y_3\in[2,2.168]$,
    $y_5\in[2t_0,3.488]$, $y_6=2t_0$.
    \refno{312132053}

\subsection{  Group 3 }
\label{A.3.8bis} \oldlabel{A.3.8, Case 2-b}

$\dih_2 > 0.74$, if $y_1\in[2t_0,2.696]$,
    $y_2,y_3\in[2,2.168]$.
    \refno{751442360}

$\tau_0 -0.2529\dih > -0.2391$, if
     $\Delta(y_5^2,4,4,8,2t_0^2,y_6^2)\ge0$,
    $y_2,y_3\in[2,2.168]$,
    $y_5\in[2t_0,3.488]$.
    \refno{893059266}

$\dih + 0.5 (2.402-y_4) < \pi/2$, if $y_5\ge 2t_0$,
    $y_2,y_3\in[2,2.168]$.
    \refno{690646028}

\subsection{  Group 4} 
\label{app:178} \oldlabel{A.3.9}

$\dih>1.78$, if
    $y_4=3.2$, $y_1\in[2t_0,2\sqrt2]$, $y_2+y_3\le4.6$.
    \refno{161665083}

\subsection{  Group 5} 
\oldlabel{A.4.4.1}

The following inequalities hold for flat quarters. In these
inequalities the fourth edge is the diagonal.

    $$
    \begin{array}{lll}
    - \dih_2& + 0.35 y_2 - 0.15 y_1 - 0.15 y_3 + 0.7022 y_5 - 0.17 y_4 > -0.0123,\\
    \dih_2& - 0.13 y_2 + 0.631 y_1 +
        0.31 y_3 - 0.58 y_5 + 0.413 y_4 + 0.025 y_6 > 2.63363,\\
    -\dih_1& + 0.714 y_1 - 0.221 y_2 - 0.221 y_3 +
        0.92 y_4 - 0.221 y_5 - 0.221 y_6 > 0.3482,\\
    \dih_1& - 0.315 y_1 + 0.3972 y_2 + 0.3972 y_3 - \\
       &\quad\quad 0.715 y_4 +  0.3972 y_5 + 0.3972 y_6 > 2.37095,\\
    - \sol& - 0.187 y_1 - 0.187 y_2 -
       0.187 y_3 + 0.1185 y_4 + 0.479 y_5 + 0.479 y_6 > 0.437235\dag,\\
    \sol& + 0.488 y_1 + 0.488 y_2 +
       0.488 y_3 - 0.334 y_5 - 0.334 y_6 > 2.244\dag,\\
    - \hat\sigma& - 0.145 y_1 - 0.081 y_2 - 0.081 y_3 -
        0.133 y_5 - 0.133 y_6 > -1.17401,\\
    -\hat\sigma& -0.12 y_1 -0.081 y_2 -0.081 y_3 -0.113 y_5 -0.113
        y_6+
        0.029 y_4 > -0.94903,\\
    \hat\sigma& + 0.153 y_4 + 0.153 y_5 + 0.153 y_6 < 1.05382,\\
    \hat\sigma& +0.419351 \sol +0.19 y_1+ 0.19 y_2 + 0.19 y_3 < 1.449,\\
    \hat\sigma& +0.419351 \sol < -0.01465 + 0.0436 y_5 + 0.0436 y_6
        +0.079431 \dih,\\
    \hat\sigma& < 0.0114,\\
    \hat\tau& > 1.019\,\pt,\\
    \end{array}
    $$
    \oldlabel{4.4.1}
\refno{867513567}

\subsection{  Group 6} 
\label{app:nugamma1}\oldlabel{A.4.5}

In a quadrilateral cluster, with a given edge $y_4$ as the
diagonal, the other diagonal will be denoted $y'_4$.

The following relations for upright quarters (scored by $\nu$)
hold.  (We use the inequalities of \cite[Section 4]{part3} for
upright quarters in quad clusters, which are scored by a different
function.) In these inequalities the upright diagonal is the first
edge. We include in this group, the inequalities
\cite[$\A_2$]{part4}, \cite[$\A_3$]{part4} for upright quarters.

$$
\begin{array}{lll}
y_1&>2t_0,\\
y_1&<2\sqrt2,\\
 \dih_1& - 0.636 y_1 + 0.462 y_2 + 0.462 y_3 - 0.82 y_4 + 0.462 y_5 +
 0.462 y_6 > 1.82419,\\
 - \dih_1& + 0.55 y_1 - 0.214 y_2 - 0.214 y_3 + 1.24 y_4 - 0.214 y_5
 - 0.214 y_6 > 0.75281,\\
 \dih_2& + 0.4 y_1 - 0.15 y_2 + 0.09 y_3 + 0.631 y_4 - 0.57 y_5 + 0.23 y_6
  > 2.5481,\\
 - \dih_2& - 0.454 y_1 + 0.34 y_2 + 0.154 y_3 - 0.346 y_4 +
0.805 y_5 > -0.3429,\\
 \dih_3& + 0.4 y_1 - 0.15 y_3 + 0.09 y_2 + 0.631 y_4 - 0.57 y_6 + 0.23 y_5
  > 2.5481,\\
 - \dih_3& - 0.454 y_1 + 0.34 y_3 + 0.154 y_2 - 0.346 y_4 +
0.805 y_6 > -0.3429,\\
 \sol& + 0.065 y_2 + 0.065 y_3 + 0.061 y_4 - 0.115 y_5 -
0.115 y_6 > 0.2618,\\
 - \sol& - 0.293 y_1 - 0.03 y_2 - 0.03 y_3 + 0.12 y_4 +
0.325 y_5 + 0.325 y_6 > 0.2514,\\
 -\nu& - 0.0538 y_2 - 0.0538 y_3 -0.083 y_4 - 0.0538 y_5 -
0.0538 y_6 > -0.5995,\\
 \nu& \le 0,\quad
 \text{(see \cite[3.13.3]{formulation} and \cite[3.13.4]{formulation})}\\
 \tau_\nu& - 0.5945 \,\pt > 0.\\  
\end{array}
$$
\refno{498839271} \oldlabel{A.4.5.1}

$$
\begin{array}{lll}
    \nu& -4.10113 \dih_1< -4.3223,\\
    \nu& -0.80449 \dih_1< -0.9871,\\
    \nu& -0.70186 \dih_1< -0.8756,\\
    \nu& -0.24573 \dih_1< -0.3404,\\
    \nu& -0.00154 \dih_1< -0.0024,\\
    \nu& +0.07611 \dih_1<  0.1196.\\
\end{array}
$$
\refno{\cite[$\A_2$]{part4}} \oldlabel{A.4.5.2}

$$
\begin{array}{lll}
    \tau_\nu& +4.16523 \dih_1>  4.42873,\\
    \tau_\nu& +0.78701 \dih_1>  1.01104,\\
    \tau_\nu& +0.77627 \dih_1>  0.99937,\\
    \tau_\nu& +0.21916 \dih_1>  0.34877,\\
    \tau_\nu& +0.05107 \dih_1>  0.11434,\\
    \tau_\nu& -0.07106 \dih_1> -0.07749.\\
\end{array}
$$
\refno{\cite[$\A_3$]{part4}} \oldlabel{A.4.5.3}

\subsection{  Group 7} 
\label{app:nugamma}\oldlabel{A.4.5.4}

 The following additional inequalities
are known to hold if the upright diagonal has height at most
$2.696$. $\nu_\Gamma$ denotes the restriction of $\nu$ to a
simplex of compression type.

\begin{eqnarray}
y_1& < 2.696,\\
\dih_1& - 0.49 y_1 + 0.44 y_2 + 0.44 y_3 - 0.82 y_4 + 0.44 y_5 + 0.44 y_6 > 2.0421,\\
-\dih_1& + 0.495 y_1 - 0.214 y_2 - 0.214 y_3 + 1.05 y_4 - 0.214
y_5 -
 0.214 y_6 > 0.2282,\label{eqn:label}\\ 
\dih_2& + 0.38 y_1 - 0.15 y_2 + 0.09 y_3 + 0.54 y_4 - 0.57 y_5 + 0.24 y_6 > 2.3398,\\
-\dih_2& - 0.375 y_1 + 0.33 y_2 + 0.11 y_3 - 0.36 y_4 + 0.72 y_5 +
    0.034 y_6 > -0.36135,\\
\sol&+ 0.42 y_1 + 0.165 y_2 + 0.165 y_3 - 0.06 y_4 - 0.135 y_5 -
    0.135 y_6 > 1.479,\\
-\sol& - 0.265 y_1 - 0.06 y_2 - 0.06 y_3 + 0.124 y_4 + 0.296 y_5 +
    0.296 y_6 > 0.0997,\\
-\nu& + 0.112 y_1 - 0.142 y_2 - 0.142 y_3 - 0.16 y_4 -
    0.074 y_5 - 0.074 y_6 > -0.9029,\\
\nu& +0.07611 \dih_1<  0.11,\\
\nu_\Gamma& - 0.015 y_1 -0.16 (y_2+y_3+y_4)-0.0738 (y_5+y_6)
    > -1.29285,\\
\tau_\nu& -0.07106 \dih_1> -0.06429,\\
\tau_\nu& > 0.0414.
\end{eqnarray}
\refno{319046543} \oldlabel{A.4.5.4}

\begin{remark} In connection with the Inequality (\ref {eqn:label}), we
occasionally use the stronger constant $0.2345$ instead of
$0.2282$.  To justify this constant, we have checked using
interval arithmetic that the bound $0.2345$ holds if $y_1\le2.68$
or $y_4\le2.475$. Further interval calculations show that the
anchored simplices can be erased if they share an upright diagonal
with such a quarter.
\end{remark}

\bigskip

\subsection{  Group 8\dag} 
\oldlabel{A.4.8}
 \label{app:dih}


We give lower and upper bounds on  dihedral angles.  The domains
that we list are not disjoint. In general we consider an edge as
belonging to the most restrictive domain that the information of
the following charts permit us to conclude that it lies in.

The following chart summarizes the bounds.  The dihedral angle is
computed along the first edge.  The chart is divided into three
sections.  In the first, there is no upright diagonal. In the
second, the upright diagonal is the first edge. In the third, the
upright diagonal is the third edge. The bounds in the second
section have been established in \cite[$\A_8$]{part4}.

In the first group $y_1,y_2,y_3\in[2,2t_0]$. In the third row, the
dihedral bound $1.624$ holds for $y_6$ on the larger interval
$[2t_0,3.02]$.  In the seventh row, the dihedral bound $1.507$
holds for $y_5,y_6$ on the larger interval $[2t_0,3.02]$.
$$
    \begin{array}{llllll}
    y_5  & y_6 & y_4 & \dih_{\mn} & \dih_{\mx} \\
    \left[ 2 , 2  t_0   \right]   &  \left[ 2 , 2 t_0\right] & [2t_0,2\sqrt2] & 1.153 & 2.28\\
    \left[2,2t_0\right] & [2,2t_0] & \ge2\sqrt2 & 1.32 & 2\pi\\
%
    \left[2,2t_0\right]&[2t_0,2\sqrt2]&[2,2t_0]&0.633&1.624\\
    \left[2,2t_0\right]&[2t_0,2\sqrt2]&[2t_0,2\sqrt2]&1.033&1.929\\
    \left[2,2t_0\right]&[2t_0,2\sqrt2]&\ge2t_0&1.033&2\pi\\
    \left[2,2t_0\right]&[2t_0,2\sqrt2]&\ge2\sqrt2&1.259&2\pi\\
%
    \left[2t_0,2\sqrt2\right]&[2t_0,2\sqrt2]&[2,2t_0]&0.817&1.507\\
    \left[2t_0,2\sqrt2\right]&[2t_0,2\sqrt2]&[2t_0,2\sqrt2]&1.07&1.761\\
    \left[2t_0,2\sqrt2\right]&[2t_0,2\sqrt2]&\ge2t_0&1.07&2\pi\\
    \left[2t_0,2\sqrt2\right]&[2t_0,2\sqrt2]&\ge2\sqrt2&1.23&2\pi\\
\end{array}
$$

In the second group $y_1\in[2t_0,2\sqrt2]$, $y_2,y_3\in[2,2t_0]$.
$$
\begin{array}{lllll}
y_5  & y_6 & y_4 & \dih_{\min} & \dih_{\max} \\
%
    \left[2,2t_0\right]&[2,2t_0]&[2,2t_0]&0.956&2.184\\
    \left[2,2t_0\right]&[2,2t_0]&[2t_0,2\sqrt2]&1.23&\pi\\
    \left[2,2t_0\right]&[2,2t_0]&\ge2t_0&1.23&\pi\\
    \left[2,2t_0\right]&[2,2t_0]&\ge2\sqrt2&1.416&\pi\\
\end{array}
$$

In the third group $y_1,y_2\in[2,2t_0]$, $y_3\in[2t_0,2\sqrt2]$.
$$
\begin{array}{lllll}
y_5  & y_6 & y_4 & \dih_{\min} & \dih_{\max} \\
%
    \left[2,2t_0\right]&[2,2t_0]&[2,2t_0]&0.633&1.624\\
    \left[2,2t_0\right]&[2,2t_0]&\ge2t_0&1.033&2\pi\\
    \left[2,2t_0\right]&[2t_0,2\sqrt2]&[2,2t_0]&0&1.381\\
    \left[2,2t_0\right]&[2t_0,2\sqrt2]&\ge2t_0&0.777&2\pi\\
\end{array}
$$
\refno{853728973}

\subsection{  Group 9\dag} 
\label{app:9} \oldlabel{A.4.9}

(The verifications in this section that do not involve the score
were verified by S. Ferguson.)

Depending on the lengths of the edges $y_5$, $y_6$, $y_4$, there
are additional inequalities that hold. Generally these
inequalities are for a fragment of a subregion. If $v_2$, $v_1$,
$v_3$ are consecutive corners, we form a simplex
$S(y_1,y_2,y_3,y_4,y_5,y_6) = (0,v_1,v_2,v_3)$.  There is no score
intrinsically associated with this simplex, unless the subregion
is triangular.  But we can state various useful bounds on the
dihedral angles along the first edge $(0,v_1)$. A triple $(a,b,c)$
preceding the inequality gives bounds on the edges
$(y_5,y_6,y_4)$, respectively.

In this group of inequalities, $y_1,y_2,y_3\in[2,2t_0]$.  Set
$\vor_x=\vor$ if the simplex has type $A$, and set
$\vor_x=\vor_0$, otherwise.
$$\begin{array}{lll}
&(y_5,y_6,y_4)\\
&([2,2t_0], [2,2t_0], \ge2\sqrt2)\\
&\quad\quad\dih-0.372 y_1 +0.465 y_2 +0.465 y_3 + 0.465 y_5 + 0.465 y_6 >4.885,\\
&([2,2t_0], [2t_0,2\sqrt2], [2t_0,2\sqrt2])\\
&\quad\quad0.291 y_1 -0.393 y_2 -0.586 y_3 +0.79 y_4 -0.321 y_5
-0.397 y_6 -\dih
    <  -2.47277,\\
&([2,2t_0], [2t_0,2\sqrt2], \ge2t_0)\\
&\quad\quad0.291 y_1 -0.393 y_2 -0.586 y_3  -0.321 y_5 -0.397 y_6
-\dih
    <  -4.45567,\\
&([2,2t_0], [2t_0,2\sqrt2], \ge2\sqrt2)\\
&\quad\quad0.291 y_1 -0.393 y_2 -0.586 y_3  -0.321 y_5 -0.397 y_6
-\dih
    <  -4.71107,\\
&\quad\quad
\dih-0.214 y_1 +0.4 y_2 +0.58 y_3 +0.155 y_5 + 0.395 y_6 > 4.52345,\\
&([2t_0,2\sqrt2], [2t_0,2\sqrt2], [2,2t_0])\\
&\quad\quad\tau_x  > D(3,2),\quad
\text{see \cite[$\A_{16}$]{part4}}\\
    &\quad\quad \vor_x < Z(3,2),\quad
    \text{ see \cite[$\A_{16}$]{part4}}\\
    &\quad\quad - \sol -0.492 y_1 - 0.492 y_2 -0.492 y_3 +0.43 y_4 +\\
    &\quad\quad\quad0.038 y_5+0.038 y_6 < -2.71884,\\
    &\quad\quad - \vor_x -0.058 y_1 -0.105 y_2  -0.105 y_3 -0.115 y_4 -\\
    &\quad\quad\quad0.062 y_5 -0.062 y_6 > -1.02014,\\
    &\quad\quad \vor_0 +0.419351 \sol < 0.3085,\\
    &\quad\quad 0.115 y_1 -0.452 y_2  -0.452 y_3 +0.613 y_4 -0.15 y_5 -0.15 y_6 -\dih  < -2.177,\\
&([2t_0,2\sqrt2], [2t_0,2\sqrt2], [2t_0,2\sqrt2])\\
&\quad\quad0.115 y_1 -0.452 y_2 -0.452 y_3 +0.618 y_4 -0.15 y_5 -0.15 y_6 -\dih  < -2.17382,\\
&\quad\quad \vor_0< -0.121,\quad\text{see \cite[$\A_{18}$]{part4}}\\
&\quad\quad\tau_0 > D(3,3)=0.21301,\quad \text{see \cite[$\A_{17}$]{part4}}\\
&([2t_0,2\sqrt2], [2t_0,2\sqrt2], \ge2t_0)\\
&\quad\quad0.115 y_1 -0.452 y_2 -0.452 y_3  -0.15 y_5 -0.15 y_6 -\dih < -3.725,\\
&([2t_0,2\sqrt2], [2t_0,2\sqrt2], \ge2\sqrt2)\\
&\quad\quad0.115 y_1 -0.452 y_2 -0.452 y_3  -0.15 y_5 -0.15 y_6 -\dih < -3.927.\\
%
\end{array}
$$\refno{529738375}
\oldlabel{A.4.9.1}

\subsection{  Group 10} 
\label{app:10}

In the next group $y_1\in[2t_0,2\sqrt2]$, $y_2,y_3\in[2,2t_0]$.
Set $\vor_x=\vor$, if the simplex has type $C$ or $C'$.  Set
$\vor_x=\vor_0$, otherwise.  Other than this, retain the notation
and conventions from Appendix \ref{app:9}.
$$
\begin{array}{lll}
&(y_5,y_6,y_4)\\
&([2,2t_0], [2,2t_0], [2t_0,2\sqrt2])\\
&\quad\quad\vor_x < 0,\quad \text{(cf. Appendix \ref{app:455})}\\
&\quad\quad 0.47 y_1 -0.522 y_2 -0.522 y_3 +0.812 y_4
        -0.522 y_5 -0.522 y_6 -\dih  < -2.82998,\\
&([2,2t_0], [2,2t_0], \ge2t_0)\\
&\quad\quad0.47 y_1 -0.522 y_2 -0.522 y_3  -0.522 y_5 -0.522 y_6
-\dih
        < -4.8681,\\
&([2,2t_0], [2,2t_0], \ge2\sqrt2)\\
&\quad\quad0.47 y_1 -0.522 y_2 -0.522 y_3  -0.522 y_5  -0.522 y_6
-\dih
        < -5.1623.\\
\end{array}
$$\refno{456320257}
\oldlabel{A.4.9.2}

\subsection{  Group 11} 

In the next group $y_1,y_2\in[2,2t_0]$, $y_3\in[2t_0,2\sqrt2]$.
Otherwise, retain the notation and conventions of Appendix
\ref{app:9}.
$$
\begin{array}{lll}
&(y_5,y_6,y_4)\\
&([2,2t_0], [2,2t_0], \ge2t_0)\\
&\quad\quad-0.4 y_3 +0.15 y_1 -0.09 y_2 -0.631 y_6-0.23 y_5-\dih < -3.9788,\\
&([2,2t_0], [2t_0,2\sqrt2], [2,2t_0])\\
&\quad\quad0.289 y_1 -0.148 y_2 -1.36 y_3 +0.688 y_4
    -0.148 y_5 -1.36 y_6 -\dih  < -6.3282,\\
&([2,2t_0], [2t_0,2\sqrt2], \ge2t_0)\\
&\quad\quad0.289 y_1 -0.148 y_2 -0.723 y_3 -0.148 y_5 -0.723 y_6
-\dih
    < -4.85746.\\
\end{array}
$$\refno{664959245}
\oldlabel{A.4.9.3}

\subsection{  Group 12} 
\oldlabel{A.4.12} \label{app:pent}

\begin{enumerate}
\item
    $$\nu < -0.055\quad\text{and }\tau_\nu>0.092,
    $$
    \oldlabel{A.4.12.1}
if $y_1\in[2.696,2\sqrt2]$, $y_2,y_6\in[2.45,2t_0]$.
\refno{704795925}

\item
    $$
    \hat\sigma < -0.039\quad\text{and }\hat\tau>0.094,
    $$
    \oldlabel{A.4.12.2}
if $y_2\in[2.45,2t_0]$, $y_4\in[2t_0,2\sqrt2]$. \refno{332919646}

\item
    $$\vor < -0.197\quad\text{and }\tau_V >0.239,
    $$
    \oldlabel{A.4.12.3}
if $y_1\in[2.696,2\sqrt2]$, $y_2\in[2.45,2t_0]$,
$y_6\in[2.45,2t_0]$, and $y_4\in[2t_0,2\sqrt2]$. (These simplices
have type $S'_C$ in the sense of \cite[2.10]{part4}, and this
justifies the use of the function $\vor$.) \refno{335795137}

\item
    $$\vor_0 < -0.089\quad\text{and }\tau_0 > 0.154,
    $$
    \oldlabel{A.4.12.6}
if $y_1\in [2.45,2t_0]$, $y_5,y_6\in[2t_0,2\sqrt2]$.
\refno{605071818\dag}

\item
    $$\vor_0 < -0.089\quad\text{and }\tau_0 > 0.154,
    $$
    \oldlabel{A.4.12.7}
if $y_1\in [2.45,2t_0]$, $y_4,y_5\in[2t_0,2\sqrt2]$.
\refno{642806938\dag}
\end{enumerate}

\subsection{  Group 13} 
\label{ineq:1.678}\oldlabel{A.6}


$\octavor < \octavor_0 - 0.017$, if $y_1\in[2t_0,2.696]$, and
        $\eta_{126}\ge\sqrt2$.
\refno{104506452}

    $$
    \begin{array}{lll}
    \dih>1.678, \text { if } y_4\ge 3.0,  \text{ and }
        y_2+y_3+y_5+y_6\le8.77,
    \end{array}
    $$
\refno{601083647}

\subsection{Group 14} 
\oldlabel{A.6.3}


$\Gamma < 0.3138-0.157 y_5$,  if  $y_5\in [2,2.138]$,
    $y_4\in[2t_0,2.6]$,
\refno{543730647}

$\Gamma < -0.06$,  if  $y_2\in[2.121,2.145]$,
    $y_4\in[2t_0,2\sqrt2]$, $y_5\in[2.22,2.238]$.
\refno{163030624}

$\Gamma < 10^{-6}+ 1.4-0.1 y_1 -0.15 (y_2+y_3+y_5+y_6)$,
     if $y_4\in[2,2\sqrt2]$, $y_1,y_2,y_3\in[2,2.2]$,
    $y_5,y_6\in[2,2.35]$,
        \refno{181462710}


\subsection{  Group 15} 
\label{app:A.7.1}

$\vor < \vor_0$ if $y_4\in[2.7,2\sqrt2]$,
$y_1,y_2,y_3\in[2,2.14]$. \refno{463544803}

$\vor _0 < -0.064$, if $y_4 \in[2t_0,2.72]$,
$\eta_{456}\ge\sqrt2$. \refno{399326202}

$\vor_0 < 1.0612 - 0.08(y_1+y_2+y_3) - 0.142(y_5+y_6)$, if $y_4
\in[2.7,2\sqrt2]$. \refno{569240360}

$\vor_0 < -0.0713$, if $y_4\in[2.59,2.64]$, $y_5\in[2.47,2t_0]$,
$y_6\in[2.1,3.51]$. \refno{252231882}

$\vor_0 < -0.06$, if $y_1,y_2,y_3\in[2,2.13]$, $y_4\in[2.7,2.74]$,
$\eta_{456}\ge\sqrt2$. \refno{472436131}

$\vor_0 < -0.058$, if $y_4\in[2t_0,2.747]$, $\eta_{456}\ge\sqrt2$.
\refno{913534858}

$\vor_0 < -0.0498$, if $y_4\in[2t_0,2.77]$, $\eta_{456}\ge\sqrt2$.
\refno{850226792}

\subsection{  Group 16} 
\label{app:A.7.2}

We assume that $y_1,y_2,y_3\le2.14$ in the following inequalities.

$-\Gamma -0.145 y_1 -0.08 (y_2+y_3) -0.133 (y_5+y_6) > -1.146$,
    if $y_4\in[2t_0,2\sqrt2]$.
\refno{594246986}

$-\Gamma -0.145 y_1 -0.081 (y_2+y_3) -0.16 (y_5+y_6) > -1.255$,
    if $y_4\in[2t_0,2\sqrt2]$, and $y_5,y_6\in[2,2.3]$.
\refno{381970727}

$-\Gamma -0.03 y_1 -0.03 (y_2+y_3) -0.094 (y_5+y_6) > -0.5361$.
    if $y_4\in[2t_0,2\sqrt2]$, $y_5+y_6\ge4.3$.
\refno{951798877}

$-\Gamma -0.03 y_1 -0.03 (y_2+y_3) -0.16 (y_5+y_6) >
-0.82-10^{-6}$.
    if $y_4\in[2t_0,2\sqrt2]$, $y_5+y_6\le4.3$.
\refno{923397705}

$\Gamma  < -0.053$ if $y_4\in[2t_0,2\sqrt2]$, $y_5\ge2.35$.
\refno{312481617}

$\Gamma < -0.041$, if $y_4\in[2t_0,2\sqrt2]$, $y_5\in[2.25,2t_0]$.
\refno{292760488}

$\Gamma +0.419351 \sol < 0.079431\dih + 0.0436 (y_5+y_6) -
0.0294$, if $\eta_{456}\le \sqrt2$, $y_4\in[2t_0,2\sqrt2]$.
\refno{155008179}

$\Gamma < 1.1457 - 0.1 (y_1+y_2+y_3) - 0.17 y_5 - 0.11 y_6$,
    if $y_1,y_2,y_3\in[2,2.13]$, $y_5\in[2,2.1]$,
    $y_6\in[2.27,2.43]$, $y_4\in[2t_0,2.67]$.
\refno{819450002}

$1.69 y_4 + y_5 + y_6 > 9.0659$, if $y_4\in[2t_0,2.7]$
    and $\eta_{456}\ge\sqrt2$.
\refno{495568072}

$1.69 y_4 + y_5 + y_6 > 9.044$, if $y_4\in[2t_0,2.77]$
    and $\eta_{456}\ge\sqrt2$.
\refno{838887715}

$y_5+y_6 > 4.4$, if $\eta_{456}\ge\sqrt2$, and $y_4\le2.72$.
\refno{794413343}

\subsection{  Group 17} 
\label{app:17} \oldlabel{A.7.3}

In these inequalities $y_5,y_6\in[2t_0,2\sqrt2]$ and
$y_4\in[2,2t_0]$. We distinguish three cases:
    \begin{enumerate}
    \item The simplex has
        type $A$, so that $y_5,y_6\in[2t_0,2.77]$,
    and the scoring is $\vor$, the analytic Voronoi function.
    \item $y_6 \in [2.77,2\sqrt2]$.
    \item The edges have lengths $y_5,y_6\in[2t_0,2.77]$, and
    $\eta_{456}\ge\sqrt2$. In the last two cases, the scoring is by
    the truncated Voronoi function $\vor_0$.
    \end{enumerate}

$-\vor -0.058 y_1 - 0.08 y_2 - 0.08 y_3 - 0.16 y_4
        -0.21 (y_5+y_6) > -1.7531$,
        if $y_1,y_2,y_3\le2.14$ and $y_5,y_6\in[2t_0,2.77]$.
\refno{378020227}

$-\vor_0 - 0.058 y_1 - 0.1 y_2 - 0.1 y_3 - 0.165 y_4
        -0.115 y_6 -0.12 y_5 > -1.38875$,
        if $\eta_{456}\ge\sqrt2$,
         $y_6\in[2.77,2\sqrt2]$, and
        $y_1,y_2,y_3\le2.14$.
\refno{256893386}

$y_4+y_5+y_6> 7.206$,
    if $y_5,y_6\in[2t_0,2.77]$ and
        $\eta_{456}\ge\sqrt2$. \refno{749955642}

$-\vor_0 - 0.058 y_1 - 0.05 y_2 - 0.05 y_3 - 0.16 y_4
        -0.13 y_6 -0.13 y_5 > -1.24547$,
        if $y_1,y_2,y_3\le2.14$, $y_5,y_6\in[2t_0,2.77]$, and
        $\eta_{456}\ge\sqrt2$.
\refno{653849975}

$\vor_0 < -0.077$, if $y_4\in[2t_0,2\sqrt2]$,
$y_5\in[2.77,2\sqrt2]$. \refno{480930831}

$\vor +0.419351 \sol < 0.289$, if $y_4,y_5\in[2t_0,2.77]$.
\refno{271703736}

$\vor_0 < 1.798 -0.1 (y_1+y_2+y_3) - 0.19 y_4 - 0.17 (y_5+y_6)$,
if $y_5,y_6\in[2.7,2.77]$. \refno{900212351}

\subsection{  Group 18 } 
\label{app:A.7.4}

The next inequality is used as an estimate when there is an
upright diagonal enclosed over a simplex of the third type in
Appendix \ref{app:17}.

$\vor< -0.078/2$, $y_1\in[2t_0,2\sqrt2]$, $y_4\in[2t_0,2.6961]$
for simplices of type $C$. \refno{455329491}

If there is no enclosed upright diagonal, then we can sometimes
use the following instead.

$\vor(S,\sqrt2) < -0.078$, if $y_5,y_6\in[2t_0,2.6961]$,
$\eta_{456}\ge\sqrt2$. \refno{857241493}

\subsection{  Group 19 } 
\oldlabel{A.2.3}


The interval calculations here show that the set of separated
vertices (\ref{definition:admissible:excess}) can be generalized
to include   opposite vertices of a quadrilateral unless the edge
between those vertices forms a flat quarter.   Consider a vertex
of type $(3,1,1)$ with $a(3)=1.4\,\pt$.  By the arguments in the
text, we may assume that the dihedral angles of the exceptional
regions at those vertices are at least $1.32$ (see
\cite[3.11.4]{part4}). Also, the three quasi-regular tetrahedra at
the vertex squander at least $1.5\,\pt$ by a linear programming
bound, if the angle of the quad cluster is at least $1.55$.  Thus,
we assume that the dihedral angles at opposite vertices of the
quad cluster are at most $1.55$.   A linear program also gives
$\tau+0.316\dih>0.3864$ for a quasi-regular tetrahedron.

  If we give bounds of the form
$\tau_x +0.316\dih> b$, for the part of the quad cluster around a
vertex, where $\tau_x$ is the appropriate squander function, then
we obtain
    $$\sum\tau_x > -0.316(2\pi-1.32) + b + 3 (0.3864)$$
for a lower bound on what is squandered.  If the two opposite
vertices give at least  $2(1.4)\,\pt + 0.1317$, then the inclusion
of  two opposite vertices in the separated set of vertices is
justified.  (Recall that $t_4=0.1317$.)  The following
inequalities give the desired result.

$\tau_\mu+0.316\dih > 0.5765$, if $\dih\le1.55$,
$y_4\in[2t_0,2\sqrt2]$. \refno{912536613}

$\tau_0+0.316\dih > 0.5765$ if $\dih\le1.55$, $y_4\ge2\sqrt2$.
\refno{640248153}

$\tau_\nu+0.316\dih_2 >0.2778$, if $y_1\in[2t_0,2\sqrt2]$,
\refno{594902677}

\subsection{  Group 20\dag} 
\oldlabel{A.2.7} \label{app:tri3}
If the circumradius of a quasi-regular tetrahedron is $\ge1.41$,
then by \cite[Section 9.17]{part1}, $\tau>1.8\,\pt$, and many of
the inequalities hold (without further interval arithmetic
calculations).

In Sections \ref{app:tri3} and \ref{app:tri4}, let
$S_1,\ldots,S_5$ be 5 simplices arranged around a common edge
$(0,v)$, with $|v|\in[2,2t_0]$. Let $y_i(S_j)$ be the edges, with
$y_1(S_j)=|v|$ for all $j$, $y_3(S_j)= y_2(S_{j+1})$, and
$y_5(S_j)= y_6(S_{j+1})$. where the subscripts $j$ are extended
modulo 5. In Sections \ref{app:tri3} and \ref{app:tri4},
$\sum\dih(S_j)\le2\pi$. Set $\piF = 2\xiV+\xiG$ if
$\hat\sigma=\vor_0$ in the cases $(y_4\ge2.6, y_1\ge2.2)$ and
$(y_4\ge2.7)$. Set $\piF=0$, otherwise.

$\tau(S_1)+\tau(S_2)+\tau(S_4) > 1.4\,\pt$,
    if $y_4(S_3),y_4(S_5)\ge2\sqrt2$.
    \refno{551665569}

$\tau(S_1)+\tau(S_2)+\tau(S_3) > 1.4\,\pt$,
    if $y_4(S_4),y_4(S_5)\ge2\sqrt2$.
    \refno{824762926}

$\tau(S_1)+\tau(S_2)+(\hat\tau(S_3)-\piF) +\tau(S_4)> 1.4\,\pt
    +D(3,1)$,
    if $y_4(S_3)\in[2t_0,2\sqrt2]$,
    $y_4(S_5)\ge2t_0$, $\dih(S_5)>1.32$,
    \refno{675785884}

$\tau(S_1)+\tau(S_2)+\tau(S_3) +(\hat\tau(S_4)-\piF)> 1.4\,\pt
    +D(3,1)$,
    if $y_4(S_4)\in[2t_0,2\sqrt2]$,
    $y_4(S_5)\ge2t_0$, $\dih(S_5)>1.32$.
    \refno{193592217}

\subsection{  Group 21\dag} 
\label{app:tri4} \oldlabel{A.2.8}


As in Section \ref{app:tri3}, the quasi-regular tetrahedra are
generally compression scored. Define $\piF$ as in Section
\ref{app:tri3}. The constraint $\sum_{(5)}\dih(S_j)=2\pi$ is
assumed.

$\tau(S_1)+\tau(S_2)+\tau(S_3) +\tau(S_4)> 1.5\,\pt$,
    if $y_4(S_5)\ge2\sqrt2$.
    \refno{325738864}

$\tau(S_1)+\tau(S_2)+\tau(S_3) +\tau(S_4)+
    (\hat\tau(S_5)-\piF) > 1.5\,\pt
    +D(3,1)$,
    if $y_4(S_5)\in[2t_0,2\sqrt2]$.
    \refno{314974315}

\subsection{  Group 22} 
\oldlabel{A.4.4.2}

If there are four quasi-regular tetrahedra $\{S_1,\ldots,S_4\}$ at
the central vertex $v$ of the flat quarter $Q$, and if there are
only five standard regions at $v$, then
    $$
    \hat\sigma(Q) +\sum_{(4)} \sigma(S_i) < 0.114.
    $$
    \oldlabel{4.4.2}
\refno{867359387}

\subsection{  Group 23} 
\label{app:455} \oldlabel{A.4.5}
\bigskip

Let $\vor_x=\vor$ for simplices of types $C$ and $C'$.  Let
$\vor_x=\vor$, otherwise.
$$
\vor_x
    < \begin{cases} 0,&y_1\in[2t_0,2\sqrt2]\\
            -0.05,&y_1\in[2t_0,2.696],\\
            -0.119,&y_1\in[2t_0,2.696],\quad\eta_{126}\ge\sqrt2,\\
    \end{cases}
$$
 if $y_4\in[2t_0,2\sqrt2]$.
The bound of $0$ is established in Appendix \ref{app:10}.
    (Even when simplex has type $C$, the
    bound $-0.05$ is based
    on the upper bound $\vor_0$.)
\refno{365179082}

\subsection{  Group 24} 
\oldlabel{A.4.6}

    \begin{eqnarray}
        \sigma_R(D) <
            \begin{cases}
                0,& y_1\in[2t_0,2\sqrt2],\\
                -0.043, & y_1\in[2t_0,2.696],\\
            \end{cases}
    \label{eqn:group24}
    \end{eqnarray}
 for quad regions $R$ constructed from  an anchored
simplex $S$ and adjacent special simplex $S'$. Assume that
$y_4(S)=y_4(S')\in[2\sqrt2,3.2]$, and that the other edges have
lengths in $[2,2t_0]$. The bound $0$ is found in \cite[Lemma
3.13]{formulation}. The bound $-0.043$ is obtained from
deformations, reducing the inequality to the following interval
calculations.

    \bigskip
$\vor_0 < -0.043/2$, if $y_6=2t_0$, $y_1\in[2t_0,2.696]$.
\refno{368244553\dag}

$\vor_0(S)+\vor_0(S(2,y_2,y_3,y_4,2,2))<-0.043$,
    if $y_1\in[2t_0,2.696]$, $y_4\in[2\sqrt2,3.2]$,
        $y_4'\ge2t_0$.
\refno{820900672\dag}

$\vor_0(S)+\vor_0(S(2t_0,y_2,y_3,y_4,2,2))<-0.043$,
    if $y_1\in[2t_0,2.696]$, $y_4\in[2\sqrt2,3.2]$,
        $y_4'\ge2t_0$.
\refno{961078136\dag}

Under certain conditions, Inequality \ref {eqn:group24} can be
improved.
$$
\begin{array}{lll}
     \vor_0(S') &< -0.033,
    \text{ if }\dih(S')\le1.8,\ y_1\in[2,2.12],\ y_4\in[2\sqrt2,3.2],\\
\vor_0(S) &< -0.058,
    \text{ if }\dih(S)\le2.5,\ y_1\in[2t_0,2.696],\ y_4\in[2\sqrt2,3.2],\\
\vor_0(S) &< -0.073,
    \text{ if }\eta_{126}\ge\sqrt2,
        \ y_1\in[2t_0,2.696],\ y_4\in[2\sqrt2,3.2].\dag\\
\end{array}
$$
(The last of these was verified by S. Ferguson.) \refno{424186517}

These combine to give
$$
\vor_0(S)+\vor_0(S') < \begin{cases}  -0.091,&\text{ or }\\
        -0.106,&
        \end{cases}
$$
for the combination of special simplex and anchored simplex under
the stated conditions.

\subsection{  Group 25 (pentagons)} 

There are a few inequalities that arise for pentagonal regions.

\begin{proposition} If the pentagonal region has no flat quarters
and no upright quarters, the subregion $F$ is a pentagon. It
satisfies
    $$
    \begin{array}{lll}
     \vor_0 &< -0.128,\\
    \tau_0 &> 0.36925.
    \end{array}
    $$
\end{proposition}

\begin{proof}  The proof is by deformations and interval calculations. If
a deformation produces a new flat quarter, then the result follows
from \cite[$\A_{13}$]{part4} and Inequality \ref {app:hexquad}. So
we may assume that all diagonals remain at least $2\sqrt2$. If all
diagonals remain at least 3.2, the result follows from the
tcc-bound on the pentagon \cite[Section 5.5]{part4}.  Thus, we
assume that some diagonal is at most $3.2$. We deform the cluster
into the form
    $$(a_1,2,a_2,2,a_3,2,a_4,2,a_5,2),\quad |v_i|=a_i\in\{2,2t_0\}.$$
Assume that $|v_1-v_3|\le3.2$.  If $\max(a_1,a_3)=2t_0$, the
result follows from \cite[$\A_{13}$]{part4} and
Section~\ref{app:hexquad}, Equations \ref{eqn:hexquadsig} and
\ref{eqn:hexquadtau}.

Assume $a_1=a_3=2$. There is a diagonal of the quadrilateral of
length at most $3.23$ because
    $$\Delta(3.23^2,4,4,3.23^2,4,3.2^2)<0.$$
  The result now follows from the following interval arithmetic
  calculations.

(These inequalities are closely related to
\cite[$\A_{21}$]{part4}.)

$\vor_0(S(2,2,2,y_4,2,2))+\vor_0(S(2,2,2,y_4',2,2))
    +\vor_0(S(2,2,2,y_4,y_4',2))
        < -0.128$, if $y_4\in[2\sqrt2,3.2]$,
        $y_4'\in[2\sqrt2,3.23]$.
    \refno{587781327\dag}

$\tau_0(S(2,2,2,y_4,2,2))+\tau_0(S(2,2,2,y_4',2,2))
    +\tau_0(S(2,2,2,y_4,y_4',2))
        > 0.36925$, if $y_4\in[2\sqrt2,3.2]$,
        $y_4'\in[2\sqrt2,3.23]$.
    \refno{807067544\dag}

$\tau_0(2,2,y_3,y_4,2,2) < \tau_0(2t_0,2,y_3,y_4,2,2)$,
    if $y_4 \in[2\sqrt2,3.06]$.
\refno{986970370\dag}

$\vor_0(2,2,y_3,y_4,2,2) > \vor_0(2t_0,2,y_3,y_4,2,2)$,
    if $y_4 \in[2\sqrt2,3.06]$.
\refno{677910379\dag}

$\vor_0< -0.128$, if $y_1=y_2=y_4=2$, $y_5,y_6\in[3.06,3.23]$.
\refno{276168273}

$\tau_0> 0.36925$, if $y_1=y_2=y_3=y_4=2$,
$y_5,y_6\in[3.06,3.23]$. \refno{411203982}

$\tau_0> 0.31$, if $y_1=y_2=y_4=2$, $y_3=2t_0$,
$y_5,y_6\in[3.06,3.23]$. \refno{860823724}

$\vor_0< -0.137-(y_5-2\sqrt2)0.14$,
    if $y_1=y_4=2$, $y_5\in[2\sqrt2,3.23]$, $y_6\in[3.06,3.23]$.
\refno{353116955}

$\tau_0> 0.31+(y_5-2\sqrt2)0.14$,
    if $y_1=y_4=2$, $y_5\in[2\sqrt2,3.23]$, $y_6\in[3.105,3.23]$.
\refno{943315982}

$\tau_0> 0.31+(y_5-2\sqrt2)0.14+(y_6-3.105)0.19$,
    if $y_1=y_4=2$, $y_5\in[2\sqrt2,3.23]$, $y_6\in[3.06,3.105]$.
\refno{941799628}

$\vor_0< 0.009+(y_5-2\sqrt2)0.14$,
    if $y_1=y_4=2$, $y_5\in[2\sqrt2,3.23]$, $y_2=y_6=2$.
\refno{674284283}

$\tau_0> 0.05925+(y_5-2\sqrt2)0.14$,
    if $y_1=y_4=2$, $y_5\in[2\sqrt2,3.23]$, $y_2=y_6=2$.
\refno{775220784}

$\tau_0> 0.05925$,
    if $y_3=2t_0$, $y_5=y_6=2$, $y_4\in[2\sqrt2,3.23]$.
\refno{286076305}

$\tau_0> -(y_4-3.105)0.19$,
    if $y_1=2t_0$, $y_2=y_3=y_5=y_6=2$, $y_4\in[3.06,3.105]$.
\refno{589319960}

\end{proof}

\subsection{  Group 26\dag} 
\oldlabel{A.4.6.3}

Let $Q$ be a quadrilateral region with parameters
    $$(a_1,2t_0,a_2,2,a_3,2,a_4,2t_0),\quad a_i\in\{2,2t_0\}.$$
Assume that $|v_2-v_4|\in[2\sqrt2,3.2]$,
    $|v_1-v_3|\in[3.2,3.46]$. Note that
$$\Delta(4,4,8,2t_0^2,2t_0^2,3.46^2)<0.$$

$\vor_0(Q) < -0.168$. \refno{302085207}

$\tau_0(Q) > 0.352$. \refno{411491283}

\subsection{  Group 27} 
\oldlabel{A.4.7.1} \label{app:471}

Consider a pentagonal region. If the pentagonal region has one
flat quarter and no upright quarters, there is a quadrilateral
region $F$.  It satisfies
    $$
    \begin{array}{lll}
    \vor_0 &< -0.075,\\
    \tau_0 &> 0.176.
    \end{array}
    $$
    \oldlabel{4.6.4}
Break the cluster into two simplices $S=S(y_1,\ldots,y_6)$,
$S'=S(y'_1,y_2,y_3,y_4,y'_5,y'_6)$, by drawing a diagonal of
length $y_4$. Assume that the edge $y'_5\in[2t_0,2\sqrt2]$.  Let
$y_4'$ be the length of the diagonal that crosses $y_4$.
    $$
    \begin{array}{lll}
    \vor_0 &< 2.1327-0.1 y_1 -0.15 y_2 -0.08 y_3 -0.15 y_5\\
            &\qquad -0.15 y_6 - 0.1 y'_1 - 0.17 y'_5 -0.16 y'_6,\\
        &\quad\text{if }\dih(S)<1.9,\ \dih(S')<2.0,\ y_1\in[2,2.2],\
            y_4\ge2\sqrt2,\\
    \vor_0 & < 2.02644 - 0.1 y_1 -0.14 (y_2+y_3)-0.15 (y_5+y_6)
            -0.1 y'_1 - 0.12 (y_5'+y_6'), \\
        &\quad\text{if }y_1\in[2,2.08],\quad y_4\le3.\\
    \vor_0 &+0.419351 \sol < 0.4542 + 0.0238 (y_5+y_6+y_6'),\\
        &\quad\text{if }\ y_4,y_4'\ge2\sqrt2.\\
    \end{array}
    $$
    \oldlabel{A.4.7.1}
The inequalities above are verified in smaller pieces:

$\vor_0 < 1.01 - 0.1 y_1 -0.05 y_2 -0.05 y_3 -0.15 y_5 -0.15 y_6$,
if $\dih\le1.9$, $y_4\ge 2\sqrt2$, and $y_1\le 2.2$.
\refno{131574415}

$\vor_0 < 1.1227 - 0.1 y_1 -0.1 y_2 -0.03 y_3 -0.17 y_5 -0.16
y_6$,
    if $\dih\le2$, $y_4\ge 2\sqrt2$, $y_5\in[2t_0,2\sqrt2]$,
    and $y_2+y_3\le 4.67$.
\refno{929773933}

$\vor_0 < 1.0159 - 0.1 y_1 - 0.08 (y_2+y_3) + 0.04 y_4 -0.15
(y_5+y_6)$,
    if $y_4\in[2\sqrt2,3]$, $y_1\in[2,2.08]$.
\refno{223261160}

$\vor_0 < 1.01054 - 0.1 y_1 -0.06 (y_2+y_3) -0.04 y_4 -0.12
(y_5+y_6)$,
    if $y_4\in[2\sqrt2,3]$, $y_5\in[2t_0,2\sqrt2]$.
\refno{135018647}

Let $Q$ be the quadrilateral subcluster expressed as a union of
two simplices $S$ and $S'$, as above
    $$\vor_0(Q)+0.419351 \sol(Q) < 0.4542 + 0.0238 (y_5+y_6+y_6'),$$
if $y_5'\in[2t_0,2\sqrt2]$ and both diagonals
have length $\ge2\sqrt2$. \refno{559676877}

\subsection{  Group 28\dag} 
\oldlabel{A.4.10}


If $v$ is a vertex of an exceptional cluster and there are exactly
$4$ quasi-regular tetrahedra along $(0,v)$, then there are 5
vertices $v_1,\ldots,v_5$ adjacent to $v$. We have
    $$\sum_{(5)}(|v-v_i|+|v_i|) > 20.42.
    $$\refno{615073260}
    \oldlabel{4.10.2}

If $v_1$ and $v_5$ are the two vertices on an exceptional cluster,
and if the additional hypothesis $|v_1-v_5|\ge2\sqrt2$ holds, then
    $$\sum_{(5)}(|v-v_i|+|v_i|) > 20.76.
    $$
    \refno{844430737}
    \oldlabel{4.10.3}

\subsection{  Group 29} 
\oldlabel{A.4.12} \label{app:pentB}

    $$
    \vor_0 < -0.136\quad\text{and }\tau_0 > 0.224,
    $$
    \oldlabel{A.4.12.4}
for a combination of anchored simplex $S$ and special simplex
$S'$, with $y_1(S)\in[2.696,2\sqrt2]$,
$y_2(S),y_6(S)\in[2.45,2t_0]$, $y_4(S)\in[2\sqrt2,3.2]$, and with
cross-diagonal at least $2t_0$. This inequality can be verified by
proving the following inequalities in lower dimension. In the
first four $y_1\in[2.696,2\sqrt2]$, $y_2,y_6\in[2.45,2t_0]$,
$y_4\in[2\sqrt2,3.2]$, and $y_4'\ge2t_0$ (the cross-diagonal).
\bigskip

$\vor_0(S(y_1,\ldots,y_6))+\vor_0(S(2,y_2,y_3,y_4,2,2))<-0.136$
\refno{967376139\dag}

$\vor_0(S(y_1,\ldots,y_6))+\vor_0(S(2t_0,y_2,y_3,y_4,2,2))<-0.136$
\refno{666869244\dag}

$\tau_0(S(y_1,\ldots,y_6))+\tau_0(S(2,y_2,y_3,y_4,2,2))>0.224$
\refno{268066802\dag}

$\tau_0(S(y_1,\ldots,y_6))+\tau_0(S(2t_0,y_2,y_3,y_4,2,2))>0.224$
\refno{508108214\dag}

$\vor_0(S(y_1,\ldots,y_6)) < -0.125$, if $y_1\in[2.696,2\sqrt2]$,
    $y_2,y_6\in[2.45,2t_0]$, $y_5=2t_0$.
\refno{322505397\dag}

$\vor_0(S(y_1,\ldots,y_6)) < 0.011$, if $y_1\in[2.696,2\sqrt2]$,
       $y_5=2t_0$.
\refno{736616321\dag}

$\tau_0(S(y_1,\ldots,y_6)) > 0.17$, if $y_1\in[2.696,2\sqrt2]$,
    $y_2,y_6\in[2.45,2t_0]$, $y_5=2t_0$.
\refno{689417023\dag}

$\tau_0(S(y_1,\ldots,y_6)) > 0.054$, if $y_1\in[2.696,2\sqrt2]$,
       $y_5=2t_0$.
\refno{748466752\dag}

\subsection{  Group 30}
\label{app:hexA}
$$\vor_0 < -0.24\text{ and }\tau_0 > 0.346,
    $$
    \oldlabel{A.4.12.5}
for an anchored simplex $S$ and simplex $S'$ with edge parameters
$(3,2)$ in a hexagonal cluster, with $y_2(S)=y_2(S')$,
$y_3(S)=y_3(S')$, $y_4(S)=y_4(S')$, $y_1(S)\in[2.696,2\sqrt2]$,
$y_4(S)\in[2\sqrt2,3.2]$, $y_2(S),y_6(S)\in[2.45,2t_0]$, and
$$\max(y_5(S'),y_6(S'))\in[2t_0,2\sqrt2],\quad
\min(y_5(S'),y_6(S'))\in[2,2t_0].$$ This breaks into separate
interval calculations for $S$ and $S'$.

This inequality  results from the following four inequalities:

$\vor_0(S) < -0.126$ and $\tau_0(S) > 0.16$ \refno{369386367\dag}

$\vor_0(S') < -0.114$ and $\tau_0(S') >0.186$ (There are two cases
for each, depending on which of $y_5,y_6$ is longer.)
\refno{724943459\dag}

\subsection{  Group 31}
\label{app:hexB}
$$\vor_0 < -0.149\quad\text{and }\tau_0 > 0.281,
    $$
for a quadrilateral subcluster with both diagonals $\ge2\sqrt2$,
and parameters
$$(a_1,b_1,a_2,b_2,a_3,b_3,a_4,b_4),$$
with $a_4\in[2.45,2t_0]$, and $b_4\in[2t_0,2\sqrt2]$.
\refno{836331201\dag}

\subsection{  Group 32}
\label{app:hexC}
    $$\vor_0 < -0.254\quad\text{and }\tau_0 > 0.42625,
    $$
for a combination of anchored simplex $S$ and quadrilateral
cluster $Q$.  It is assumed that $y_1(S)\in[2.696,2\sqrt2]$,
$y_2(S),y_6(S)\in[2.45,2t_0]$. The adjacent quadrilateral
subcluster is assumed to have both diagonals $\ge2\sqrt2$, and
parameters
$$(a_1,b_1,a_2,b_2,a_3,b_3,a_4,b_4),$$
with $b_4\in[2\sqrt2,3.2]$. The verification of this inequality
reduces to separate inequalities for the anchored simplex and
quadrilateral subcluster. For the anchored simplex we use the
bounds $\vor_0(S')<-0.126$, $\tau_0(S')>0.16$ that have already
been established above.  We then show that the quad cluster
satisfies

$\vor_0 < -0.128$ and $\tau_0 > 0.26625$. \refno{327474205\dag}

For this, use deformations to reduce either to the case where the
diagonal is $2\sqrt2$, or to the case where $b_1=b_2=b_3=2$,
$a_2,a_3\in\{2,2t_0\}$.  When the diagonal is $2\sqrt2$, the flat
quarter can be scored by \cite[$\A_{13}$]{part4}:
    $(\vor_0<0.009,\tau_0>0.05925)$.
(There are two cases depending on which direction the diagonal of
length $\sqrt2$ runs.)

\twocolumn
\section{Appendix Hexagonal Inequalities}

There are a number of inequalities that have been particularly
designed for standard regions that are hexagons.  This appendix
describes those inequalities.  They are generally inequalities
involving more than six variables, and because of current
technological limitations on interval arithmetic, we were not able
to prove these inequalities directly with interval arithmetic.

Instead we give various lemmas that deduce the inequalities from
inequalities in a smaller number of variables (small enough to
prove by interval arithmetic.)

\subsection{Statement of results} 

There are a number of inequalities that hold in special situations
when there is a hexagonal region.  After stating all of them, we
will turn to the proofs.

\begin{enumerate}
\item\label{app:hex1}
 If there are no flat quarters and no upright quarters (so that
there is a single subregion $F$), then
    \begin{eqnarray}
    \vor_0 &< -0.212\\
    \tau_0 &> 0.54525.
    \end{eqnarray}
    \oldlabel{eqn:4.6.1}

\item \label{app:hex2}
If there is one flat quarter and no upright quarters, there is
a pentagonal subregion $F$.  It satisfies
    $$
    \begin{array}{lll}
    \vor_0 &< -0.221\\
    \tau_0 &> 0.486.
    \end{array}
    $$

\item\label{app:hex3}
If there are two flat quarters and no upright quarters, there is a
quadrilateral subregion $F$.  It satisfies
    $$
    \begin{array}{lll}
    \vor_0 &< -0.168,\\
    \tau_0 &> 0.352.
    \end{array}
    $$
These are twice the constants appearing in \ref{app:hex11};

\item\label{app:hex4}
If there is an edge of length between $2t_0$ and $2\sqrt2$
running between two opposite corners of the hexagonal cluster, and
if there are no flat or upright quarters on one side, leaving a
quadrilateral region $F$, then $F$ satisfies
    $$
    \begin{array}{lll}
    \vor_0 &< -0.075,\\
    \tau_0 &> 0.176.
    \end{array}
    $$

\item\label{app:hex5}
If the hexagonal cluster has an upright diagonal with context
$(4,2)$, and if there are no flat quarters
(Figure~\ref{fig:hex42}, then the hexagonal cluster $R$ satisfies
    $$
    \begin{array}{lll}
    \sigma_R &< -0.297,\\
    \tau_R &> 0.504.
    \end{array}
    $$
\begin{figure}[htb]
  \centering
  \includegraphics{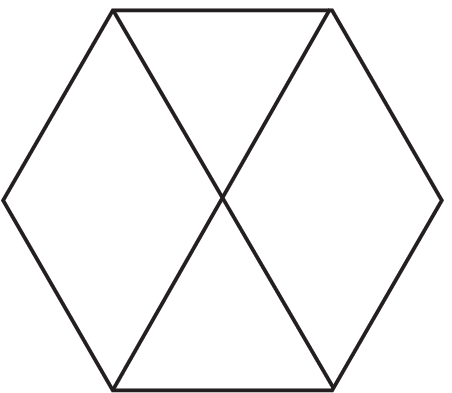}
  \caption{A hexagonal cluster with context $(4,2)$.}
  \label{fig:hex42}
\end{figure}

\item\label{app:hex6}
If the hexagonal cluster has an upright diagonal with context
$(4,2)$, and if there is one unmasked flat quarter
(Figure~\ref{fig:hex42a}, let $\{F\}$ be the set of four
subregions around the upright diagonal. (That is, take all
subregions except for the flat quarter.) In the following
inequality and Inequality \ref{app:hex7}, let $\sigma_R^+$ be
defined as $\sigma_R$ on quarters, and $\vor_x$ on other anchored
simplices.  $\tau_R^+$ is the adapted squander function.
    $$
    \begin{array}{lll}
    \sum_{(4)}\sigma_R^+ &< -0.253,\\
    \sum_{(4)}\tau_R^+ &> 0.4686.\\
    \end{array}
    $$
\begin{figure}[htb]
  \centering
  \includegraphics{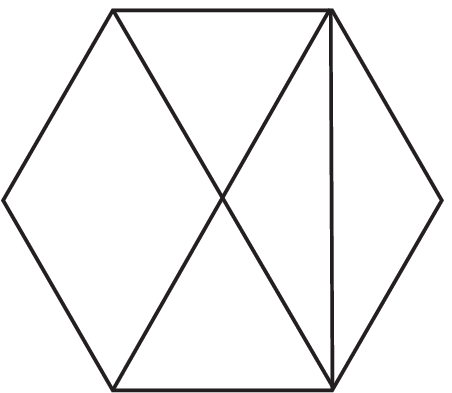}
  \caption{A hexagonal cluster with context $(4,2)$.}
  \label{fig:hex42a}
\end{figure}

\item\label{app:hex7}
If the hexagonal cluster has an upright diagonal with context
$(4,2)$, and if there are two unmasked flat quarters
(Figure~\ref{fig:hex42b}, let $\{F\}$ be the set of four
subregions around the upright diagonal. (That is, take all
subregions except for the flat quarters.)
    $$
    \begin{array}{lll}
    \sum_{(4)}\sigma_R^+ &< -0.2,\\
    \sum_{(4)}\tau_R^+ &> 0.3992.\\
    \end{array}
    $$
\begin{figure}[htb]
  \centering
  \includegraphics{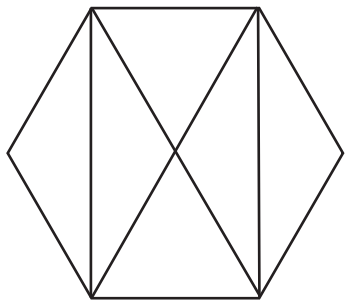}
  \caption{A hexagonal cluster with context $(4,2)$.}
  \label{fig:hex42b}
\end{figure}

\item\label{app:hex8}
If the hexagonal cluster has an upright diagonal in context
$(4,1)$, and if there are no flat quarters, let $\{F\}$ be the set
of four subregions around the upright diagonal. Assume that the
edge opposite the upright diagonal on the anchored simplex has
length at least $2\sqrt2$.
 (See Figure~\ref{fig:hex41}.)
    $$
    \begin{array}{lll}
    \vor_{0,R}(D)+\sum_{(3)}\sigma(Q) &< -0.2187\\
    \tau_{0,R}(D)+\sum_{(3)}\tau(Q) &> 0.518.
    \end{array}
    $$
\begin{figure}[htb]
  \centering
  \includegraphics{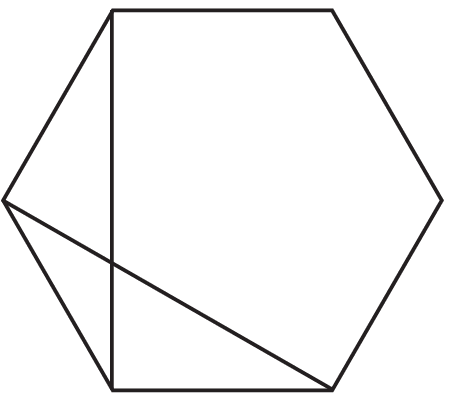}
  \caption{A hexagonal cluster with context $(4,1)$.}
  \label{fig:hex41}
\end{figure}

\item\label{app:hex9}
In this same context, let $F$ be the pentagonal subregion along
the upright diagonal. It satisfies
    \begin{eqnarray}
    \vor_0 &< -0.137,\\
    \tau_0 &> 0.31.
    \end{eqnarray}

\item\label{app:hex10}
If the hexagonal cluster has an upright diagonal in context
$(4,1)$, and if there is one unmasked flat quarter, let $\{F\}$ be
the set of four subregions around the upright diagonal. Assume
that the edge opposite the upright diagonal on the anchored
simplex has length at least $2\sqrt2$.
 (There are five subregions, shown in Figure~\ref{fig:hex41a}.)
    $$
    \begin{array}{lll}
    \vor_{0,R}(D)+\sum_{(3)}\sigma(Q) &< -0.1657,\\
    \tau_{0,R}(D)+\sum_{(3)}\tau(Q) &> 0.384.
    \end{array}
    $$
\begin{figure}[htb]
  \centering
  \includegraphics{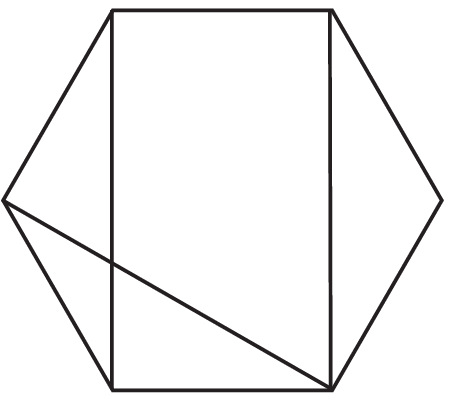}
  \caption{A hexagonal cluster with context $(4,1)$.}
  \label{fig:hex41a}
\end{figure}

\item\label{app:hex11}
In this same context, let $F$ be the quadrilateral subregion in
Figure~\ref{fig:hex41a}. It satisfies
    $$
    \begin{array}{lll}
    \vor_0 &< -0.084,\\
    \tau_0 &> 0.176.
    \end{array}
    $$

\end{enumerate}

\subsection{Proof of inequalities}
\label{app:hexquad}

\begin{proposition}  Inequalities \ref{app:hex1} --\ref{app:hex11}  are valid.
\end{proposition}

 We prove the inequalities in reverse order \ref{app:hex11}--
\ref{app:hex1}. The bounds $\vor_0<0.009$ and $\tau_0>0.05925$
from \cite[$\A_{13}$]{part4} for what a flat quarters with
diagonal $\sqrt8$  will be used repeatedly.  Some of the proofs
will make use of tcc-bounds, which are described in
\cite[Section~5.2]{part4}.

\begin{proof}
({\it Inequality \ref{app:hex10} and Inequality \ref{app:hex11}.})
Break the quadrilateral cluster into two simplices $S$ and $S'$
along the long edge of the anchored simplex $S$.  The anchored
simplex $S$ satisfies $\tau(S)\ge0$, $\sigma(S)\le0$.  The other
simplex satisfies $\tau_0(S')>0.176$ and $\vor_0(S')<-0.084$ by an
interval calculation \ref{app:p11}. This gives Inequality
\ref{app:hex11}.  For Inequality \ref{app:hex10}, we combine these
bounds with the linear programming bound on the four anchored
simplices around the upright diagonal.  From the inequalities
    \cite[$\A_2$]{part4} --
    \cite[$\A_{7}$]{part4},
    \cite[$\A_{22}$]{part4},
    \cite[$\A_{24}$]{part4},
we find that they score $<-0.0817$ and squander $>0.208$.  Adding
these to the bounds from Inequality \ref{app:hex11}, we obtain
Inequality \ref{app:hex10}.
\end{proof}

\begin{proof}
{\it (Inequality (\ref{app:hex8}) and (\ref{app:hex9}).)} The
pentagon is a union of an anchored simplex and a quadrilateral
region. LP-bounds similar to those in the previous paragraph and
based on the inequalities of \cite{part4} show that the loop
scores at most $-0.0817$ and squanders at least $0.208$.  If we
show that the quadrilateral satisfies
    \begin{eqnarray}
    \vor_0&<-0.137,\label{eqn:hexquadsig}\\
    \tau_0&>0.31,\label{eqn:hexquadtau}
    \end{eqnarray}
then Inequalities (\ref{app:hex8}) and (\ref{app:hex9}) follow. If
by deformations a diagonal of the quadrilateral drops to
$2\sqrt2$, then the result follows from Inequalities \ref{app:p8}
and \cite[$\A_{13}$]{part4}. By this we may now assume that the
quadrilateral has the form
    $$(a_1,2,a_2,2,a_3,2,a_4,b_4),\quad a_2,a_3\in\{2,2t_0\}.$$
If the diagonals drop under $3.2$ and $\max(a_2,a_3)=2t_0$, again
the result follows from Inequalities \ref{app:p8} and
\cite[$\A_{13}$]{part4}. If the diagonals drop under $3.2$ and
$a_2=a_3=2$, then the result follows from \cite[$\A_{19}$]{part4}.
So finally we attain by deformations $b_4=2\sqrt2$ with both
diagonals greater than $3.2$. But this does not exist, because
    $$\Delta(4,4,4,3.2^2,4,8,3.2^2)<0.$$
\end{proof}

\begin{proof}
{\it (Inequality \ref{app:hex5}, Inequality \ref{app:hex6}, and
Inequality \ref{app:hex7}.)} Inequalities \ref{app:hex7} are
derived in \cite[Section~5.11]{part4}. Inequalities
\ref{app:hex5}, \ref{app:hex6} are LP-bounds based on
    \cite[$\A_2$]{part4}, \cite[$\A_7$]{part4}, \cite[$\A_{22}$]{part4}.

{\it Proof of \ref{app:hex4}.} Deform as in \cite{part4}.  If at
any point a diagonal of the quadrilateral drops to $2\sqrt2$, then
the result follows from \cite[$\A_{13}$]{part4} and Inequality
\ref{app:hex11}:
    $$
    \begin{array}{lll}
    \vor_0 &< 0.009 -0.084 = -0.075,\\
    \tau_0 &> 0 + 0.176 = 0.176,\\
    \end{array}
    $$
Continue deformations until the quadrilateral has the form
    $$(a_1,2,a_2,2,a_3,2,a_4,b_4),\quad a_2,a_3\in\{2,2t_0\}.$$
There is necessarily a diagonal of length $\le3.2$, because
    $$\Delta(4,4,3.2^2,8,4,3.2^2)<0.$$
Suppose the diagonal between vertices $v_2$ and $v_4$ has length
at most $3.2$.  If $a_2=2t_0$ or $a_3=2t_0$, the result follows
from \cite[$\A_{13}$]{part4} and Inequality \ref{app:hex11}. Take
$a_2=a_3=2$. Inequality \ref{app:hex4} now follows from
\cite[$\A_{19}$]{part4}.
\end{proof}

\begin{proof}
{\it (Inequality \ref{app:hex3}).}
 We prove that the quadrilateral
satisfies
    $$
    \begin{array}{lll}
    \vor_0 &< -0.168\\
    \tau_0 &> 0.352.\\
    \end{array}
    $$
There are two types of quadrilaterals.  In (a), there are two flat
quarters whose central vertices are opposite corners of the
hexagon.  In (b), the flat quarters share a vertex.
  We consider case (a) first.

Case (a). We deform the quadrilateral as in \cite{part4}.  If at
any point there is a diagonal of length at most 3.2, the result
follows from Inequality \ref{app:hex10} and Inequality
\ref{app:hex11}. Otherwise, the deformations give us a
quadrilateral
    $$(a_1,2,a_2,2t_0,a_3,2,a_4,2), \quad a_i \in\{2,2t_0\}.$$
The tcc approximation now gives the result (see
\cite[5.3)]{part4}.

Case (b). Label the vertices of the quadrilateral
 $v_1,\ldots,v_4$, where $(v_1,v_2)$ and $(v_1,v_4)$
are the diagonals of the flat quarter.  Again, we deform the
quadrilateral.  If at any point of the deformation, we find that
$|v_1-v_3|\le 3.2$, the result follows from Inequalities
\ref{app:hex10}, \ref{app:hex11}. If during the deformation
$|v_2-v_4|\le2\sqrt2$, the result follows from
\cite[$\A_{13}$]{part4} and the interval calculations
\ref{app:p3}. If the diagonal $(v_2,v_4)$ has length at least
$3.2$ throughout the deformation, we eventually obtain a
quadrilateral of the form
    $$(a_1,2t_0,a_2,2,a_3,2,a_4,2t_0),\quad a_i\in\{2,2t_0\}.$$
But this does not exist:
$$\Delta(4,4,3.2^2,(2t_0)^2,(2t_0)^2,3.2^2)<0.$$

We may assume that $|v_2-v_4|\in[2\sqrt2,3.2]$.  The result now
follows from interval calculations \ref{app:p3}.
\end{proof}

\begin{proof}
{\it (Inequality \ref{app:hex2}).} This case requires more effort.
We show that
    $$
    \begin{array}{lll}
    \vor_0 &< -0.221\\
    \tau_0 &> 0.486\\
    \end{array}
    $$
Label the corners $(v_1,\ldots,v_5)$ cyclically with $(v_1,v_5)$
the diagonal of the flat quarter in the hexagonal cluster. We use
the deformation theory of \cite{part4}.  The proof appears in
steps $(1),\ldots,(6)$.

(1) If during the deformations, $|v_1-v_4|\le3.2$ or
$|v_2-v_5|\le3.2$, the result follows from Inequalities
\ref{app:hexquad} and \ref{app:hex11}.  We may assume this does
not occur.

(2) If an edge $(v_1,v_3)$, $(v_2,v_4)$, or $(v_3,v_5)$ drops to
$2\sqrt2$, continue with deformations that do not further decrease
this diagonal. If $|v_1-v_3|=|v_3-v_5|=2\sqrt2$, then the result
follows from
    \cite[$\A_{13}$]{part4} and interval calculations \ref{app:p2a}.

If we have $|v_1-v_3|=2\sqrt2$, deform the figure to the form
    $$(a_1,2,a_2,2,a_3,2,a_4,2,a_5,2t_0),\quad a_2,a_4,a_5\in\{2,2t_0\}.$$
Once it is in this form, break the flat quarter $(0,v_1,v_2,v_3)$
from the cluster and deform $v_3$ until $a_3\in\{2,2t_0\}$.  The
result follows from an interval calculation \ref{app:p2b}.

We handle a boundary case of the preceding calculation separately.
After breaking the flat quarter off, we have the cluster
    $$(a_1,2\sqrt2,a_3,2,a_4,2,a_5,2t_0),\quad a_3,a_4,a_5\in\{2,2t_0\}.$$
If $|v_1-v_4|=3.2$, we break the quadrilateral cluster into two
pieces along this diagonal and use interval calculations
\ref{app:p2c} to conclude the result. This completes the analysis
of the case $|v_1-v_3|=2\sqrt2$.

(3) If $|v_2-v_4|\le3.2$, then deform until the cluster has the
form
    $$(a_1,2,a_2,2,a_3,2,a_4,2,a_5,2t_0),\quad a_1,a_3,a_5\in\{2,2t_0\}.$$
Then cut along the special simplex to produce a quadrilateral.
Disregarding cases already treated by the interval calculations
\ref{app:p2c}, we can deform it to
    $$(a_1,2,a_2,2\sqrt2,a_4,2,a_5,2t_0),\quad a_i\in\{2,2t_0\},$$
with diagonals at least $3.2$. The result now follows from the
interval calculations \ref{app:p2d}.

In summary of (1), (2), (3), we find that by disregarding cases
already considered, we may deform the cluster into the form
    $$(a_1,2,a_2,2,a_3,2,a_4,2,a_5,2t_0),\quad a_i\in\{2,2t_0\},$$
$|v_1-v_3|>2\sqrt2$, $|v_3-v_5|>2\sqrt2$, $|v_2-v_4|>3.2.$

(4) Assume $|v_1-v_3|,|v_3-v_5|\le3.2$. If
$\max(a_1,a_3,a_5)=2t_0$, we invoke interval calculations
\ref{app:p2b}  and \cite[$\A_{13}$]{part4} to prove the
inequalities.
 So we may assume $a_1=a_3=a_5=2$.  The result now follows from
interval calculations \ref{app:p2e}. This completes the case
$|v_1-v_3|,|v_3-v_5|\le3.2$.

(5) Assume $|v_1-v_3|,|v_3-v_5|\ge3.2$. We deform to
    $$(a_1,2,a_2,2,a_3,2,a_4,2,a_5,2t_0),\quad a_i\in\{2,2t_0\}.$$
If $a_2=2t_0$ and  $a_1=a_3=2$, then the simplex does not exist by
Section  \cite[5.6]{part4}. Similarly, $a_4=2t_0$, $a_5=a_3=2$
does not exist. The tcc bound gives the result except when
$a_2=a_4=2$. The condition $|v_2-v_4|\ge3.2$ forces $a_3=2$. These
remaining cases are treated with the interval calculations
\ref{app:p2f}.

(6) Assume $|v_1-v_3|\le3.2$ and $|v_3-v_5|\ge3.2$. This case
follows from deformations, interval calculations \ref{app:p2b},
and \ref{app:p2g}. This completes the proof of Inequalities
\ref{app:hex2}.
\end{proof}

\smallskip
\begin{proof}
{\it (Inequality \ref{app:hex1}).}
 Label the corners of the hexagon
$v_1,\ldots,v_6$.  The proof to this inequality is similar to the
other cases.  We deform the cluster by the method of IV until it
breaks into pieces that are small enough to be estimated by
interval calculations. If a diagonal between opposite corners has
length at most $3.2$, then the hexagon breaks into two
quadrilaterals and the result follows from Inequality
\ref{app:hexquad}.

If a flat quarter is formed during the course of deformation, then
the result follows from Inequality \ref{app:hex2} and
\cite[$\A_{13}$]{part4}. Deform until the hexagon has the form
    $$(a_1,2,a_2,2,\ldots,a_6,2),\quad a_i\in\{2,2t_0\}.$$
We may also assume that the hexagon is convex (see
\cite[Section~4.11]{part4}).

If there are no special simplices, we consider the tcc-bound. The
tcc-bound implies Inequality \ref {app:hex1}, except when $a_i=2$,
for all $i$. But if this occurs, the perimeter of the convex
spherical polygon is $6\arc(2,2,2)=2\pi$.  Thus, there is a pair
of antipodal points on the hexagon.  The hexagon degenerates to a
lune with vertices at the antipodal points.  This means that some
of the angles of the hexagon are $\pi$.   One of the tccs has the
form $C(2,1.6,\pi)$, in the notation of \cite[Lemma~4.11]{part4}.
With this extra bit of information, the tcc bound implies
Inequality \ref {app:hex1}.

If there is one special simplex, say $|v_5-v_1|\in[2\sqrt2,3.2]$,
we remove it. The score of the special simplex is
\cite[Inequalities $\A_{13}$]{part4}
    $$
    \begin{array}{lll}
    \vor_0 &< 0,\quad \tau_0>0.05925,\quad\text{if }\max(|v_1|,|v_5|)=2t_0,\\
    \vor_0 &< 0.0461,\quad \tau_0>0,\quad\text{if }|v_1|=|v_5|=2,\\
    \end{array}
    $$
The resulting pentagon can be deformed. If by deformations, we
obtain $|v_2-v_5|=3.2$ or $|v_1-v_4|=3.2$, the result follows from
Inequalities \ref{app:hexquad} and the following two interval
calculations

\oldlabel{A.4.6.1.a}

$\vor_0 < -0.212 -0.0461 +0.137$, if $y_4=3.2$,
$y_5\in[2\sqrt2,3.2]$. \refno{725257062\dag}

$\tau_0 > 0.54525 - 0 - 0.31$. if $y_4=3.2$,
$y_5\in[2\sqrt2,3.2]$. \refno{977272202\dag}

If $|v_5-v_1|=2\sqrt2$, we use Inequality \ref{app:hex2} and
\cite[$\A_{13}$]{part4} unless $|v_1|=|v_5|=2$.  If
$|v_1|=|v_5|=2$, we use the interval calculations \ref{app:p1b}.
If a second special simplex forms during the deformations, the
result follows from the interval calculations \ref{app:p1c}.

The final case of Inequality \ref{app:hex1}  to consider is that
of two special simplices. We divide this into two cases. (a) The
central vertices of the specials are $v_2$ and $v_6$.  (b) The
central vertices are opposite $v_1$ and $v_4$. In case (a), the
result follows by deformations and interval calculations
\ref{app:p1d}. In case (b), the result follows by deformations and
interval calculations \ref{app:p1e}. This completes the proof of
Inequalities \ref{app:hex1} and the proof of the Proposition.
\end{proof}

\subsection{  Group\dag} 
\label{app:p1b} \oldlabel{A.4.6.1.b}

Let $R$ be a pentagonal region with parameters
    $$(2,2,a_2,2,a_3,2,a_4,2,2,2\sqrt2),\quad a_2,a_3,a_4\in\{2,2t_0\}.$$
Since this is a pentagonal region, we may discard any edge-length
combinations that produce $\Delta<0$. Assume $|v_1-v_3|\ge3.2$,
$|v_3-v_5|\ge3.2$. Under these conditions the following
inequalities hold.

$\vor_0+0.0461 < -0.212$ \refno{583626763}

$\tau_0> 0.54525$ \refno{390951718}

\subsection{Group\dag} 
\label{app:p1c} \oldlabel{A.4.6.1.c}

Let $R$ be a pentagonal region with parameters
    $$(a_1,2,a_2,2,a_3,2,a_4,2,a_5,b_5),\quad a_i\in\{2,2t_0\}.$$
Since this is a pentagonal region, we may discard any edge-length
combinations that produce $\Delta<0$. Assume $|v_1-v_3|$,
$|v_3-v_5|$, $b_5=|v_5-v_1|\in[2\sqrt2,3.2]$. Under these
conditions the following inequalities hold.

$\vor_0 + 0.0461 < -0.212$ \refno{621852152}

$\tau_0 > 0.54525$ \refno{207203174}

\subsection{  Group\dag} 
\label{app:p1d} \oldlabel{A.4.6.1.d}

Let $R$ be a hexagonal region with parameters
    $$(a_1,2,a_2,2,a_3,2,a_4,2,a_5,2,a_6,2),\quad a_i\in\{2,2t_0\}.$$
Since this is a hexagonal region, we may discard any edge-length
combinations that produce $\Delta<0$. Assume
$|v_1-v_5|,|v_1-v_3|\in[2\sqrt2,3.2]$, $|v_1-v_4|\ge3.2$,
$|v_3-v_5|\ge2\sqrt2$. Under these conditions the following
inequalities hold.

$\vor_0 < -0.212$. \refno{368258024}

$\tau_0 > 0.54525$. \refno{564618342}

\subsection{  Group\dag} 
\label{app:p1e} \oldlabel{A.4.6.1.e}

Let $R$ be a hexagonal region with parameters
    $$(a_1,2,a_2,2,a_3,2,a_4,2,a_5,2,a_6,2),\quad a_i\in\{2,2t_0\}.$$
Since this is a hexagonal region, we may discard any edge-length
combinations that produce $\Delta<0$. Assume
$|v_2-v_6|,|v_3-v_5|\in[2\sqrt2,3.2]$, $|v_2-v_5|\in[3.2,3.78]$,
$|v_3-v_6|\ge3.2$. Under these conditions the following
inequalities hold.

$\vor_0 < -0.212$. \refno{498774382}

$\tau_0 > 0.54525$. \refno{544865225}

\subsection{Group\dag}
\label{app:p2a} 
$\vor_0 < -0.221 - 2(0.009)$, if $y_4=2\sqrt2$,
    $y_5\in[2\sqrt2,3.2]$, and $y_6\in[2t_0,2\sqrt2]$.
\refno{234734606\dag}

$\tau_0 > 0.486 - 2 (0.05925)$, if $y_4=2\sqrt2$,
    $y_5\in[2\sqrt2,3.2]$, $y_6\in[2t_0,2\sqrt2]$.
\refno{791682321\dag}

\subsection{  Group\dag} 
\label{app:p2b} \oldlabel{A.4.6.2.b}

Let $Q$ be a quadrilateral region with parameters
    $$(a_1,b_1,a_2,2,a_3,2,a_4,2t_0),\quad a_2,a_3,a_4\in\{2,2t_0\}.$$
Since this is a quadrilateral region, we may discard any
edge-length combinations that produce $\Delta<0$. Assume
$|v_1-v_3|\ge3.2$, $|v_2-v_4|\ge2\sqrt2$.
$b_1=|v_1-v_2|\in[2\sqrt2,3.2]$.  Under these conditions the
following inequalities hold.

$\vor_0(Q) < -0.221 -0.009$. \refno{995351614}

$\tau_0(Q) > 0.486 - 0.05925$. \refno{321843503}

\subsection{Group\dag}
\label{app:p2c}
$\vor_0< -0.19-(y_5-\sqrt2)0.14$, if $y_4=2$,
$y_5\in[2\sqrt2,3.2]$, and  $y_6\in[3.2,3.47]$. \refno{354217730}

$\tau_0 > 0.281$, if $y_4=2$, $y_5\in[2\sqrt2,3.2]$,
$y_6\in[3.2,3.23]$. \refno{595674181}

$\vor_0< -0.11$, if $y_4=2$, $y_5=2t_0$, $y_6=3.2$.
\refno{547486831\dag}

$\tau_0> 0.205$, if $y_4=2$, $y_5=2t_0$, $y_6=3.2$.
\refno{683897354\dag}

$\vor_0 < 0.009 + (y_5-2\sqrt2)0.14$, if
    $y_5\in[2\sqrt2,3.2]$, $y_4=y_6=2$.
\refno{938003786}

\subsection{  Group\dag} 
\label{app:p2d} \oldlabel{A.4.6.2.d}

Let $Q$ be a quadrilateral region with parameters
    $$(a_1,2,a_2,2\sqrt2,a_3,2,a_4,2t_0),\quad a_i\in\{2,2t_0\}.$$
Since this is a quadrilateral region, we may discard any
edge-length combinations that produce $\Delta<0$. Assume
$|v_1-v_3|\ge3.2$, $|v_2-v_4|\ge3.2$. Under these conditions the
following inequalities hold.

$\vor_0(Q) < -0.221 - 0.0461$. \refno{109046923}

$\tau_0(Q) > 0.486$. \refno{642590101}

\subsection{  Group\dag} 
\label{app:p2e} \oldlabel{A.4.6.2.e}

Let $R$ be a pentagonal region with parameters
    $$(2,2,a_2,2,2,2,a_4,2,2,2t_0),\quad a_2,a_4\in\{2,2t_0\}.$$
Since this is a pentagonal region, we may discard any edge-length
combinations that produce $\Delta<0$. Assume
$|v_1-v_3|\in[2\sqrt2,3.2]$, $|v_3-v_5|\in[2\sqrt2,3.2]$. Under
these conditions the following inequalities hold.

$\vor_0(R) < -0.221$. \refno{160800042}

$\tau_0(R) > 0.486$. \refno{690272881}

\subsection{  Group\dag} 
\label{app:p2f} \oldlabel{A.4.6.2.f}

Let $R$ be a pentagonal region with parameters
    $$(a_1,2,2,2,2,2,2,2,a_5,2t_0),\quad a_1,a_5\in\{2,2t_0\}.$$
Since this is a pentagonal region, we may discard any edge-length
combinations that produce $\Delta<0$. Assume $|v_1-v_3|\ge3.2$,
$|v_3-v_5|\ge3.2$. Assume
    $$
    \begin{array}{lll}
    \dih(0,v_3,v_4,v_5)&+\dih(0,v_3,v_5,v_1)+
        \dih(0,v_3,v_1,v_2)\\
        &\ge\dih(S(2,2,2,3.2,2,2))\\
        &=\arccos(-53/75).
    \end{array}
    $$
Under these conditions the following inequalities hold.

$\vor_0(R) < -0.221$. \refno{713930036}

$\tau_0(R) > 0.486$. \refno{724922588}

\subsection{  Group\dag} 
\oldlabel{A.4.6.2.g} \label{app:p2g}

Let $R$ be a pentagonal region with parameters
    $$(2,2,a_2,2,2,2,a_4,2,a_5,2t_0),\quad a_2,a_4,a_5\in\{2,2t_0\}.$$
Since this is a pentagonal region, we may discard any edge-length
combinations that produce $\Delta<0$. Assume
$|v_1-v_3|\in[2\sqrt2,3.2]$, $|v_1-v_4|\ge3.2$, $|v_3-v_5|\ge3.2$.
Under these conditions the following inequalities hold.

$\vor_0(R) < -0.221$. \refno{821730621}

$\tau_0(R) > 0.486$. \refno{890642961}

\subsection{Group\dag}
\label{app:p3}

$\vor_0 < -0.168 - 0.009$, $y_4=2\sqrt2$,
$y_5,y_6\in[2t_0,2\sqrt2]$. \refno{341667126\dag}

$\tau_0 > 0.352 -0.05925$, $y_4=2\sqrt2$,
$y_5,y_6\in[2t_0,2\sqrt2]$. \refno{535906363\dag}

\subsection{Group}
\label{app:p8} 
$\vor_0 < -0.146$, if $y_5,y_6\in[2\sqrt2,3.2]$. \refno{516537931}

$\tau_0(S(y_1,\ldots,y_6))+
    \tau_0(S(y_1,y_2',y_3,2,y_5,2)) > 0.31$, if $y_5,y_6\in[2\sqrt2,3.2]$,
    and $y_2'\in\{2,2t_0\}$.
\refno{130008809}

\subsection{Group}
\label{app:p11} 
$\vor_0<-0.084$, if $y_5\in[2t_0,2\sqrt2]$, $y_6\in[2\sqrt2,3.2]$.
\refno{531861442}

$\vor_0 < -0.084 - (y_5-2\sqrt2)0.1$, if $y_1=y_3=y_4=2$,
$y_6=2t_0$,
    $y_5\in[2\sqrt2,3.2]$.
\refno{292827481}

$\vor_0 < 0.009 + (y_5-2\sqrt2)0.1$, if $y_5\in[2\sqrt2,3.2]$,
    $y_1=y_3=y_4=y_6=2$.
\refno{710875528}

$\tau_0>0.176$, if $y_5\in[2t_0,2\sqrt2]$, $y_6\in[2\sqrt2,3.2]$.
\refno{286122364}

\newpage\twocolumn

    \renewcommand{\thefootnote}{}
    \footnote{This paper was rewritten on 3/3/02 from the version of 7/31/98.}
    \footnote{Research supported in part by the NSF}

\end{document}